\documentclass[11pt]{amsart}
\newcommand\mtop{1in}
\newcommand\mbottom{1in}
\newcommand\mleft{1.2in}
\newcommand\mright{1.2in}
\usepackage{amssymb}
\usepackage{amsmath}
\usepackage{amsthm}
\usepackage{stmaryrd}
\usepackage{mathrsfs}
\usepackage{color}
\usepackage[colorlinks=true,linkcolor=mblue,citecolor=mblue]{hyperref}
\usepackage{dsfont}
\usepackage[all]{xy}
    \SelectTips{cm}{10}  % Better arrowheads; shamelessly stolen from Anton Geraschenko's notes
\usepackage{float}
\usepackage{calc}
\usepackage{tikz}
\usepackage{mathdots}
\usepackage{colonequals}
\usepackage{adjustbox}
\usepackage{caption}
\usepackage{longtable}
\providecommand{\mparwidth}{1in}
\providecommand{\mtop}{1in}
\providecommand{\mbottom}{1in}
\providecommand{\mleft}{1.2in}
\providecommand{\mright}{1.2in}
\usepackage[top = \mtop, bottom = \mbottom, left = \mleft, right=\mright, marginparwidth=\mparwidth]{geometry}
\usepackage{fancyhdr}
\usepackage{setspace}
\usepackage{mathtools}
\usepackage{scalefnt}
\usepackage{microtype}
\usepackage{pifont}
\usepackage{ifthen}
\usepackage{expl3}
\ExplSyntaxOn
\let\exp_last_unbraced:NNNf\relax
\ExplSyntaxOff
\usepackage[clockwise]{rotating}
\usetikzlibrary{arrows.meta}

\usepackage{enumitem}
\setlist[enumerate,1]{leftmargin=6ex,topsep=-5pt,label=(\arabic*)}
\setlist[enumerate]{itemsep=7pt}
\setlist[itemize,1]{leftmargin=4ex,topsep=-1em}
\setlist[itemize]{itemsep=5pt}
\setlist[itemize,2]{label=$\circ$}
\setlist[itemize,3]{label={\scalefont{0.6}\color{gray}$\blacktriangleright$}}
\setlist[itemize,4]{label=$\ast$}
\setlist{nolistsep}

\DeclareMathOperator*{\colim}{colim}

\usepackage{ifthen}
\usepackage{pifont}
\usepackage[new]{old-arrows}
\usepackage[nolabel]{showlabels}

\begin{document}

\newcommand{\theoremnumstyle}{section}
%% Actual formatting
\parskip=0.2in \parindent=0in
\allowdisplaybreaks % Allow align environments to split over page breaks
\raggedbottom % Don't leave awkward spaces in the middle of pages

% From Atanas Atanasov: \replacecommand is similar to \providecommand but it initializes the command as necessary even if a previous definition exists.
\newcommand{\replacecommand}[2]{\providecommand{#1}{}\renewcommand{#1}{#2}}

\renewcommand{\arraystretch}{1.5}

%% Miscellaneous \newcommand's
\renewcommand{\l}{\overset}
\newcommand{\into}{\hookrightarrow}
\newcommand{\onto}{\twoheadrightarrow}
\newcommand{\tto}{\longrightarrow}
\newcommand{\too}[1]{\l{#1}\to}
\newcommand{\ttoo}[1]{\l{#1}\longrightarrow}
\newcommand{\intoo}[1]{\l{#1}\into}
\newcommand{\ontoo}[1]{\l{#1}\onto}
\newcommand{\mapstoo}[1]{\l{#1}\mapsto}
\newcommand{\bto}{\leftarrow}
\newcommand{\btto}{\longleftarrow}
\newcommand{\btoo}[1]{\l{#1}\bto}
\newcommand{\bttoo}[1]{\l{#1}\longleftarrow}
\newcommand{\binto}{\hookleftarrow}
\newcommand{\bonto}{\twoheadleftarrow}
\newcommand{\bintoo}[1]{\l{#1}\binto}
\newcommand{\bontoo}[1]{\l{#1}\bonto}
% stupid hack to get the open immersion symbol; see http://tex.stackexchange.com/questions/66723/nudging-overset-characters-downwards
\newcommand{\ointo}{\hspace{3pt}\text{\raisebox{-1.5pt}{$\overset{\circ}{\vphantom{}\smash{\text{\raisebox{1.5pt}{$\into$}}}}$}}\hspace{3pt}}
\newcommand{\lu}{\underset}
\newcommand{\bimplies}{\impliedby}
\newcommand{\ints}{\cap}
\newcommand{\intss}{\bigcap}
\newcommand{\union}{\cup}
\newcommand{\unions}{\bigcup}
\newcommand{\djunion}{\sqcup}
\newcommand{\djunions}{\bigsqcup}
\newcommand{\propersubset}{\subsetneq}
\newcommand{\propersupset}{\supsetneq}
\newcommand{\contains}{\supset}
\newcommand{\semidirect}{\rtimes}
\newcommand{\isom}{\cong}
\newcommand{\normal}{\triangleleft}
\replacecommand{\dsum}{\oplus}
\newcommand{\dsums}{\bigoplus}
\newcommand{\tensor}{\otimes}
\newcommand{\tensors}{\bigotimes}
\newcommand{\cotensor}{{\,\scriptstyle\square}}
\let\originalbar\bar
\renewcommand{\bar}[1]{{\overline{#1}}}
\newcommand{\rlim}{\mathop{\varinjlim}\limits}
\newcommand{\llim}{\mathop{\varprojlim}\limits}
\newcommand{\x}{\times}
\providecommand{\st}{\hspace{2pt} : \hspace{2pt}}
\newcommand{\vv}{\vspace{10pt}}
\newcommand{\til}{\widetilde}
\renewcommand{\hat}{\widehat}
\newcommand{\hhat}{{\hat{\ }}}
\newcommand{\iy}{\infty}
\newcommand{\hteq}{\simeq}
\newcommand{\dd}[2]{\frac{\partial #1}{\partial #2}}
\newcommand{\sm}{\wedge} % stop getting confused about the word "wedge"
\newcommand{\noqed}{\renewcommand{\qedsymbol}{}}
\newcommand{\adjoint}{\dashv}
\newcommand{\wreath}{\wr}
\newcommand{\heart}{\heartsuit}

%http://tex.stackexchange.com/questions/23432/how-to-create-my-own-math-operator-with-limits
\newcommand{\bigast}{\mathop{\vphantom{\sum}\mathchoice%
  {\vcenter{\hbox{\huge *}}}
  {\vcenter{\hbox{\Large *}}}{*}{*}}\displaylimits}

%% \operatorname
\renewcommand{\dim}{\operatorname{dim}}
\newcommand{\diam}{\operatorname{diam}}
\newcommand{\coker}{\operatorname{coker}}
\newcommand{\im}{\operatorname{im}}
\newcommand{\disc}{\operatorname{disc}}
\newcommand{\Pic}{\operatorname{Pic}}
\newcommand{\Der}{\operatorname{Der}}
\newcommand{\ord}{\operatorname{ord}}
\newcommand{\nil}{\operatorname{nil}}
\newcommand{\rad}{\operatorname{rad}}
\newcommand{\ssum}{\operatorname{sum}}
\newcommand{\codim}{\operatorname{codim}}
\newcommand{\cchar}{\operatorname{char}}
\newcommand{\sspan}{\operatorname{span}}
\newcommand{\rank}{\operatorname{rank}}
\newcommand{\Aut}{\operatorname{Aut}}
\newcommand{\Div}{\operatorname{Div}}
\newcommand{\Gal}{\operatorname{Gal}}
\newcommand{\Hom}{\operatorname{Hom}}
\newcommand{\Mor}{\operatorname{Mor}}
\newcommand{\Vect}{\operatorname{Vect}}
\newcommand{\Fun}{\operatorname{Fun}}
\newcommand{\Iso}{\operatorname{Iso}}
\newcommand{\Map}{\operatorname{Map}}
\newcommand{\Ho}{\operatorname{Ho}}
\newcommand{\Mod}{\operatorname{Mod}}
\newcommand{\Tot}{\operatorname{Tot}}
\newcommand{\cofib}{\operatorname{cofib}}
\newcommand{\fib}{\operatorname{fib}}
\newcommand{\hocofib}{\operatorname{hocofib}}
\newcommand{\hofib}{\operatorname{hofib}}
\newcommand{\Maps}{\operatorname{Maps}}
\newcommand{\Sym}{\operatorname{Sym}}
\newcommand{\Diff}{\operatorname{Diff}}
\newcommand{\Tr}{\operatorname{Tr}}
\newcommand{\Frac}{\operatorname{Frac}}
\renewcommand{\Re}{\operatorname{Re}} % don't like the curly default ones
\renewcommand{\Im}{\operatorname{Im}}
\newcommand{\gr}{\operatorname{gr}}
\newcommand{\tr}{\operatorname{tr}}
\newcommand{\End}{\operatorname{End}}
\newcommand{\Mat}{\operatorname{Mat}}
\newcommand{\Proj}{\operatorname{Proj}}
\newcommand{\Th}{\operatorname{Thom}}
\newcommand{\Thom}{\operatorname{Thom}}
\newcommand{\Spec}{\operatorname{Spec}}
\newcommand{\Ext}{\operatorname{Ext}}
\newcommand{\Cotor}{\operatorname{Cotor}}
\newcommand{\Tor}{\operatorname{Tor}}
\newcommand{\vol}{\operatorname{vol}}

%% new operatornames since preamble8
\newcommand{\Set}{\operatorname{Set}}
\newcommand{\Top}{\operatorname{Top}}
\newcommand{\Fin}{\operatorname{Fin}}
\newcommand{\Spaces}{\operatorname{Spaces}}
\newcommand{\Sp}{\operatorname{Sp}}
\newcommand{\Spectra}{\operatorname{Spectra}}
\newcommand{\Spt}{\operatorname{Spt}}
\newcommand{\Comod}{\operatorname{Comod}}
\newcommand{\Spf}{\operatorname{Spf}}
\newcommand{\tmf}{\mathit{tmf}}
\newcommand{\Tmf}{\mathit{Tmf}}
\newcommand{\TMF}{\mathit{TMF}}
\newcommand{\Null}{\operatorname{Null}}
\newcommand{\Fil}{\operatorname{Fil}}
\newcommand{\Sq}{\operatorname{Sq}}
\newcommand{\Stable}{\operatorname{Stable}}
\newcommand{\Poly}{\operatorname{Poly}}
\newcommand{\Cat}{\operatorname{Cat}}
\newcommand{\Orb}{\operatorname{Orb}}
\newcommand{\Exc}{\operatorname{Exc}}
\newcommand{\Part}{\operatorname{Part}}
\newcommand{\Comm}{\operatorname{Comm}}
\newcommand{\Res}{\operatorname{Res}}
\newcommand{\Thick}{\operatorname{Thick}}

\newcommand{\Sm}{\operatorname{Sm}}
\newcommand{\Var}{\operatorname{Var}}
\newcommand{\Frob}{\operatorname{Frob}}
\newcommand{\Rep}{\operatorname{Rep}}
\newcommand{\Ch}{\operatorname{Ch}}
\newcommand{\Shv}{\operatorname{Shv}}
\newcommand{\Corr}{\operatorname{Corr}}
\newcommand{\Span}{\operatorname{Span}}
\newcommand{\Sch}{\operatorname{Sch}}
\newcommand{\ev}{\operatorname{ev}}
\newcommand{\Homog}{\operatorname{Homog}}
\newcommand{\conn}{\operatorname{conn}}
\newcommand{\type}{\operatorname{type}}
\newcommand{\num}{\operatorname{num}}
\newcommand{\Aff}{\operatorname{Aff}}
\newcommand{\Psh}{\operatorname{Psh}}
\newcommand{\sk}{\operatorname{sk}}
\newcommand{\Cart}{\operatorname{Cart}}

\newcommand{\Br}{\operatorname{Br}}
\newcommand{\BW}{\operatorname{BW}}
\newcommand{\Cl}{\operatorname{Cl}}
\newcommand{\Conf}{\operatorname{Conf}}
\newcommand{\Alg}{\operatorname{Alg}}
\newcommand{\CAlg}{\operatorname{CAlg}}
\newcommand{\Lie}{\operatorname{Lie}}
\newcommand{\Coalg}{\operatorname{Coalg}}
\newcommand{\Ab}{\operatorname{Ab}}
\newcommand{\Ind}{\operatorname{Ind}}
\newcommand{\ind}{\operatorname{ind}}
\newcommand{\Fix}{\operatorname{Fix}}
\newcommand{\ho}{\operatorname{ho}}
\newcommand{\coeq}{\operatorname{coeq}}
\newcommand{\CMon}{\operatorname{CMon}}
\newcommand{\Sing}{\operatorname{Sing}}
\newcommand{\Inj}{\operatorname{Inj}}
\newcommand{\StMod}{\operatorname{StMod}}
\newcommand{\Loc}{\operatorname{Loc}}
\newcommand{\Free}{\operatorname{Free}}
\newcommand{\Art}{\operatorname{Art}}
\newcommand{\Gpd}{\operatorname{Gpd}}
\newcommand{\Def}{\operatorname{Def}}
\newcommand{\Hyp}{\operatorname{Hyp}}
\newcommand{\Pre}{\operatorname{Pre}}
\newcommand{\Lat}{\operatorname{Lat}}
\newcommand{\Coords}{\operatorname{Coords}}
\newcommand{\cone}{\operatorname{cone}}
\newcommand{\Spc}{\operatorname{Spc}}
\newcommand{\QCoh}{\operatorname{QCoh}}
\newcommand{\height}{\operatorname{ht}}
\newcommand{\Sub}{\operatorname{Sub}}
\newcommand{\Cone}{\operatorname{Cone}}
\newcommand{\Cocone}{\operatorname{Cocone}}
\newcommand{\Ran}{\operatorname{Ran}}
\newcommand{\Lan}{\operatorname{Lan}}
\newcommand{\LieAlg}{\operatorname{LieAlg}}
\newcommand{\Com}{\operatorname{Com}}
\newcommand{\CoAlg}{\operatorname{CoAlg}}
\newcommand{\Prim}{\operatorname{Prim}}
\newcommand{\Coh}{\operatorname{Coh}}
\newcommand{\FormalGrp}{\operatorname{FormalGrp}}
\newcommand{\Fact}{\operatorname{Fact}}
%END

%% \A etc.
\newcommand{\A}{\mathbb{A}}
\replacecommand{\C}{\mathbb{C}}
\newcommand{\CP}{\mathbb{C}\mathrm{P}}
\newcommand{\E}{\mathbb{E}}
\newcommand{\F}{\mathbb{F}}
\replacecommand{\G}{\mathbb{G}}
\renewcommand{\H}{\mathbb{H}}
\newcommand{\K}{\mathbb{K}}
\newcommand{\M}{\mathbb{M}}
\newcommand{\N}{\mathbb{N}}
\renewcommand{\O}{\mathcal{O}}
\renewcommand{\P}{\mathbb{P}}
\newcommand{\Q}{\mathbb{Q}}
\newcommand{\R}{\mathbb{R}}
\newcommand{\RP}{\mathbb{R}\mathrm{P}}
\newcommand{\V}{\vee}
\newcommand{\T}{\mathbb{T}}
\providecommand{\U}{\mathscr{U}}
\newcommand{\Z}{\mathbb{Z}}
\renewcommand{\k}{\Bbbk}
\newcommand{\g}{\mathfrak{g}}
\newcommand{\m}{\mathfrak{m}}
\newcommand{\n}{\mathfrak{n}}
\newcommand{\p}{\mathfrak{p}}
\newcommand{\q}{\mathfrak{q}}
\renewcommand{\t}{\mathfrak{t}}

%% Large parentheses, etc.
\newcommand{\pa}[1]{\left( {#1} \right)}
\newcommand{\br}[1]{\left[ {#1} \right]}
\newcommand{\cu}[1]{\left\{ {#1} \right\}}
\newcommand{\ab}[1]{\left| {#1} \right|}
\newcommand{\an}[1]{\left\langle {#1}\right\rangle}
\newcommand{\fl}[1]{\left\lfloor {#1}\right\rfloor}
\newcommand{\ceil}[1]{\left\lceil {#1}\right\rceil}
\newcommand{\tf}[1]{{\textstyle{#1}}}
\newcommand{\patf}[1]{\pa{\textstyle{#1}}}

%% Weird constructions
\renewcommand{\mp}{\ \raisebox{5pt}{\text{\rotatebox{180}{$\pm$}}}\ }
\renewcommand{\d}[1]{\ss \mathrm{d}#1}
\newcommand{\imod}{\hspace{-7pt}\pmod}
%\renewcommand{\check}[1]{\overset{\smile}{#1}}

%% Better versions of existing commands
\renewcommand{\epsilon}{\varepsilon}
\renewcommand{\phi}{{\mathchoice{\raisebox{2pt}{\ensuremath\varphi}}{\raisebox{2pt}{\!\! \ensuremath\varphi}}{\raisebox{1pt}{\scriptsize$\varphi$}}{\varphi}}}
\newcommand{\ph}{{\color{white}.\!}}
\newcommand{\tspacer}{{\ensuremath{\color{white}\Big|\!}}}
\newcommand{\chii}{\raisebox{2pt}{\ensuremath\chi}}

% from http://mbork.pl/2009-04-27_Fun_with_quantifiers_%28en%29
\let\originalchi=\chi
\renewcommand{\chi}{{\!{\mathchoice{\raisebox{2pt}{
$\originalchi$}}{\!\raisebox{2pt}{
$\originalchi$}}{\raisebox{1pt}{\scriptsize$\originalchi$}}{\originalchi}}}}

\let\originalforall=\forall
\renewcommand{\forall}{\ \originalforall}

\let\originalexists=\exists
\renewcommand{\exists}{\ \originalexists}

\let\realcheck\check
\newcommand{\vH}{\realcheck{H}}

% quote block: \begin{qu}{leftmargin}{rightmargin}{  ...  }
\newenvironment{qu}[2]
{\begin{list}{}
	  {\setlength\leftmargin{#1}
	  \setlength\rightmargin{#2}}
	  \item[]\footnotesize}
		  {\end{list}}

%%%%%%%%%%%%%%%%%%%%%%%%%
%% newcommands that require certain packages

\newcommand{\itext}{\shortintertext} % requires mathtools
\renewcommand{\u}{\underbracket[0.7pt]} % requires mathtools
\newcommand{\margin}[1]{\marginpar{\raggedright \scalefont{0.7}#1}} % requires scalefnt

% require package pigpen and xy
\newcommand{\pullback}{\ar@{}[rd]|<<{\text{\pigpenfont A}}}
\newcommand{\pushout}{\ar@{}[rd]|>>{\text{\pigpenfont I}}}

% requires xy
\newcommand{\longleftrightarrows}{\xymatrix@1@C=16pt{
\ar@<0.4ex>[r] & \ar@<0.4ex>[l]
}}
\newcommand{\longrightrightarrows}{\xymatrix@1@C=16pt{
\ar@<0.4ex>[r]\ar@<-0.4ex>[r] & 
}}
\newcommand{\mapstto}{\,\xymatrix@1@C=16pt{
\ar@{|->}[r] & 
}\,}
\newcommand{\mapsttoo}[1]{\xymatrix@1@C=16pt{
\ar@{|->}[r]^-{#1} & 
}}
\newcommand{\rightrightrightarrows}{\xymatrix@1@C=16pt{
\ar[r]\ar@<0.8ex>[r]\ar@<-0.8ex>[r] & 
}}
\newcommand{\longleftleftarrows}{\xymatrix@1@C=16pt{
 & \ar@<0.4ex>[l]\ar@<-0.4ex>[l]
}}
\newcommand{\leftleftleftarrows}{\xymatrix@1@C=16pt{
 & \ar[l]\ar@<0.8ex>[l]\ar@<-0.8ex>[l]
}}
\newcommand{\leftleftleftleftarrows}{\xymatrix@1@C=16pt{
 & \ar@<0.8ex>[l]\ar@<0.3ex>[l]\ar@<-0.3ex>[l]\ar@<-0.8ex>[l]
}}
\newcommand{\lcircle}{\ar@(ul,dl)} % arrow circle to the left of the node
\newcommand{\rcircle}{\ar@(ur,dr)}
\newcommand{\intto}{\ \xymatrix@1@C=16pt{
\ar@{^(->}[r] & 
}}

% Colors
% Now this preamble doesn't cause errors when there is no \usepackage{color}
% (you just can't use these commands)
\makeatletter
\@ifundefined{mathds}{
	\newcommand{\Id}{Id}
	}{
	\newcommand{\Id}{\mathds{1}} % requires mathds
	}
\@ifundefined{color}{}{
	\definecolor{darkgreen}{RGB}{0,70,0}
	\definecolor{dgreen}{RGB}{0,100,0}
	\definecolor{purple}{RGB}{120,00,120}
	\definecolor{gray}{RGB}{100,100,100}
	\definecolor{mgreen}{RGB}{0,150,0}
	\definecolor{llgray}{RGB}{230,230,230}
	\definecolor{lgreen}{RGB}{100,200,100}
	\definecolor{mgray}{RGB}{150,150,150}
	\definecolor{lgray}{RGB}{190,190,190}
	\definecolor{maroon}{RGB}{150,0,0}
	\definecolor{lblue}{RGB}{120,170,200}
	\definecolor{mblue}{RGB}{65,105,225}
	\definecolor{dblue}{RGB}{0,56,111}
	\definecolor{orange}{RGB}{255,165,0}
	\definecolor{brown}{RGB}{177,84,15}
	\definecolor{rose}{RGB}{135,0,52}
	\definecolor{gold}{RGB}{177,146,87}
	\definecolor{dred}{RGB}{135,19,19}
	\definecolor{mred}{RGB}{194,28,28}
	\newcommand{\edit}[1]{\itshape{\color{gray}#1}\upshape}
	\newcommand{\fixme}[1]{{\color{maroon}\it{#1}}}
	\newcommand{\citeme}[1]{{\color{orange}\textit{#1}}}
	\newcommand{\later}[1]{{\color{rose}#1}}
	\newcommand{\corr}[1]{{\color{red}\itshape #1}}
	\newcommand{\question}[1]{\itshape{\color{blue}#1}\upshape}
}
\@ifundefined{substack}{}{
    \newcommand{\attop}[1]{{\let\textstyle\scriptstyle\let\scriptstyle\scriptscriptstyle\substack{#1}}}
}
\makeatother

%  vim:ft=tex

\newtheoremstyle{gloss}{\topsep}{\topsep}{}{0pt}{\bfseries}{}{\newline}{\newline
*{\bf #3} }
\theoremstyle{gloss}
\newtheorem*{defstar}{Definition}

\newtheoremstyle{newplain}{20pt}{0pt}{\it}{0pt}{\bfseries}{.}{1ex}{}
\theoremstyle{newplain}

% Number by section (1.1) by default, can override with
% \newcommand{\theoremnumstyle}{} (must be empty)
\ifthenelse{\isundefined\theoremnumstyle}
	{\newtheorem{theorem}{Theorem}[section] 
	\numberwithin{equation}{section}} % Number equations by section, like (1.1) instead of (1)
	{\ifthenelse{\equal\theoremnumstyle{}}
		{\newtheorem{theorem}{Theorem}}
		{\newtheorem{theorem}{Theorem}[section]
		\numberwithin{equation}{section}
		}
	}

\newtheorem{corollary}[theorem]{Corollary}
\newtheorem{claim}[theorem]{Claim}
\newtheorem{lemma}[theorem]{Lemma}
\newtheorem{proposition}[theorem]{Proposition}
\newtheorem{fact}[theorem]{Fact}

\newtheoremstyle{newtextthm}{20pt}{0pt}{}{0pt}{\bfseries}{.}{1ex}{}
\theoremstyle{newtextthm}
\newtheorem{definition}[theorem]{Definition}
\newtheorem{example}[theorem]{Example}
\newtheorem{problem}[theorem]{Problem}
\newtheorem{remark}[theorem]{Remark}
\newtheorem{notation}[theorem]{Notation}

\newtheorem*{theoremstar}{Theorem}
\newtheorem*{lemmastar}{Lemma}
\newtheorem*{corstar}{Corollary}
\newtheorem*{corollarystar}{Corollary}
\newtheorem*{propositionstar}{Proposition}
\newtheorem*{claimstar}{Claim}
\newtheorem*{examplestar}{Example}

% random
\newcommand{\argforrandom}{}
\theoremstyle{newtextthm}
\newtheorem{helperforrandom}[theorem]{\argforrandom}
\newtheorem*{helperforrandomstar}{\argforrandom}
\newenvironment{random}[1]{\renewcommand{\argforrandom}{#1}\begin{helperforrandom}}{\end{helperforrandom}}
\newenvironment{randomstar}[1]{\renewcommand{\argforrandom}{#1}\begin{helperforrandomstar}}{\end{helperforrandomstar}}

\newenvironment{exercise}[1]{\hspace{1pt}\nn \large {\sc #1.}\hv \normalsize
\vspace{10pt}\\ }{} 
\newcommand{\subthing}[1]{\hv\large(#1)\hv\hv \normalsize }

\renewcommand{\showlabelfont}{\tiny\color{red}}

% For some reason doing the obvious thing was slow; this compiles the
% problematic line just once.
\newbox\deltabox
\setbox\deltabox = \hbox{\scalefont{0.2}$\Delta$}
\newcommand{\tensorD}{\l{\copy\deltabox}\tensor}
\newcommand{\cotensorD}{\l{\copy\deltabox}\cotensor}
\newlength\len
\newcommand{\ED}[2]{\setbox3=\hbox{$E^{#1}_{#2}$}\setbox4=\hbox{$E$}\setlength\len{-0.5\wd3+0.5\wd4}\l{\raisebox{-2pt}
{\hspace{1.25\len}\copy\deltabox}}{E^{#1}_{#2}}}

\newbox\lbox
\setbox\lbox = \hbox{\scalefont{0.2}$L$}
\newbox\rbox
\setbox\rbox = \hbox{\scalefont{0.2}$R$}
\newcommand{\tensorL}{\l{\copy\lbox}\tensor}
\newcommand{\tensorR}{\l{\copy\rbox}\tensor}
\newcommand{\DL}[1]{\l{\!\!{\copy\lbox}}{D^{#1}_\Gamma}}
\newcommand{\DR}[1]{\l{\!\!{\copy\rbox}}{D^{#1}_\Gamma}}
\newcommand{\CL}[1]{\l{\!\!{\copy\lbox}}{C^{#1}_\Gamma}}
\newcommand{\CLR}[1]{\l{\!\!{\copy\rbox}}{C^{#1}_\Gamma}}
\newcommand{\CLPhi}[1]{\l{\!\!{\copy\lbox}}{C^{#1}_\Phi}}
\newcommand{\DLPhi}[1]{\l{\!\!{\copy\lbox}}{D^{#1}_\Phi}}
\newcommand{\barDL}[1]{\l{\!\!{\copy\lbox}}{\originalbar{D}^{#1}_\Gamma}}
\newcommand{\DD}[1]{\l{\!\!{\copy\deltabox}}{D^{#1}_\Gamma}}
\newcommand{\CD}[1]{\l{\!\!{\copy\deltabox}}{C^{#1}_\Gamma}}
\newcommand{\DDPhi}[1]{\l{\!\!{\copy\deltabox}}{D^{#1}_\Phi}}
\renewcommand{\cotensor}{\,\text{\scalefont{0.7}$\square$}}
\newcommand{\Imu}{I_{\Phi,\mu}}

\newcommand{\ann}[1]{\langle {#1}\rangle}
\newcommand{\si}{\l{\textit{st}}\isom}
\newcommand{\cc}{\ ;\ }

\newcommand{\EU}{{}^{U\hspace{-3pt}}E}

\renewcommand{\showlabelfont}{\tiny\slshape\color{mgreen}}

\author{Eva Belmont}
\title{Localizing the $E_2$ page of the Adams spectral sequence}
\lhead{}\chead{\it\small Localizing the $E_2$
page of the Adams spectral sequence}\rhead{}\lfoot{}\cfoot{\thepage}\rfoot{}

\maketitle

\begin{abstract}
There is only one nontrivial localization of $\pi_*S_{(p)}$ (the
chromatic localization at $v_0=p$), but there are infinitely many nontrivial
localizations of the Adams $E_2$ page for the sphere. The first non-nilpotent
element in the $E_2$ page after $v_0$ is $b_{10}\in
\Ext_A^{2p(p-1)-2}(\F_p,\F_p)$. We work at $p=3$ and study
$b_{10}^{-1}\Ext_P(\F_3,\F_3)$ (where $P$ is the algebra of dual reduced
powers), which agrees with the infinite summand $\Ext_P(\F_3,\F_3)$ of
$\Ext_A(\F_3,\F_3)$ above a line of slope ${1\over 23}$. We compute up to the
$E_9$ page of an Adams spectral sequence in the category $\mathrm{Stable}(P)$
converging to $b_{10}^{-1}\Ext_P(\F_3,\F_3)$, and conjecture that the spectral
sequence collapses at $E_9$. We also give a complete calculation of
$b_{10}^{-1}\Ext_P^*(\F_3,\F_3[\xi_1^3])$.
\end{abstract}
\setcounter{tocdepth}{1}
{\parskip=0in \tableofcontents}

\section{Introduction}
For a $p$-local finite spectrum $X$, the Adams spectral sequence
$$ E_2^{**} = \Ext_A^*(\F_p,H_*X)\implies \pi_*X^\hhat_p $$
is one of the main tools for computing (the $p$-completion of) the homotopy
groups of $X$. If one understands the $A$-comodule structure of $H_*X$, it is
possible to compute the $E_2$ page algorithmically in a finite range of
dimensions. However, for many spectra $X$ of interest such as the sphere
spectrum, there is no chance of determining the $E_2$ page completely.
The motivating goal behind this work is to compute an infinite part of
the Adams $E_2$ page $\Ext_A^*(\F_3,\F_3)$ for the sphere at $p=3$.
Specifically, we wish to compute the $b_{10}$-periodic part, where $b_{10}\in
\Ext^{2,2p(p-1)}(\F_p,\F_p)$ converges to $\beta_1\in \pi_{2p(p-1)-2}S$. We show
that there is a plane above which $\Ext_A^*(\F_p,\F_p)$ is $b_{10}$-periodic,
where the third grading $f$ (in addition to internal degree $t$ and homological
degree $s$) is related to the collapse of the Cartan-Eilenberg spectral sequence
at odd primes $p$ (see \eqref{CESS-collapse}).

%In particular, if $M$ is an evenly graded $A$-comodule, $P =
%\F_p[\xi_1,\xi_2,\dots]$ is the algebra of dual reduced powers, and $E =
%\F_p[\tau_0,\tau_1,\dots]/(\tau_i^2)$, then the Cartan-Eilenberg spectral
%sequence associated to the extension

The only known localization of the Adams $E_2$ page for the sphere is
\begin{equation}\label{a_0-inverted-Ext-sphere} a_0^{-1}\Ext_A^*(\F_p,\F_p)
\isom \F_p[a_0^{\pm 1}]
\end{equation}
where $a_0 = [\tau_0]$ converges to $p\in \pi_0^*S$; this follows from Adams'
fundamental work \cite{adams-periodicity} on the structure of the $E_2$ page.
This localization agrees with $\Ext_A^*(\F_p,\F_p)$ above a line of slope
${1\over 2p-2}$ (in the $(t-s,s)$ grading). Our proposed localization
$b_{10}^{-1}\Ext_A^*(\F_p,\F_p)$ agrees with $\Ext_A^*(\F_p,\F_p)$ above a plane
whose fixed-$f$ cross section is a line of slope ${1\over p^3-p-1}$. While the
only $a_0$-periodic elements lie in the zero-stem (corresponding to chromatic
height zero), the $b_{10}$-periodic region encompasses nonzero classes in
arbitrarily high stems, including some elements in chromatic height 2, such as
$b_{10}$ itself. Though we do not give a complete calculation of
$b_{10}^{-1}\Ext_A^*(\F_p,\F_p)$, we will see that it is much more complicated
than $a_0^{-1}\Ext_A^*(\F_p,\F_p)$. Thus in some sense, one may think of
$b_{10}^{-1}\Ext_A^*(\F_p,\F_p)$ as a richer and more revealing version of the
classical calculation.

In a different sense, however, these two localizations come from different
worlds. Inverting $a_0$ is the Adams $E_2$ avatar of $p$-localization on
($p$-local) homotopy (rationalization). Equivalently, the sphere has chromatic
type zero, and $a_0$ is just the algebraic name for the chromatic height-0
operator $v_0$. On the other hand, inverting $b_{10}$ is not the shadow of any
homotopy-theoretic localization: by the Nishida nilpotence theorem, $\beta_1$ is
nilpotent in homotopy, so $\beta_1^{-1}\pi_*^sS = 0$. While $v_0=p$ is the only
chromatic periodicity operator acting on the sphere, $a_0$ and $b_{10}$ are just
the first two out of infinitely many non-nilpotent elements in
$\Ext_A^*(\F_p,\F_p)$. Palmieri \cite{palmieri-book} describes a more
complicated analogue of the classical theory of periodicity and nilpotence that
operates only on Adams $E_2$ pages, almost all of which (except the $v_n$
operators) is destroyed by the time one reaches the Adams $E_\iy$ page.

Recall that the odd-primary dual Steenrod algebra has a presentation $A =
\F_p[\xi_1,\xi_2,\dots]\tensor E[\tau_0,\tau_1,\dots]$ where $E[x] =
\F_p[x]/x^2$ denotes an exterior algebra, $|\xi_n| = 2(p^n-1)$ and $|\tau_n| =
2p^n-1$. Let $P = \F_p[\xi_1,\xi_2,\dots]$ be the Steenrod reduced powers
algebra, and let $E$ be the quotient Hopf algebra $E[\tau_0,\tau_1,\dots]$. If
$M$ is an evenly graded $A$-comodule, there is an isomorphism
\begin{equation}\label{CESS-collapse} \Ext_A^{s,t}(\F_p,M)\isom
\Ext_P^{s,t-f}(\F_p,\Ext_E^{f,*}(\F_p,M))
\end{equation}
which arises from the collapse of the Cartan-Eilenberg spectral sequence at odd
primes $p$. In light of this, we recast our goal as follows:
\begin{random}{Goal}
Compute $b_{10}^{-1}\Ext^*_P(\F_3,M)$ for $P$-comodules $M$.
\end{random}
In particular, we are most interested in $M = \Ext_E^*(\F_p,\F_p)$. In this
paper, we focus on the $f=0$ summand $\Ext_E^0(\F_p,\F_p) \isom \F_p$. We show:
\begin{theorem}\label{main-result}
Let $D = \F_p[\xi_1]/\xi_1^3$ and let $R = b_{10}^{-1} \Ext_D(\F_p,\F_p) =
E[h_{10}] \tensor \F_p[b_{10}^{\pm 1}]$. There is a spectral sequence
$$ E_2 \isom R[w_2,w_3,\dots] \implies b_{10}^{-1}\Ext_P(\F_p,\F_p) $$
where $w_n$ has filtration 1, internal homological degree 1, and internal
topological degree $2(3^n+1)$. For degree reasons, $d_r = 0$ for $r\geq 2$
unless $r\equiv 4\pmod 9$ or $r\equiv 8\pmod 9$. The first nontrivial
differential is
$$ d_4(w_n) = b_{10}^{-4}h_{10}w_2^2w_{n-1}^3. $$
Furthermore, we give a complete description of the $d_8$'s.
\end{theorem}
We conjecture that the spectral sequence collapses at $E_9$, and show that this
is equivalent to the following conjecture.

\begin{random}{Conjecture}\label{tilW-conj}
There is an isomorphism
$$b_{10}^{-1}\Ext_P^*(\F_3,\F_3) \isom b_{10}^{-1}
\Ext_D^*(\F_3,\til{W})$$
where $\til{W} = \F_p[\til{w}_2,\til{w}_3,\dots]$ with $|\til{w}_n| = 2(3^n-5)$
and coaction given by $\psi(\til{w}_n) = 1\tensor \til{w}_n + \xi_1\tensor
\til{w}_2^2 \til{w}_{n-1}^3$ for $n\geq 3$. (These generators are related to the
generators of Theorem \ref{main-result} by $\til{w}_n = b_{10}^{-1}w_n$.)
\end{random}
Adams' theorem \eqref{a_0-inverted-Ext-sphere} has the more general form
$$ a_0^{-1}\Ext_A^*(\F_p,M) \isom a_0^{-1}\Ext_{E[\tau_0]}(\F_p,M) $$
for an $A$-comodule $M$ (see also \cite{may-milgram}). In particular, the
localized cohomology depends only on the $E[\tau_0]$-comodule structure on $M$.
The analogous statement for $b_{10}$-localization (that $b_{10}^{-1}
\Ext_P^*(\F_3,M)$ depends only on the $D$-comodule structure of $M$) is not
true. In general, we propose the following:
\begin{random}{Conjecture}
There is a functor $\mathscr{E}: \Comod_P\to \Comod_D$ such that
$$ b_{10}^{-1}\Ext_P^*(\F_3,M) \isom
b_{10}^{-1}\Ext_D(\F_3,\mathscr{E}(M)) $$
and, as vector spaces, $R\tensor \mathscr{E}(M)$ agrees with the
$E_2$ page of the spectral sequence described below in Theorem
\ref{K(xi_1)-MPASS-M} with $\Gamma=P$.
\end{random}

Our best complete result is the following; it is proved in Section \ref{section:D_{1,iy}} using different methods.
\begin{theorem}\label{D_iy}
There is an isomorphism
$$b_{10}^{-1}\Ext_P(\F_3,\F_3[\xi_1^3]) \isom b_{10}^{-1}\Ext_D(\F_3,
\F_3[h_{20},b_{20},w_3,w_4,\dots]/h_{20}^2)$$
where $D$ acts trivially on all the generators on the right.
\end{theorem}

\subsection{Main tool}\label{section:main-tool}
Our main tool (the spectral sequence mentioned in Theorem \ref{main-result}) is
as follows. It is a special case of the construction discussed in
\cite{CESS-paper}.
\begin{theorem}\label{K(xi_1)-MPASS-M}
Let $D = \F_p[\xi_1]/\xi_1^p$ and let $\Gamma$ be a Hopf algebra over $\F_p$
with a surjection of Hopf algebras $\Gamma\to D$. Let $B_\Gamma =
\Gamma\cotensor_D \F_p$. For a $\Gamma$-comodule $M$, there is a spectral
sequence
\begin{align}
\label{E_1-general} E_1^{s,t} \isom b_{10}^{-1}\Ext_D^t(\F_p,\bar{B}_\Gamma^{\tensor
s}\tensor M) \implies b_{10}^{-1}\Ext_\Gamma(\F_p,M)
\end{align}
(where $\bar{B}_\Gamma$ is the coaugmentation ideal $\coker(\F_p\to B_\Gamma)$).
At $p=3$, $b_{10}^{-1} \Ext_D^*(\F_p,B_\Gamma)$ is flat as a $b_{10}^{-1}
\Ext_D^*(\F_p,\F_p)$-module, and
\begin{align}
\label{E_2-general} E_2^{**} & \isom b_{10}^{-1}
\Ext^*_{b_{10}^{-1}\Ext^*_D(\F_p,B_\Gamma)}(R,b_{10}^{-1}\Ext^*_D(\F_p,M)).
\end{align}
\end{theorem}
We work at $p=3$ throughout. The main focus is the case $\Gamma = P$ and $B_P =
P\cotensor_D \F_3\equalscolon B$; this is the spectral sequence of Theorem
\ref{main-result}. We also apply this for two quotients of $P$---for a spectral
sequence comparison argument in Section \ref{section:d_4} and for the proof of
Theorem \ref{D_iy} in Section \ref{section:D_{1,iy}}. Convergence is proved in
Appendix A in the case that $\Gamma$ is a quotient of $P$.

In \cite{CESS-paper}, we show that the following three constructions of
\eqref{E_1-general} coincide at the $E_1$ page.
\begin{enumerate} 
\item The first construction is a $b_{10}$-localized Cartan-Eilenberg-type
spectral sequence associated to the sequence of $P$-comodule algebras
$B\to P \to D$. (Note that the inclusion $B\to P$ is not a map of coalgebras; see
\cite[\S2.3]{CESS-paper} for a precise construction in this case.)
\item The second construction is an Adams spectral sequence internal to the
category $\Stable(P)$. See \cite[Chapter 14]{margolis} or \cite[\S9.6]{HPS} for
a definition of $\Stable(\Gamma)$ for a Hopf algebra $\Gamma$ over $\F_p$, or
\cite[\S4]{BHV} for a more modern viewpoint; the idea is that it is a variation
of the derived category of $\Gamma$-comodules designed to satisfy
$\Hom_{\Stable(\Gamma)}(\F_p,x^{-1} M) = x^{-1} \Hom_{\Stable(\Gamma)}(\F_p,M)$. In
particular, if $M$ is a $\Gamma$-comodule, then $\Hom_{\Stable(\Gamma)}(\F_p,M) =
\Ext_\Gamma^*(\F_p,M)$. 
The Adams spectral sequence in this setting was first studied by Margolis
\cite{margolis} and Palmieri \cite{palmieri-book}, and so we call this the
\emph{Margolis-Palmieri Adams spectral sequence (MPASS)}.

In particular, let $K(\xi_1)\colonequals b_{10}^{-1}B = \colim(B\too{b_{10}}B
\too{b_{10}} \dots )$ (where the colimit is taken in $\Stable(P)$); then our
spectral sequence is the $K(\xi_1)$-based Adams spectral sequence.
\item The third construction is obtained by $b_{10}$-localizing the filtration
spectral sequence on the normalized $P$-cobar complex $C_P^*(\F_p,\F_p)\colonequals
\bar{P}^{\tensor *}$ in which $[a_1|\dots|a_n]\in
F^sC_P^*$ if at least $s$ of the $a_i$'s lie in $\ker(P\to D) = \bar{B}P$.
\end{enumerate}
Our dominant viewpoint will be via the framework of (2), but the other two
formulations will be useful at key moments.
By a ``$b_{10}$-localized'' spectral sequence, we mean the spectral sequence
whose $E_r$ page is obtained by $b_{10}$-localizing the original $E_r$ page. It
is not automatic that this converges to the $b_{10}$-localization of the
original spectral sequence; this is what is checked in Appendix A.

The essential reason we focus on $p=3$ is that for the analogous construction at
$p>3$, the flatness condition does not hold. (This comes from the Adams spectral
sequence flatness condition applied in the setting of (2).)

\subsection{Outline}
In Section \ref{section:overview}, we prove some basic results about the
structure of the spectral sequence converging to $b_{10}^{-1}\Ext_P(k,k)$ and
introduce definitions and notation that will be used extensively in the
computational sections. In Section \ref{section:periodicity-line}, we apply
vanishing line results to describe a region in which $b_{10}^{-1}
\Ext_P(\F_p,M)$ agrees with $\Ext_P(\F_p,M)$. Sections
\ref{section:cooperations} and \ref{section:B-hopf-algebroid} are devoted to
computing the $E_2$ page of the $K(\xi_1)$-based MPASS converging to
$b_{10}^{-1}\Ext_P(\F_3,\F_3)$. In Section \ref{section:d_4} we determine $d_4$,
the first nontrivial differential after the $E_2$ page. In Section
\ref{section:d_8}, we determine $d_8$ and show that our conjecture that the
spectral sequence collapses at $E_9$ would imply the desired form of
$b_{10}^{-1}\Ext_P^*(\F_3,\F_3)$ in Conjecture \ref{tilW-conj}. In Section
\ref{section:D_{1,iy}} we prove Theorem \ref{D_iy}. In Appendix A we show
convergence of the MPASS in the cases of interest, and also show convergence of
an auxiliary spectral sequence needed for Section \ref{section:d_4}.

\subsection{Acknowledgements}
I am grateful to Haynes Miller, my graduate advisor, for suggesting this as a
thesis project and for providing invaluable guidance at every step along the
way. I would also like to thank Dan Isaksen and Zhouli Xu for
helpful conversations about this work, and Hood Chatham for productive
conversations and for the spectral sequences \LaTeX\ package used to draw the
charts in Appendix B.

\section{Overview of the MPASS converging to $b_{10}^{-1}\Ext_P(k,k)$} 
\label{section:overview}
In every section except Sections \ref{section:periodicity-line} and
\ref{section:cooperations} we will set $p=3$ and let $k = \F_3$.
We will denote exterior and truncated polynomial algebras,
respectively, by $E[x] = k[x]/x^2$ and $D[x] = k[x]/x^p$. Let $D = D[\xi_1]$.

If $M$ is a $P$-comodule and $E = b_{10}^{-1} M$, we adopt the notation of
\cite{palmieri-book} and write:
\begin{align*}
\pi_{**}(M) & \colonequals M_{**}\colonequals \Hom_{\Stable(P)}^{**}(k,M) =
\Ext^{**}_P(k,M)
\\M_{**}M & \colonequals \Hom_{\Stable(P)}^{**}(k,M\tensor M) = \Ext^{**}_P(k,M\tensor M)
\\\pi_{**}(E) & \colonequals E_{**} \colonequals \Hom_{\Stable(P)}^{**}(k,E)=
b_{10}^{-1} \Ext^{**}_P(k,M)
\\E_{**}E & \colonequals \Hom_{\Stable(P)}^{**}(k,E\tensor E) =
b_{10}^{-1}\Ext^{**}_P(k,M\tensor M).
\end{align*}
Here $M\tensor M$ is given the diagonal $P$-comodule structure: $\psi(a\tensor
b) = \sum a'b'\tensor a'' \tensor b''$ where $\psi(a) = \sum a'\tensor a''$ and
$\psi(b) = \sum b'\tensor b''$. Define
\begin{align*}
B & \colonequals P\cotensor_D k %k[\xi_1^p,\xi_2,\xi_3,\dots]
\\K(\xi_1) & \colonequals b_{10}^{-1}B \colonequals \colim(B\too{b_{10}}
B\too{b_{10}} \dots )
\end{align*}
where the colimit is taken in $\Stable(P)$.
Due to the general machinery of Adams spectral sequences in $\Stable(P)$ (see
\cite[\S1.4]{palmieri-book}), we have a $K(\xi_1)$-based spectral sequence
$$ E_1^{s,t,u} = \pi_{t,u}(K(\xi_1)\tensor \bar{K(\xi_1)}^{\tensor s}) =
b_{10}^{-1}\Ext_P(k,B\tensor \bar{B}^{\tensor s}) \implies b_{10}^{-1}\Ext_P(k,k) $$
which we call the \emph{$K(\xi_1)$-based Margolis-Palmieri Adams spectral
sequence (MPASS)}.
Here $\bar{(\ )}$ denotes coaugmentation ideal.
By the shear isomorphism (Lemma \ref{shear}) and the change of rings theorem, we
may write $E_1^{s,t,u} = b_{10}^{-1}\Ext_D^*(k, \bar{B}^{\tensor s})$.
If
$K(\xi_1)_{**}K(\xi_1)$ is flat over $K(\xi_1)_{**}$, then the $E_2$ page
\eqref{K(xi_1)-E_2-page} has
the form
\begin{equation}\label{K(xi_1)-E_2-page} \Ext_{K(\xi_1)_{**}K(\xi_1)}(K(\xi_1)_{**},K(\xi_1)_{**}).\end{equation}
The differential $d_r$ is a map $E_r^{s,t,u}\to E_r^{s+r,t-r+1,u}$. Here $s$ is
the MPASS filtration, $t$ is the internal homological degree, and $u$ is the
internal topological degree. Furthermore, we will often find it convenient to
work with the degree
$$ u':= u-6(s+t) $$
which has the property that $u'(b_{10})=0$. In this grading, the differential
$d_r$ is a map $E_r^{s,u'}\to E_r^{s+r,u'-6}$.

%The goal of the next two sections is to compute $K(\xi_1)_{**}K(\xi_1)$ and show
%that it is flat over $K(\xi_1)_{**}$. In fact, the flatness is a more general fact
%(see \citeme{somewhere}).
The coefficient ring $K(\xi_1)_{**}$ is easy to compute using the change of
rings theorem:
\begin{align*}
K(\xi_1)_{**} & = b_{10}^{-1}\Ext^*_P(k,B) =
b_{10}^{-1}\Ext^*_P(k,P\cotensor_D k)
\\ & = b_{10}^{-1}\Ext^*_D(k,k) = E[h_{10}]\tensor k[b_{10}^{\pm 1}]
\end{align*}
where $h_{10}$ is in homological degree 1 and $b_{10}$ is in homological degree
2. It will be useful to have notation for this coefficient ring:
\begin{equation} R\colonequals E[h_{10}]\tensor k[b_{10}^{\pm}]. \end{equation}
Using the shear isomorphism (Lemma \ref{shear}) and the change of rings theorem, we have 
\begin{equation}\label{cooperations-change-of-rings}K(\xi_1)_{**}K(\xi_1) \isom b_{10}^{-1}\Ext_P(k,B\tensor B) \isom
b_{10}^{-1}\Ext_P(k,P\cotensor_D B) \isom b_{10}^{-1}\Ext_D(k,B).\end{equation}

\begin{random}{Notation}\label{xi-antipode}
We have chosen to define $B$ as a left $P$-comodule. It can be written
explicitly as $\F_p[\bar{\xi}_1^p,\bar{\xi}_2,\bar{\xi}_3,\dots]$. To simplify
the notation, from now on we will redefine the symbol $\xi_n$ to mean the
antipode of the usual $\xi_n$. Thus, going forward, we will have $\Delta(\xi_n)
= \sum_{i+j=n} \xi_i\tensor \xi_j^{p^i}$, and
$$B = \F_p[\xi_1^p,\xi_2,\xi_3,\dots].$$
\end{random}

In Section \ref{section:B-hopf-algebroid}, we will show that the flatness
condition holds and
$K(\xi_1)_{**}K(\xi_1)$ is isomorphic, as a Hopf algebra over $R$, to an
exterior algebra on generators
$$ e_n = [\xi_1]\xi_n-[\xi_1^2]\xi_{n-1}^3 \in \Ext_D(k,B).$$
This implies that the $E_2$ page is isomorphic to a polynomial algebra
over $R$ on classes $w_n\colonequals [e_n]$ of degree $(s,t,u) =
(1,1,2(3^n+1))$.

The generator $w_2$ is a permanent cycle, and converges to
$g_0=\an{h_{10},h_{10},h_{11}}\in
\Ext_P^*(k,k)$. We will see in Section \ref{section:d_4} that the other $w_n$'s
support differentials, so it is less easy to see how these generators connect to
familiar elements in the Adams $E_2$ page. One heuristic comes from looking at
the images of these classes in $P/(\xi_1^3,\xi_2^9,\xi_3^9,\dots)$: in that
setting, $w_n = \an{h_{10},h_{10},h_{n-1,1}}$ and $h_{10}w_n = b_{10}h_{n-1,1}$.

Let $W_+ = k[b_{10}^{\pm 1}][w_2,w_3,\dots]$ and $W_- = W_+\{ h_{10} \}$.
Then $E_2 = W_+\dsum W_-$, and using simple degree arguments, we will show that
higher differentials take $W_+$ to $W_-$ and vice versa.

\begin{lemma}\label{s-possibilities}
Suppose $x\in E_2^{s(x),u'(x)}$ is nonzero. If $u'(x)\equiv 0\pmod 4$, then
$x\in W_+$ and $s \equiv -u' \pmod 9$.
Otherwise, $u'(x)\equiv 2\pmod 4$, in which case $x\in W_-$ and $s\equiv
7-u'\pmod 9$.
\end{lemma}
\begin{proof}
This can be read off the following table of degrees.
\begin{center}\renewcommand{\arraystretch}{1.1}
\begin{tabular}{|c|c|c|c|}
\hline element & $s$ & $u'$ & $t$
\\\hline\hline $h_{10}$ & 0 & $-2$ & $1$
\\\hline $b_{10}$ & 0 & 0 & $2$
\\\hline $w_n$ & 1 & $2(3^n-5)$ & $1$
\\\hline
\end{tabular}
\\\qedhere
\end{center}
\end{proof}
\begin{proposition}\label{d_4-d_8}
If $r\geq 2$ and $r\not\equiv 4\pmod 9$ or $r\not\equiv 8\pmod 9$, then $d_r =
0$. Furthermore, 
\begin{align*}
d_{4 + 9n}(W_+) & \subseteq W_- & d_{4 + 9n}(W_-) & =0
\\d_{8+9n}(W_+) & = 0 & d_{8+9n}(W_-) & \subseteq W_+.
\end{align*}
\end{proposition}
\begin{proof}
This is a degree argument, so we simplify to considering $d_r(x)$ where $x$ is a
monomial. First notice that $s(d_r(x)) + t(d_r(x)) = s(x) + t(x) + 1$. If $x\in
W_+$, then $s+t$ is even; if $x\in W_-$, then $s+t$ is odd. Thus
$d_r(W_+)\subseteq W_-$ and $d_r(W_-)\subseteq W_+$.

If $x\in W_+^{s,u'}$, then $d_r(x)\in W_-^{s+r,u'-6}$. If $d_r(x)\neq 0$, Lemma
\ref{s-possibilities} implies $s+u'\equiv 0\pmod 9$ and $s+r+u'-6\equiv 7\pmod
9$, so $r\equiv 4\pmod 9$. Similarly, if $x\in W_-^{s,u'}$, then $d_r(x)\in
W_+^{s+r,u'-6}$, which implies $r\equiv 8\pmod 9$ if $d_r(x)\neq 0$.
\end{proof}
In Section \ref{section:d_8}, we show that if $d_r(x) = h_{10}y$ is the first nontrivial
differential on $x\in W_+$, and $d_4(y)=h_{10}z$, then $d_{r+4}(h_{10}x) =
b_{10}z$. Combined with our complete calculation of $d_4$ in Section
\ref{section:d_4}, this determines the spectral sequence through the $E_9$ page.
We conjecture that the spectral sequence collapses at $E_9$.
The idea is that there is an operator $\partial:W_+\to W_+$ defined by
$\partial(x) = {1\over h_{10}}d_r(x)$ where $d_{r'}(x)=0$ for $r<r'$, and that
the spectral sequence essentially operates by taking Margolis homology of this
operator: if $x\in E_2$ supports a nontrivial $d_r$, then $d_r(x) = h_{10}
\partial(x)$, and $d_{r+4}(h_{10}x) = b_{10}\partial^2(x)$.

\begin{remark}
It is tempting to expect that Conjecture \ref{tilW-conj} comes from a map $k
\to P\cotensor_D\til{W}$, which would induce a map
$b_{10}^{-1}\Ext_P^*(k,k)\to b_{10}^{-1}\Ext_P^*(k,P\cotensor_D\til{W}) \isom
b_{10}^{-1}\Ext_D^*(k,\til{W})$ by the change of rings theorem. However, this is
not the case: $k\to P\cotensor_D \til{W}$ would factor through
$P\cotensor_D k$, which would mean that the map in $b_{10}^{-1}\Ext_P^*(k,-)$
would factor through $b_{10}^{-1}\Ext_P^*(k,P\cotensor_D k)\isom R$.
\end{remark}

\section{Identifying the $b_{10}$-periodic region}\label{section:periodicity-line}
%\renewcommand{\thesection}{B}
%\setcounter{theorem}{0}
%\setcounter{equation}{0}
%In this section we will use the following notation:
%\begin{itemize} 
%\item For a $P$-comodule $M$, $\underline{M}$ is the resolution of $M$ thought
%of as an object in $\Stable(P)$.
%\item $K(\xi_t^{p^s}) =
%b_{ts}^{-1}\underline{P\cotensor_{D[\xi_t^{p^s}]/(\xi_t^{p^s})^p}k}$
%\item $s(\xi_t^{p^s}) = {p\over 2}|\xi_t^{p^s}| = p^{s+1}(p^t-1)$ (see Section
%\ref{section:vanishing-lines})
%\end{itemize}\ 
In this section, let $p$ be an odd prime and let $k = \F_p$. The following
characterization of a $b_{10}$-periodic region in Ext is a consequence of
results of Palmieri that generalize the vanishing line theorems of Miller and
Wilkerson \cite{miller-wilkerson} to the stable category of comodules.

\begin{proposition} \label{b10-periodicity-line}
The localization map 
$\Ext^{s,t}_P(k,M)\to b_{10}^{-1}\Ext_P^{s,t}(k,M)$ is an isomorphism in the
range $s>{1\over p^3-p-1}(t-s) + c'$ for some constant $c'$.
\end{proposition}

Our main input is the following theorem, which Palmieri states for the Steenrod
dual $A$ instead of the algebra $P$ of dual reduced powers, as we do below.
The necessary changes in the case of $P$ follow immediately from the discussion in 
\cite[\S2.3.2]{palmieri-book}.%
\footnote{The only difference is that, over $A$,
one must also take into account the objects $Z(n)$ corresponding to $\tau_n$'s
as opposed to $\xi_t^{p^s}$'s, which do not come into
play over $P$.} 

Following Palmieri \cite[Notation 2.2.8]{palmieri-book}, define the \emph{slope} of $\xi_t^{p^s}$ to be:
$$s(\xi_t^{p^s}) = {1\over 2}p|\xi_t^{p^s}| = p^{s+1}(p^t-1).$$
Let $D[x]$ denote the truncated polynomial algebra $k[x]/x^p$. We have
$\Ext_D[\xi_t^{p^s}](k,k) = E[h_{ts}]\tensor k[b_{ts}]$. Let $K(\xi_t^{p^s}) =
b_{ts}^{-1} (P\cotensor_{D[\xi_t^{p^s}]}k)$, where the localization is defined by
taking a colimit of multiplication by $b_{ts}$ in $\Stable(P)$.
\begin{theorem} 
[{\cite[Theorem 2.3.1]{palmieri-book}}]\label{palmieri-vanishing-line}
Suppose $X$ is an object in $\Stable(P)$ satisfying the following conditions:
\begin{enumerate} 
\item There exists an integer $i_0$ such that $\pi_{i,*}X = 0$ if $i< i_0$,
\item There exists an integer $j_0$ such that $\pi_{i,j}X = 0$ if $j-i< j_0$,
\item There exists an integer $i_1$ such that the
homology of the cochain complex $X$ vanishes in homological degree $>i_1$. (In
particular, this is satisfied if $X$ is the resolution of a bounded-below
comodule.)
\end{enumerate}
Suppose $d = s(\xi^{p^{s_0}}_{t_0})$ (with $s_0 < t_0$) has the property that
$K(\xi^{p^s}_t)_{**}(X)=0$ for all $(s,t)$ with $s<t$ and $s(\xi^{p^s}_t) < d$.
Then $\pi_{**}X$ has a vanishing line of slope $d$: for some $c$, $\pi_{i,j}X=0$
when $j<di-c$.
\end{theorem}

\begin{figure}[H]
\rotatebox{270}{\begin{minipage}{0.75\textheight}
  \centering
    \includegraphics[width=\textheight]{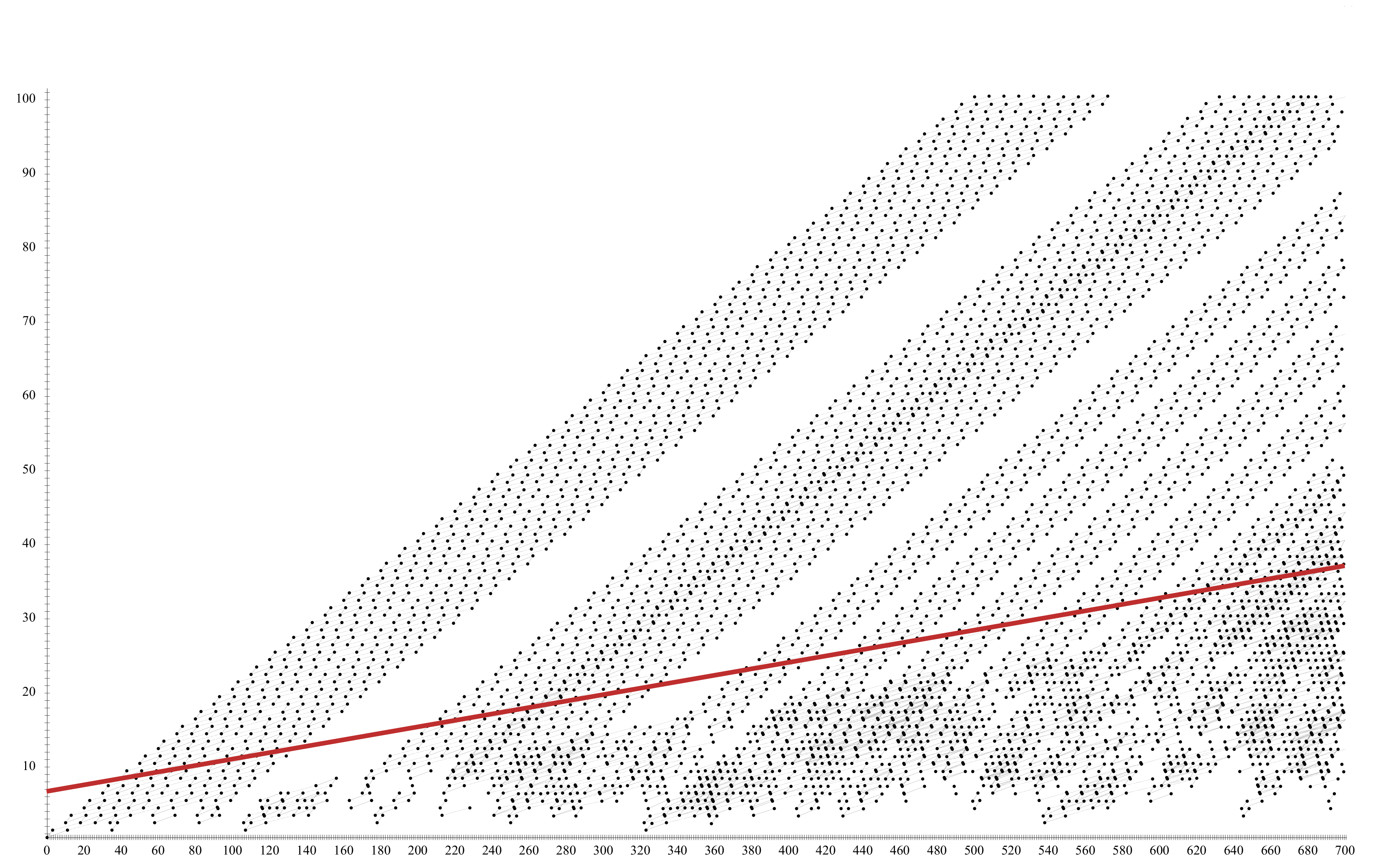}
    \captionof{figure}{%
	  Chart of $\Ext_P^*(\F_3,\F_3)$ with the line of Proposition
	  \ref{b10-periodicity-line} drawn in red: classes above the line are
	  $b_{10}$-periodic.
      \label{figure:Ext_P-chart}}
	  \end{minipage} }
\end{figure}

\begin{proof}[Proof of Proposition \ref{b10-periodicity-line}]
Let $M/b_{10}$ denote the cofiber in $\Stable(P)$ of $b_{10}\in \Ext_P^2(k,k)$,
thought of as a map $k\to k$ in $\Stable(P)$. It is not hard to check the
conditions (1)--(3) of Theorem \ref{palmieri-vanishing-line} for
$M/b_{10}$. We will apply the theorem with $d = s(\xi_2) = p^3-p$; note that
$\xi_2$ is the next $\xi_t^{p^s}$ with $s<t$ and higher slope than $\xi_1$, so
we just have to check $K(\xi_1)_{**}(M/b_{10})=0$.
This follows because the cofiber sequence
\begin{equation}\label{SES} M\too{b_{10}} M[2] \to
M/b_{10}[2] \end{equation}
gives rise to a long exact sequence in $K(\xi_1)_{**}$, and
multiplication by $b_{10}$ is an isomorphism on $K(\xi_1)_{**}(M)$ by
construction. So the theorem implies that there exists some $c$ such that
$\pi_{s,t}(M/b_{10})=0$ when $t<(p^3-p)s-c$.

Applying $\Ext$ to \eqref{SES}, we obtain
%$$ \u{\Ext_P^{s+1,t+|b_{10}|}(k,M/b_{10})}_{0\text{ if }t+|b_{10}| < (p^3-p)(s+1)-c} \to \Ext_P^{s,t}(k,k)\to \Ext_P^{s+2,t+|b_{10}|}(k,k)\to
%\u{\Ext_P^{s+2,t+|b_{10}|}(k,M/b_{10})}_{0\text{ if }t+|b_{10}| < (p^3-p)(s+2)-c}$$
$$ \Ext_P^{s+1,t+|b_{10}|}(k,M/b_{10}) \to \Ext_P^{s,t}(k,k)\to
\Ext_P^{s+2,t+|b_{10}|}(k,k)\to \Ext_P^{s+2,t+|b_{10}|}(k,M/b_{10})$$
where $|b_{10}| = 2p(p-1)$. Applying the vanishing condition for $M/b_{10}$
directly gives a region in which multiplication by $b_{10}$ is an isomorphism.
\end{proof}
In particular, at $p=3$, $b_{10}^{-1}\Ext_P(k,k)$ agrees with $\Ext_P(k,k)$
above a line of slope ${1\over 23}$ (see Figure \ref{figure:Ext_P-chart}).
In \cite[2.3.5(c)]{palmieri-book}, Palmieri gives an explicit expression for the
constant, which allows us to calculate the $y$-intercept to be $c'\approx
6.39$.

\newpage

\newpage

\section{$R$-module structure of $K(\xi_1)_{**}K(\xi_1)$ at $p>2$}
\label{section:cooperations}
In this section, we work at an arbitrary odd prime, and let $k = \F_p$ and $D
= \F_p[\xi_1]/\xi_1^p$.

In preparation for studying the $E_2$ page
$\Ext_{K(\xi_1)_{**}K(\xi_1)}(R,R)$, our goal for the next two sections is to study the Hopf algebra
$K(\xi_1)_{**}K(\xi_1)$, which in \eqref{cooperations-change-of-rings} we showed
is isomorphic to $b_{10}^{-1}\Ext_D(k,B)$.
Most of this section is devoted to giving an expression for $B$ as a
$D$-comodule. In the next section, we will obtain a more explicit description at
$p=3$, in which case we calculate the $E_2$ page.

\subsection{$D$-comodule structure of $B$}\label{section:B}
Note that $B$ is an algebra and a $P$-comodule, but not a
coalgebra. Let $\psi$ denote the $D$-coaction $B\to D\tensor B$ that comes from
composing the $P$-coaction $B\to P\tensor B$ with the surjection $P\to D$.

\begin{definition}\label{def-partial}
If we write $$\psi(x) = 1\tensor x + 
\xi_1\tensor a_1 + \xi_1^2\tensor a_2 + \dots + \xi_1^{p-1}\tensor a_{p-1}$$ for some
$a_i$'s, define $$\partial(x) : = a_1.$$
\end{definition}
For example, since $\Delta(\xi_n) = 1\tensor \xi_n + \xi_1\tensor \xi_{n-1}^p +
\dots$ (using the convention of Notation \ref{xi-antipode}), we have $\partial(\xi_n)=\xi_{n-1}^p$, and $\partial(\xi_{n-1}^p)=0$.
One can show using coassociativity that $a_k = {1\over k!}\partial^{k-1}a_1$. As
$\xi_1$ is dual to $P^0_1$ in the Steenrod algebra, the operator $\partial: P
\to P$ is dual to the operator $P^\vee \to P^\vee$ given by left
$P^0_1$-multiplication. In particular, $(P^0_1)^p=0$ implies $\partial^p = 0$.

%\begin{lemma}
%We have $a_k = {1\over k!}\partial^{k-1}a_1$.
%\end{lemma}
%\begin{proof}
%\fixme{Show $\partial a_{k-1} = k\cdot a_k$.} Iterating this, we obtain
%\begin{align*}
%\partial^k a_n  & = \partial^{k-1}((n+1)a_{n+1}) = \dots = (n+1)\dots
%(n+k)a_{n+k}.\qedhere
%\end{align*}
%\end{proof}
%\begin{corollary}
%$\psi(x) = 1\tensor x + \sum_{i=1}^{p-1} {1\over i!}\xi_1^i\tensor \partial^i x$
%\end{corollary}
\begin{lemma}\label{M-resolution}
$\Ext_D(k,M)$ is the cohomology of the chain complex $0\to
M\too{\partial}M\too{\partial^2}M\too{\partial}M \to \dots$, and
$b_{10}^{-1}\Ext_D(k,M)$ is the cohomology of the unbounded chain complex $\dots \to
M\too{\partial}M\too{\partial^2}M\too{\partial} M \to \dots$.
\end{lemma}

\begin{lemma}\label{leibniz}
We have $\partial(xy) = \partial(x)y + x\partial(y)$.
\end{lemma}
\begin{proof}
We have
\begin{align*}
\Delta(xy)  & = \Delta(x)\Delta(y) = (1\tensor x + \xi_1\tensor \partial x +
\dots)(1\tensor y + \xi_1\tensor \partial y + \dots)
\\ & = 1\tensor xy + \xi_1\tensor (y\partial x + x\partial y) + \dots. \qedhere
\end{align*}
\end{proof}

The structure theorem for modules over a PID says that modules over $D^\vee
\isom D$ decompose as sums of modules isomorphic to $\F_p[\xi_1]/\xi_1^i$ for $1\leq
i\leq p$. Dually, we have the following:
\begin{lemma}
Let $M(n)$ denote the $D$-comodule $\F_p[\xi_1]/\xi_1^{n+1}$. Then every
$D$-comodule splits uniquely as a direct sum of $D$-comodules isomorphic to $M(n)$
for $n\leq p-1$.
\end{lemma}
Note that $M(0)\isom \F_p$ and $M(p-1)\isom D$.
\begin{remark}
Since $\Ext_D^*(k,D)$ is a 1-dimensional vector space in homological degree 0
and zero otherwise, $b_{10}^{-1}\Ext_D^*(k,F)=0$ for any free $D$-comodule $F$.
If $0\leq i\leq p-2$, $\Ext_D^*(k,M(i))$ is 1-dimensional in every homological
degree.
%We will repeatedly use the fact that $\Ext_D^*(k,D)$ is a 1-dimensional
%$k$-vector space in homological degree 0 and zero otherwise, and for $i\in \{
%0,1 \}$,
%$\Ext_D^*(k,M(i))$ is 1-dimensional in homological degree $\geq 0$. As $b_{10}$
%is the generator of $\Ext_D^2(k,k)$, we have
%$b_{10}^{-1}\Ext_D^*(k,D)=0$, and $b_{10}^{-1}\Ext_D^*(k,M(i))$ is a
%1-dimensional $k$-vector space in every dimension. Furthermore, for any
%$D$-comodule $M$, the localization map $\Ext_D^*(k,M)\to
%b_{10}^{-1}\Ext_D^*(k,M)$ is an isomorphism in homological degree $>0$.
\end{remark}

The goal is to prove the following proposition.
\begin{proposition}\label{B-decomposition}
Define the indexing set $\mathscr{B}$ to be the set of monomials of the form
$\prod_{j=1}^n\xi_{i_j}^{e_j}$ such that $1\leq e_j\leq p-2$, and for $X\in
\mathscr{B}$, write $x_j(X) \colonequals \xi_{i_j}^{e_j}$ and $e_j(X) \colonequals e_j$.
Then there is a $D$-comodule isomorphism
$$ B\  \isom\ \dsums_{X\in \mathscr{B}}\,\tensors_{j=1}^n
M(e_j(X))_{x_j(X)}\ \ \dsum\ \  F$$
where $F$ is a free $D$-comodule, the tensor product is endowed with the
diagonal $D$-comodule structure, and $M(e)_{\xi_i^e}\colonequals \F_p\{
\xi_i^e, \partial \xi_i^e, \dots , \partial^e \xi_i^e
\}\isom M(e)$.
\end{proposition}
\begin{corollary}\label{cooperations-all-p}
We have an $R$-module isomorphism
\begin{align}
\label{Ext_D}b_{10}^{-1}\Ext^*_D(k,B) & \isom
b_{10}^{-1}\Ext^*_D\big(k,\dsums_{X\in \mathscr{B}} \tensors_{j=1}^n
M(e_j(X))_{x_j(X)}\big).
\end{align}
\end{corollary}
\begin{remark}
There is a formula due to Renaud \cite[Theorem 1]{renaud} that allows one to
decompose the tensor products $\tensors M(e_i)$ into a sum of the basic
comodules $M(n)$, but in general it is rather complicated; instead we will do
this in the next section only at $p=3$.
\end{remark}

If $e\leq p-1$ then $M(e)_{\xi_n^e}$ is a sub-$D$-comodule of $B$ with dimension
$e+1$. By the Leibniz rule (Lemma \ref{leibniz}) we have $$M(e+pf)_{\xi_n^{e+pf}} = \F_p\{
\xi_n^e\xi_n^{pf},
\partial(\xi_n^e)\xi_n^{pf},\dots,\partial^e(\xi_n^e)\xi_n^{pf} \}=
M(e)_{\xi_n^e}\tensor \F_p\{ \xi_n^{pf} \} $$
for $e\leq p-1$.
For any collection of $e_i \in \N$, define
\begin{equation}\label{T(xi_nxi_m)}T(\xi_{n_1}^{e_1}\dots \xi_{n_d}^{e_d})\colonequals
M(e_1)_{\xi_{n_1}^{e_1}}\tensor \dots \tensor M(e_d)_{\xi_{n_d}^{e_d}}.
\end{equation}
This is a sub-$D$-comodule spanned (as a vector space) by monomials of the form
$\partial^{k_1} (\xi_{n_1}^{e_1})\dots\partial^{k_d}(\xi_{n_d}^{e_d})$.
Clearly, $B = \sum_\attop{\text{monomials}\\\prod \xi_{n_i}^{e_i}\in B} T(\xi_{n_1}^{e_1}\dots
\xi_{n_d}^{e_d})$, but this is not a direct sum decomposition---any given
monomial appears in many different summands.
To fix this, we will study the poset of $T(X)$'s, and find that $B$ is a
direct sum of the maximal elements of that poset.

\begin{random}{Notation}
Define the notation $$\Big\langle\prod_{i\geq 1} \xi_i^{e_i}\cc\prod_{i\geq 2}
\xi_i^{f_i}\Big\rangle \colonequals
\prod \xi_i^{e_i}\ \prod \xi_{i-1}^{pf_i}.$$
\end{random}
(These are not formal products; they only make sense if
$e_i=0=f_i$ for all but finitely many $i$.) For example, we have
$\an{X\cc 1} = X$ for any monomial $X$, and 
$\an{1\cc\xi_n} =\xi_{n-1}^p = \partial(\xi_n)$.
%and if $X = \prod \xi_i^{f_i}$ and $X' = \prod \xi_{i-1}^{pf_i}$,
%\begin{align}
%\notag \partial(\an{\xi_a\xi_b\xi_c\ ;\ X})  & = \partial(\xi_a\xi_b\xi_c\cdot X')
%\\\notag & =\xi_{a-1}^p\xi_b\xi_c\cdot
%X' + \xi_a\xi_{b-1}^p\xi_c\cdot X' + \xi_a\xi_b\xi_{c-1}^p\cdot X'
%\\\label{abc}& = \an{\xi_b\xi_c\cc\xi_aX} + \an{\xi_a\xi_c\cc \xi_bX}+\an{\xi_a\xi_b\cc \xi_cX}.
%\end{align}
Expressions
$\an{\prod_{i\geq 2} \xi_i^{e_i}\cc\prod_{i\geq 2} \xi_i^{f_i}}$ 
represent elements of $B\subset P$, and conversely every element of $B$ has a
representation of this form (note that $\xi_1^p = \an{1\cc \xi_2}$).
Monomials in $B$ do not have unique expressions of the form $\an{X\cc Y}$: for
example, $\an{\xi_{n-1}^p\cc 1}=\an{1\cc\xi_n}$. 
\begin{lemma}
There is a bijection
\begin{equation} \label{eqrel-bij} \textstyle \cu{\text{monomials in
}B}\longleftrightarrow \cu{\an{\prod_{i\geq 2} \xi_i^{e_i}\cc \prod_{i\geq
2}\xi_i^{f_i}} \st e_i \leq p-1}. \end{equation}
\end{lemma}
Say that a bracket expression is \emph{admissible} if it is of the form on the
right hand side.
\begin{proof}
Given a monomial, the admissible bracket expression is the one with
the greatest number of terms on the right-hand side. 
\end{proof}
\begin{lemma}\label{move-terms}
If $X$ is a monomial with admissible bracket expression $\an{\prod
\xi_i^{e_i}\cc \prod \xi_i^{f_i}}$ and $Y$ is a monomial in $T(X)$, then $Y$
(up to invertible scalar) has admissible expression $\an{\prod\xi_i^{e_i-c_i}\cc \prod \xi_i^{f_i+c_i}}$
for a set of $c_i\geq 0$ that are zero for all but finitely many $i$.
\end{lemma}
The idea is that $Y$ is obtained from $X$ by moving terms from the left
to the right.
\begin{proof}
If $e\leq p-1$ then we have
$$ \partial^i (\xi_n^e)={e!\over (e-i)!}\xi_n^{e-i}\xi_{n-1}^{pi}. $$
By definition, $X = \prod_{i\geq 1} \xi_i^{e_i + pf_{i+1}}$ where $e_1=0$, and
\begin{align*}
Y  & = \prod \partial^{k_i}\xi_i^{e_i+pf_{i+1}} = \prod (\partial^{k_i}
\xi_i^{e_i})\xi_i^{pf_{i+1}} = \prod {e_i!\over
(e_i-k_i)!}\xi_i^{e_i-k_i+pk_{i+1}}\xi_i^{pf_{i+1}}
\\ & = \Big\langle\prod {e_i!\over (e_i-k_i)!} \xi_i^{e_i-k_i}\cc \prod
\xi_i^{k_i+f_i}\Big\rangle
\end{align*}
using the fact that $\partial \xi_i^p = 0$. So we can take $c_i=k_i$ in the
lemma statement.
\end{proof}

\begin{definition}
For monomials $X$ and $Y$, write $X\geq Y$ if $Y\in T(X)$.
\end{definition}
It is easy to check that this makes the set of monomials into a poset, and that
$X\geq Y$ if and only if $T(X)\supseteq T(Y)$.
\begin{lemma}\label{til-W}
Suppose $W$ is a monomial with admissible bracket expression $\ann{\prod
\xi_i^{e_i}\cc \prod \xi_i^{f_i}}$. Let
$\til{W} = \ann{\prod \xi_i^{c_i}\cc \prod \xi_i^{d_i}}$ where $c_i = \min \{ e_i
+ f_i, p-1 \}$ and $d_i = f_i - (c_i - e_i)$. Then $\til{W}$ is the maximal
object $\geq W$.
\end{lemma}
\begin{proof}
Let $X$ be an arbitrary monomial, written in its unique admissible bracket
expression. Then $X\geq W$ if and only $X$ can be obtained from $W$ by moving
terms in $W$ from the right to the left side of the bracket expression. Note
that $\til{W}$ is the bracket expression obtained by moving as many terms to the
left as possible while still keeping the resulting expression admissible. This
implies $\til{W}$ is maximal.
\end{proof}
%\begin{lemma}\label{til-X-til-Y}
%If $X \in T(W)$ is a nonzero monomial, then $\til{X} = \til{W}$.
%\end{lemma}
%\begin{proof}
%We just have to show that $\til{W}$ contains $X$ and is maximal with that
%property. But $T(X) \subset T(W)$, so $T(\til{W})\contains T(X)\ni X$;
%maximality is clear.
%\end{proof}
%\begin{lemma}
%If $T(X)\ints T(Y)\neq \{ 0 \}$, then $\til{X} = \til{Y}$. In particular,
%$T(\til{X}) = T(\til{Y})$ contains both $T(X)$ and $T(Y)$.
%\end{lemma}
%\begin{proof}
%Apply Lemma \ref{til-X-til-Y} to any nonzero monomial $W\in T(X)\ints T(Y)$.
%\end{proof}
Define an equivalence relation on monomials where $X\sim Y$ if $\til{X} =
\til{Y}$.
\begin{lemma} \label{B-direct-sum}
There is a direct sum decomposition $\displaystyle B \isom \dsums_\attop{\text{eq.
class}\\\text{reps. } X}T(\til{X})$.
\end{lemma}
\begin{proof}
I claim that $T(\til{X}) = \F_p\{ Y\st X\sim Y \}$; this follows from
the fact that, by definition, $T(\til{X})$ is generated by $Y$ such that $Y\leq
\til{X}$. So the direct sum decomposition comes from partitioning monomials into
their equivalence classes.
\end{proof}
Let $\mathscr{I}$ be the set of admissible bracket expressions $X$ such that
$\til{X} = X$. By Lemma \ref{til-W} we have the following.
\begin{lemma}
$\mathscr{I}$ is the set of admissible bracket expressions
$\an{\prod \xi_i^{e_i}\cc \prod \xi_i^{f_i}}$ such that $e_i \leq p-1$ and if
$e_i<p-1$ then $f_i=0$.
\end{lemma}
\begin{lemma}\label{cancel-f}
%Suppose $X = \ann{\prod \xi_i^{e_i}\cc \prod \xi_i^{f_i}}$ is in $\mathscr{I}$. Then
%$\ann{\prod \xi_i^{e_i}\cc 1}$ is also in $\mathscr{I}$, and 
%there is an isomorphism of $D$-comodules $T(\ann{\prod \xi_i^{e_i}\cc 1})\to T(X)$.
If $X = \ann{\prod \xi_i^{e_i}\cc \prod \xi_i^{f_i}}$ is an admissible
expression, there is an isomorphism of $D$-comodules $T(\ann{\prod
\xi_i^{e_i}\cc 1}) \isom T(X)$.
\end{lemma}
\begin{proof}
By Lemma \ref{move-terms}, every $Y$ in $T(X)$ has a bracket expression obtained from $X$ by moving terms
from the left to the right, so the right hand side of the bracket expression for
$Y$ is divisible by $\prod \xi_i^{f_i}$,
and so $Y$ is divisible by $u\colonequals\ann{1\cc \prod \xi_i^{f_i}} = \prod
\xi_{i-1}^{pf_i}$. So multiplication by $u$ gives a map $T(\an{\prod
\xi_i^{e_i}\cc 1}) \to T(X)$, and moreover from the above description of $Y\in
T(X)$ it is easy to see that this is a bijection. Finally, since
$\partial(u)=0$, this is an isomorphism of $D$-comodules.
\end{proof}
\begin{lemma}\label{T(X)-free-e}
If $X = \ann{\prod \xi_i^{e_i}\cc \prod \xi_i^{f_i}}$ is an admissible
expression such that $e_k=p-1$ for some $k$ then $T(X)$ is a free $D$-comodule.
\end{lemma}
\begin{proof}
By definition, we have $T(X) = \tensors
M(e_i)_{\xi_{n_i}^{e_i}}$ where the tensor product is endowed with the diagonal
$D$-comodule structure and $M(e_k)_{\xi_{n_k}^{e_k}}\isom M(p-1)\isom D$ by
assumption. After rearranging terms, it suffices to show that, for any
$D$-comodule $M$, there is a $D$-comodule isomorphism $D \tensor M \to
D\tensor M$ where the left hand side has a diagonal $D$-coaction and the right
hand side has a left coaction ($\psi(d\tensor m) = \sum d'm'\tensor
d''\tensor m''$ vs. $\psi(d\tensor m) = \sum d'\tensor d''\tensor m$ where
$\Delta(d)=\sum d'\tensor d''$ and $\psi(m)=\sum m'\tensor m''$). This
isomorphism is a variant of the shear isomorphism of Lemma \ref{shear}, and is
given by $d\tensor m\mapsto \sum dm'\tensor m''$.
\end{proof}
By Lemmas \ref{cancel-f} and \ref{T(X)-free-e}, we have:
\begin{corollary}\label{T(X)-free-f}
If $X = \ann{\prod \xi_i^{e_i}\cc \prod \xi_i^{f_i}}$ is an admissible bracket
expression in $\mathscr{I}$ such that $f_i\neq 0$ for any $i$, then $T(X)$ is
free as a $D$-comodule.
\end{corollary}
\begin{proof}[Proof of Proposition \ref{B-decomposition}]
From Lemma \ref{B-direct-sum} we have $B \isom \dsums_{X\in
\mathscr{I}}T(X)$, and by Corollary \ref{T(X)-free-f} there are free
$D$-comodules $F$ and $F'$ such that
\begin{align*}
B\ \ & \isom \dsums_{\ann{X\cc 1}\in \mathscr{I}}T(\ann{X\cc 1}) \dsum F\ =
\dsums_{\ann{X\cc 1}\in \mathscr{I}}T(X)\dsum F
\\ & \isom \dsums_\attop{\ann{X\cc 1}\text{ s.t.}\\e_i(X)\leq p-2}T(X)\dsum F'
\\ & = \dsums_{X\in \mathscr{B}} T(X)\dsum F'
\\ & \isom \dsums_{X\in\mathscr{B}} {\tensors_i}\,
M(e_i(X))_{x_i(X)} \dsum F'.\qedhere
\end{align*}
\end{proof}

We conclude this section with a useful lemma that simplifies checking relations
in certain $b_{10}$-local Ext groups of interest.
\begin{lemma}\label{B^6-lemma}
Let $I(n) = (\xi_1^{pn},\, \xi_2^{pn},\,\dots)B$. Then $I(p-1)$ is contained in the
free part of $B$ according to the decomposition in Proposition
\ref{B-decomposition}. In particular, if $x \in \Ext^*_P(k,P\cotensor_D I(p-1))$
then $x=0$ in $b_{10}^{-1}\Ext^*_P(k,P\cotensor_D B)$.
\end{lemma}
\begin{proof}
Consider an arbitrary monomial $q=\xi_n^{(p-1)p}X$ in
$I(p-1)$. If $X$ has an admissible
expression $\an{\prod \xi_i^{e_i}\cc \prod \xi_i^{f_i}}$ then $q$ has an
admissible expression $\an{\prod \xi_i^{e_i}\cc \xi_{n+1}^{p-1}\prod \xi_i^{f_i}}$.
By Lemmas \ref{B-direct-sum} and \ref{T(X)-free-e}, it suffices to show that
$\til{q} = \an{\prod \xi_i^{c_i}\cc \prod \xi_i^{d_i}}$ satisfies $c_k=p-1$ for
some $k$. Using the formula for $\til{q}$ in Lemma \ref{til-W}, we have
$c_{n+1}=p-1$.
\end{proof}
\begin{corollary}\label{B^6-cor}
Let $I(n)$ be as in Lemma \ref{B^6-lemma}. If $x\in \Ext_P^*(k,P\cotensor_D
(P\cotensor_D I(p-1)))$, then $x$ is zero in $b_{10}^{-1}
\Ext_P^*(k,P\cotensor_D (P\cotensor_D I(p-1)))$.
\end{corollary}

\section{Hopf algebra structure of $K(\xi_1)_{**}K(\xi_1)$ at $p=3$}\label{section:B-hopf-algebroid}
Henceforth we will work at $p=3$. This assumption will allow us to simplify the
formula for $K(\xi_1)_{**}K(\xi_1)$ obtained in Corollary
\ref{cooperations-all-p} and show that $K(\xi_1)_{**}K(\xi_1)$ is flat over
$K(\xi_1)_{**}$ (this is not true at higher primes), enabling us to calculate
the $E_2$ page \eqref{K(xi_1)-E_2-page} of the $K(\xi_1)$-based MPASS.
In particular, our goal is to show the following:

\begin{theorem}\label{thm:section-3-main}
At $p=3$, the ring of co-operations $K(\xi_1)_{**}K(\xi_1)$ is flat over $K(\xi_1)_{**}$, and moreover there is an isomorphism of Hopf algebras
$$ K(\xi_1)_{**}K(\xi_1) = R\tensor E[e_2,e_3,\dots] $$
for generators $e_n\in b_{10}^{-1}\Ext_D^1(k, B)$ in internal topological degree
$2(3^n+1)$. That is, $e_n$ is primitive, and $K(\xi_1)_{**}K(\xi_1)$ is exterior
as a Hopf algebra over $R=K(\xi_1)_{**}$.
\end{theorem}
Plugging this into the expression \eqref{K(xi_1)-E_2-page} for the $E_2$ page,
we obtain:
\begin{corollary}
The $E_2$ page of the $K(\xi_1)$-based MPASS for computing
$\pi_{**}(b_{10}^{-1} k)$ is
$$ E_2^{**} \isom R[w_2,w_3,\dots] $$
where $w_n = [e_n]$.
\end{corollary}

\begin{remark}
As $B$ is a $P$-comodule algebra,
there is a Hopf algebroid $(B, B\tensor B)$ in $\Stable(P)$,
where $B\tensor B$ carries the diagonal coaction of $P$ (see Section
\ref{section:overview}) and the comultiplication is given by
\begin{align*}
B\tensor B & \ttoo{-\tensor \eta\tensor -} B\tensor B\tensor B\isom
(B\tensor B)\tensor_{B} (B\tensor B).
\end{align*}
The Hopf algebroid above is given by applying
$b_{10}^{-1}\pi_{**}(-)=b_{10}^{-1}\Ext^*_P(k,-)$ to this one.
\end{remark}

\subsection{Vector space structure of $K(\xi_1)_{**}K(\xi_1)$ at $p=3$}
\label{section:vector-space-structure}
In the $p=3$ case, Corollary \ref{cooperations-all-p} reads
$$ K(\xi_1)_{**}K(\xi_1) \isom b_{10}^{-1}\Ext_D^*(k,B) \isom
b_{10}^{-1}\Ext_D^*\big(k,\hspace{-5pt}\dsums_{\substack{\text{monomials}\\\xi_{n_1}\dots
\xi_{n_d}\\n_i\neq n_j}}\tensors_{i=1}^d M(1)_{\xi_{n_i}}\big) $$
where the tensor product has a diagonal $D$-coaction.
It is easy to see directly that $M(1)\tensor M(1)\isom D\dsum
\Sigma^{0,|\xi_1|} k$. (Here we use bigraded notation for the shift for
consistency with viewing these objects in $\Stable(D)$, so $\Sigma^{0,|\xi_1|}$
denotes a shift of 0 in the homological dimension and $|\xi_1|$ in internal degree).
In particular,
$$ k\{ x,\partial x \} \tensor k\{ y,\partial y \}\, \isom\, k\{ xy,\ \partial(x)y +
x\partial(y),\ \partial(x)\partial(y)\} \dsum k\{ \partial(x)y-x\partial(y)
\}.$$
After inverting $b_{10}$, free comodules become zero, and the only basic types of
comodules are $M(0)=k$ and $M(1)$.

\begin{lemma}\label{M(0)=M(1)}
In $\Stable(D)$, we have an isomorphism
$$ b_{10}^{-1}M(1) \isom \Sigma^{-1,2|\xi_1|} b_{10}^{-1} M(0). $$
\end{lemma}
\begin{proof}
A representative for $M(1)$ in $\Stable(D)$
(i.e., an injective resolution for it) is $0\to
D\too{\partial^2}\Sigma^{0,2|\xi_1|}D\too{\partial}\Sigma^{0,3|\xi_1|}D\too{\partial^2}\Sigma^{0,5|\xi_1|}D\to \dots$, and so $b_{10}^{-1}
M(1)\colonequals \colim(M(1)\too{b_{10}}\Sigma^{2,-|b_{10}|}M(1)\to \dots)$ is represented by 
$$\dots\to
\Sigma^{0,-|\xi_1|}D\too{\partial}\hspace{-12pt}\u{D}_{\text{hom.deg.0}}\hspace{-12pt}\too{\partial^2}\Sigma^{0,2|\xi_1|}D\too{\partial}\Sigma^{0,3|\xi_1|}D\to\dots.$$
Similarly, $b_{10}^{-1}M(0)$ is represented by
\begin{align*}
\dots\to
\Sigma^{0,-2|\xi_1|}D\too{\partial^2}\hspace{-12pt}\u{D}_{\text{hom.deg.0}}\hspace{-12pt}\too{\partial}\Sigma^{0,|\xi_1|}D\too{\partial^2}\Sigma^{0,3|\xi_1|}D\to\dots.
\end{align*}
and so there is a degree-preserving isomorphism $b_{10}^{-1}M(1)\to
\Sigma^{-1,2|\xi_1|}b_{10}^{-1}M(0)$.
\end{proof}
(At arbitrary primes, the formula $b_{10}^{-1}M(n)\isom
\Sigma^{-1,(p-1)|\xi_1|}b_{10}^{-1} M(p-2-n)$ holds for the same reason.)
Therefore, if $M$ is a $D$-comodule, then $b_{10}^{-1} M\in \Stable(D)$ is a sum
of shifts of the unit object $k\isom M(0)$.
Remembering that $\Stable(D)$ was
constructed so that $\Hom_{\Stable(D)}(k,b_{10}^{-1}M) =
b_{10}^{-1}\Ext_D(k,M)$, we obtain the following K\"unneth isomorphism:
\begin{lemma}[K\"unneth isomorphism for $b_{10}^{-1}\Ext_D^*(\F_3,-)$]\label{kunneth}
If $M$ and $N$ are $D$-comodules, then $$b_{10}^{-1}\Ext^*_D(k,M\tensor N)\isom
b_{10}^{-1}\Ext^*_D(k,M)\tensor b_{10}^{-1}\Ext^*_D(k,N).$$
\end{lemma}
This only works at $p=3$, and is the essential reason we have made the
simplification of working at $p=3$.

Applying this to \eqref{Ext_D} we have the following.
\begin{corollary}\label{flatness}
We have an isomorphism
\begin{align*}
b_{10}^{-1}\Ext^*_D(k,B) & \isom \dsums_\attop{\text{monomials}\\\xi_{n_1}\dots \xi_{n_d}}
b_{10}^{-1}\Ext^*_D(k,\Sigma^{-d,2|\xi_1|}k_{\xi_{n_1}\dots\xi_{n_d}})
\end{align*}
where $\Sigma^{-d,2d|\xi_1|}k_{\xi_{n_1}\dots\xi_{n_d}}$ is the copy of
$\Sigma^{-d,2d|\xi_1|}k$ isomorphic to $\tensors_{i=1}^dM(1)_{\xi_{n_i}}$ under
Lemma \ref{M(0)=M(1)}.
In particular,
$K(\xi_1)_{**}K(\xi_1) = b_{10}^{-1}\Ext_D^*(k,B)$ is free over
$K(\xi_1)_{**} = b_{10}^{-1}\Ext_D(k,k)$.
\end{corollary}
So $b_{10}^{-1}\Ext_D(k,B)$ has $R$-module generators in bijection
with monomials of the form $\xi_{n_1}\dots \xi_{n_d}$ (where $n_i\neq n_j$ if
$i\neq j$). Now we will be more precise in choosing these generators.

\begin{lemma}\label{e(x)e(y)}
Suppose $N$ is a $D$-comodule algebra with sub-$D$-comodules
$k\{ x,\partial x \}\isom M(1)$ and $k\{ y,\partial y \}\isom M(1)$.
\begin{enumerate} 
\item The image of $\Ext_D^1(k,k\{ x,\partial x \})$ in $\Ext_D^1(k,N)$ is
generated by $$e(x) \colonequals [\xi_1]x - [\xi_1^2]\partial x.$$
\item We have $$ e(x) \cdot e(y) = b_{10}(y\partial x - x\partial y) $$
in the multiplication $\Ext^*_D(k,N)\tensor \Ext^*_D(k,N)\to \Ext^*_D(k,N)$
induced by the product structure on $N$. In particular, $e(x)^2 = 0$.
\item If the multiplication map embeds $k\{ x,\partial x \}\tensor k\{ y,\partial y \}$ in $N$ injectively, then $b_{10}^{-1}\Ext_D^2(k,k\{
x,\partial x \}\tensor k\{ y,\partial y \})\subset
b_{10}^{-1}\Ext_D^2(k,N)$ is a 1-dimensional vector space with generator $e(x)\cdot e(y)$.
\end{enumerate}
\end{lemma}
Since $\Ext^i_D(k,M) = b_{10}^{-1}\Ext^i_D(k,M)$ for $i>0$, note that this also
gives a generator of $b_{10}^{-1}\Ext^1_D(k,N)$.
\begin{proof}
Since $\Ext^1_D(k,M(1))$ is a 1-dimensional $k$-vector space, for (1) it suffices
to show that $e(x)$ is a cycle that is not a boundary. Indeed, since $dx =
[\xi_1]\partial x$ and $d(\partial x)=0$, we have $d(e(x))=-[\xi_1|\xi_1]\partial
x + [\xi_1|\xi_1]\partial x=0$, and
$e(x)$ is not in $d(C^0_D(k,k\{ x,\partial x \}))= d(k\{ x,\partial x \})$.

For (2), we use a special case of the cobar complex
multiplication formula in \cite[Proposition 1.2]{miller-localization}:
\begin{fact}
The multiplication $C^1_D(k,M)\tensor C^1_D(k,N)\to C^2_D(k,M\tensor N)$ is given by
$$ [\xi]m \tensor [\omega]n \mapsto \sum [\xi\tensor m'\omega](m''\tensor n). $$
\end{fact}
Thus the product $C^1_D(k,N)\tensor C^1_D(k,N)\to C^2_D(k,N\tensor N)\too{\mu}
C^2_D(k,N)$ takes $[\xi]m\tensor [\omega]n\mapsto \sum [\xi\tensor
m'\omega]m''n$.
Using this formula, we have:
\begin{align*}
e(x)\cdot e(y) & = [\xi_1|x]\cdot [\xi_1|y] - \br{\xi_1|x}\cdot
[\xi_1^2|\partial y] 
\\ & \hspace{40pt}- \br{\xi_1^2|\partial x}\cdot [\xi_1|y] + \br{\xi_1^2|\partial x}\cdot [\xi_1^2|\partial y] 
\\\br{\xi_1|x}\cdot [\xi_1|y] & = \sum [\xi_1|x'\xi_1]x''y =
[\xi_1|\xi_1]xy+[\xi_1|\xi_1^2](\partial x)y
\\\br{\xi_1|x}\cdot [\xi_1^2|\partial y] & =
\sum [\xi_1|x'\xi_1^2]x''\partial y = [\xi_1|\xi_1^2]x\partial y
\\\br{\xi_1^2|\partial x}\cdot [\xi_1|y] & =
\sum [\xi_1^2|(\partial x)'\xi_1](\partial x)''y =
[\xi_1^2|\xi_1](\partial x)y
\\\br{\xi_1^2|\partial x}\cdot [\xi_1^2|\partial y] & =
\sum [\xi_1^2|\xi_1^2(\partial x)'](\partial x)''\partial y =
[\xi_1^2|\xi_1^2]\partial x\partial y
\\d([\xi_1^2]xy) & =
2[\xi_1|\xi_1]xy-[\xi_1^2|\xi_1](\partial x)y-[\xi_1^2|\xi_1]x\partial y-[\xi_1^2|\xi_1^2]\partial x\partial y
\\ e(x)\cdot e(y) + d([\xi_1^2]xy) & =
[\xi_1|\xi_1^2](\partial x)y+[\xi_1^2|\xi_1](\partial x)y-[\xi_1|\xi_1^2]x\partial y-[\xi_1^2|\xi_1]x\partial y
\\ & = b_{10}((\partial x)y-x\partial y)
\end{align*}

For (3), note that there is a decomposition of $D$-comodules
\begin{align*}
k\{ x,\partial x \}\tensor k\{ y,\partial y \}  & \lu{\isom}{\too{\mu}} k\{ xy, x\partial y, (\partial x)y, (\partial x)(\partial y) \}
\\ & = k\{ xy, (\partial x)y + x\partial y, (\partial x)(\partial y) \}\dsum
k\{ (\partial x)y - x(\partial y) \}
\end{align*}
and since $\Ext^{*>0}_D(k,D)=0$, the quotient map
$$b_{10}^{-1} \Ext_D^2(k,k\{ x,\partial x \}\tensor k\{ y,\partial x \}) \isom
b_{10}^{-1} \Ext_D^2(k,k\{ x\partial y - (\partial x)y \})$$
is an isomorphism. By (2), $e(x)\cdot e(y)$ is a generator of the latter Ext
group.
\end{proof}

\begin{lemma}\label{e(x)y}
Suppose $N$ is a $D$-comodule algebra with sub-$D$-comodules $k\{ x,\partial x
\}\isom M(1)$ and $k\{ y \}\isom k$.
\begin{enumerate} 
\item The image of $\Ext_D^0(k, k\{ y \})$ in $\Ext^0_D(k,N)$ is generated by $y$.
\item We have
\begin{align*}
e(x)\cdot y  & = [\xi_1]xy - [\xi_1^2](\partial x)y = y\cdot e(x).
\end{align*}
\item If the multiplication map embeds $k\{ x,\partial x \}\tensor k\{ y \}$ in
$N$ injectively, then
$e(x)\cdot y$ is a generator of the 1-dimensional vector space $b_{10}^{-1}\Ext_D^1(k,k\{ x,\partial x
\}\tensor k\{ y \})$.
\end{enumerate}
\end{lemma}
\begin{proof}
(1) is clear. (2) follows from the cobar complex multiplication formulas
\begin{align*}
C^0_D(k,M)\tensor C^1_D(k,N) & \to C^1_D(k,M\tensor N)
 &   m\tensor [\xi]n & \mapsto [\xi](m\tensor n)
\\C^1_D(k,M)\tensor C^0_D(k,N) & \to C^1_D(k,M\tensor N)
 &   [\xi]n\tensor m & \mapsto [\xi](n\tensor m).
\end{align*}
For (3), note that $k\{ x,\partial x \}\tensor k\{ y \} = k\{ xy, (\partial x)y
\}$. Note that $(\partial x)y = \partial(xy)$. From Lemma \ref{e(x)e(y)},
$b_{10}^{-1}\Ext^1_D(k,k\{ xy, \partial(xy) \})$ is generated by $e(xy) =
[\xi_1]xy - [\xi_1^2]\partial(xy)=e(x)\cdot y$.
\end{proof}

\begin{definition}\label{def:e_n}
Define $e_n\colonequals e(\xi_n) = [\xi_1]\xi_n - [\xi_1^2]\xi_{n-1}^3$ as the
chosen generator of $b_{10}^{-1}\Ext_D^1(k,M(1)_{\xi_n})$.
\end{definition}

\begin{lemma}\label{e_n-explicit}
Under the change of rings isomorphism
$$b_{10}^{-1}\Ext_D(k,B)\isom b_{10}^{-1}\Ext_P(k,P\cotensor_D B)$$
the image of $e(x)$ in $\Ext^1_P(k,P\cotensor_D B)$ has cobar representative
$$[\xi_1](1|x)-[\xi_1^2](1|\partial x)+[\xi_1](\xi_1|\partial x)\in
\bar{P}\tensor (P\cotensor_D B).$$
\end{lemma}
\begin{proof}
The change of rings isomorphism $\Ext_D(k,M) \isom \Ext_P(k,P\cotensor_D M)$
works as follows: since $P$ is free over $D$, the functor $P\cotensor_D -$ is
exact, and so given
an injective $D$-resolution $M\to X^\bullet$ for $M$, the complex
$P\cotensor_D M\to P\cotensor_D X^\bullet$ is an injective $P$-resolution.
So we have $\Ext^i_D(k,M) \isom \Cotor^i_D(k,M) =  H^i(k\cotensor_D X^\bullet)$,
which agrees with $\Ext^i_P(k,P\cotensor_D M)\isom \Cotor^i_P(k,P\cotensor_D M) =
H^i(k\cotensor_P(P\cotensor_D X^\bullet)) \isom H^i(k\cotensor_DX^\bullet)$.

In particular, $\Ext_P(k,P\cotensor_D B)$ can be computed by applying
$k\cotensor_P -$ to the resolution
\begin{equation}\label{intermediate-resolution} P\cotensor_D\, C_D(k,B) =
\big(P\cotensor_D B \to P\cotensor_D (D\tensor B)\to P\cotensor_D (D\tensor
\bar{D}\tensor B)\to \dots\big). \end{equation}
By Lemma \ref{e(x)e(y)}, $e(x)$ has representative $[1|\xi_1]x - [1|\xi_1^2]\partial x \in D\tensor \bar{D}\tensor B$ in the $D$-cobar
\emph{resolution} for $B$, and so its representative in 
\eqref{intermediate-resolution} is
$1|1|\xi_1|x - \cdot 1|1|\xi_1^2|\partial x$.

But we wanted a representative in the cobar complex $C_P(k,P\cotensor_D B)$, so
we will write down part of an explicit map from the $P$-cobar resolution for
$P\cotensor_D B$ to \eqref{intermediate-resolution}:
$$ \xymatrix{
P\cotensor_D B\ar[d]\ar@{=}[r] & P\cotensor_D B\ar[d]
\\P\tensor (P\cotensor_D B)\ar[d]\ar[r]^-{f^0}\ar[d] & P\tensor
B\ar[d]
\\P\tensor \bar{P}\tensor
(P\cotensor_D B)\ar[r]^-{f^1}\ar[d] & P\tensor \bar{D}\tensor B\ar[d]
\\ P\tensor \bar{P}^{\tensor 2}\tensor (P\cotensor_D B)\ar[d]\ar[r] & P\tensor \bar{D}^{\tensor 2}\tensor B\ar[d]
\\\vdots & \vdots
}$$
By basic homological algebra, the map $f^*$ exists and is unique, so to
find $f^0$ and $f^1$ it suffices to find $P$-comodule maps that make the first
two squares commute. In particular, one can check that the maps
\begin{align*}
f^0(a|b|c) & = \epsilon(b)a|c
\\f^1(a|b|c|d) & = \epsilon(c)a|b|d
\end{align*}
make the diagram commute, and $z\colonequals
[1|\xi_1](1|x)+[1|\xi_1](\xi_1|\partial x) -
 [1|\xi_1^2](1|\partial x)$ is a cycle in $P\tensor \bar{P}\tensor
(P\cotensor_D B)$ such that $(k\cotensor_Pf)(z)=e(x)$.
\end{proof}

\subsection{Multiplicative structure}
\begin{proposition}
The summand
$$ b_{10}^{-1}\Ext^d_D(k,M(1)_{\xi_{n_1}}\tensor \dots \tensor
M(1)_{\xi_{n_d}})\subset b_{10}^{-1}\Ext_D^d(k,B)$$
is generated by the product $e_{n_1}\dots e_{n_d}$.
\end{proposition}
\begin{proof}
Since
\begin{align*}
b_{10}^{-1}\Ext_D^d(k,{\textstyle\tensors M(1)_{\xi_{n_i}}}) &
= \begin{cases} \Sigma^{d,0}\, b_{10}^{-1}\Ext_D^0(k,{\textstyle \tensors M(1)_{\xi_{n_i}}}) & d\text{ is even}\\
\Sigma^{d-1,0}\,b_{10}^{-1}\Ext_D^1(k,{\textstyle\tensors M(1)_{\xi_{n_i}}})& d\text{ is odd,}\end{cases}
\end{align*}
it suffices to show that $b_{10}^{-1}\Ext_D^0(k,M(1)_{\xi_{n_1}}\tensor \dots \tensor
M(1)_{\xi_{n_d}})$ is generated by $b_{10}^{-d/2}e_{n_1}\dots e_{n_d}$ when $d$
is even, and $b_{10}^{-1}\Ext_D^1(k,M(1)_{\xi_{n_1}}\tensor \dots \tensor
M(1)_{\xi_{n_d}})$ is generated by $b_{10}^{-(d-1)/2}e_{n_1}\dots e_{n_d}$ when
$d$ is odd. We proceed by induction on $d$. The base case $d=1$ is by definition.

\emph{Case 1: $d$ is even.} 
The tensor product
$M(1)_{\xi_{n_1}}\tensor\dots\tensor M(1)_{\xi_{n_{d-1}}}$ is
isomorphic to $M(1)\dsum F$ for a free summand $F$. By Lemma \ref{e(x)e(y)},
$b_{10}^{-1}\Ext^2_D(k,(M(1)_{\xi_{n_1}} \tensor \dots \tensor
M(1)_{\xi_{n_{d-1}}})\tensor M(1)_{\xi_{n_d}})$
is generated by $e(x)\cdot
e_{n_d}$ where $e(x)$ is a generator of
$b_{10}^{-1}\Ext^1_D(k,M(1)_{\xi_{n_1}}\tensor
\dots \tensor M(1)_{\xi_{n_{d-1}}})$. By the inductive hypothesis, we can take
$e(x)=b_{10}^{-(d-2)/2}e_{n_1}\dots e_{n_{d-1}}$. So
then $b_{10}^{-1}e(x)e_{n_d}=b_{10}^{-d/2}e_{n_1}\dots e_{n_d}$ is a generator for
$b_{10}^{-1}\Ext^0_D(k,M(1)_{\xi_{n_1}}\tensor \dots \tensor M(1)_{\xi_{n_d}})$.

\emph{Case 2: $d$ is odd.} 
In this case,
$M(1)_{\xi_{n_1}}\tensor\dots\tensor M(1)_{\xi_{n_{d-1}}}$ is
isomorphic to $k\dsum F$ for a free summand $F$. By Lemma \ref{e(x)y},
$b_{10}^{-1}\Ext^1_D(k,(M(1)_{\xi_{n_1}} \tensor \dots \tensor
M(1)_{\xi_{n_{d-1}}})\tensor M(1)_{\xi_{n_d}})$
is generated by $y\cdot
e_{n_d}$ where $y$ is a generator of
$b_{10}^{-1}\Ext^0_D(k,M(1)_{\xi_{n_1}}\tensor
\dots \tensor M(1)_{\xi_{n_{d-1}}})$. By the inductive hypothesis, we can take
$y =b_{10}^{-(d-1)/2}e_{n_1}\dots e_{n_{d-1}}$.
\end{proof}
Recall we defined $R = b_{10}^{-1}\Ext_D(k,k) = E[h_{10}]\tensor k[b_{10}^{\pm
1}]$.
\begin{corollary}
There is an $R$-module isomorphism
$b_{10}^{-1}\Ext_D^*(k,M(1)_{\xi_{n_1}}\tensor \dots \tensor
M(1)_{\xi_{n_d}}) \isom R\{ e_{n_1}\dots e_{n_d} \}$
where the generator $e_{n_1}\dots e_{n_d}$ is in degree $d$.
\end{corollary}
\begin{corollary}\label{exterior-R-alg}
The map $R\tensor E[e_2,e_3,\dots]\to b_{10}^{-1}\Ext_D^*(k,B)$ is an isomorphism of
$R$-algebras.
\end{corollary}

\subsection{Antipode}\label{subsection:antipode}
The antipode is the map induced on $\Ext$ by the swap map $\tau:B\tensor B
\to B\tensor B$. In order to get a useful formula for this map, we will need the
following basic properties of Hopf algebras.
\begin{fact}\label{hopf-facts}
Denote the coproduct on an element $x$ of a Hopf algebra by $\Delta(x) = \sum
x'\tensor x''$.
\begin{enumerate} 
\item (coassociativity) $\sum x'\tensor (x'')'\tensor (x'')'' = \sum
(x')'\tensor (x')'' \tensor x''$
\item $\sum c(x)'\tensor c(x)'' = \sum c(x'')\tensor c(x')$
\item $\sum c(x')x'' = \epsilon(x)$
\item $\sum \epsilon(x')\tensor x'' = 1\tensor x$
\end{enumerate}
\end{fact}

\begin{lemma}[Shear isomorphism] \label{shear}
Suppose $M$ is a left $P$-comodule, and $B\tensor M$ is given the diagonal
$P$-coaction: $\psi(b\tensor m) = \sum b'm'\tensor b''\tensor m''$ (where
$\psi(b) = \sum b'\tensor b''$ and $\psi(m) = \sum m'\tensor m''$). Then there
is an isomorphism $S_M:B\tensor M\to P\cotensor_D M$ (where $P$ coacts on the
left on $P\cotensor_D M$) sending $b\tensor m \mapsto \sum bm'\tensor m''$. It
has an inverse $S_M^{-1}: b\tensor m\mapsto \sum bc(m')\tensor m''$. 
\end{lemma}

In order to be able to apply Lemma \ref{B^6-lemma}, we
now obtain an explicit formula for the induced map $\tau'\colonequals S_{B}\circ
\tau\circ S_{B}^{-1}: P\cotensor_D B\to P\cotensor_D B$. This map is:
$$ \xymatrix{
B\tensor B\ar[r]^{\tau} & B\tensor B\ar[d]^-{S_{B}}
\\P\cotensor_D B\ar[u]^-{S_{B}^{-1}}\ar@{.>}[r]^{\tau'} & P\cotensor_D B
}\hspace{30pt} \xymatrix{
\sum xc(y')| y''\ar@{|->}[r] & \sum y''|xc(y')\ar@{|->}[d]
\\x|y\ar@{|->}[u] & \sum y''\cdot x'c(y')'|x''c(y')''
} $$
Using Fact \ref{hopf-facts} we have:
\begin{align*}
\tau'(x\tensor y)  & =  \sum y''\cdot x'c(y')'|x''c(y')'' 
\\ & = \sum x'y''c((y')'')|x''c((y')')
\\ & = \sum x'(y'')''c((y'')')|x''c(y')  & \text{coassociativity}
\\ & = \sum x' \epsilon(y'')|x''c(y')
\\ & = \sum x'|x''c(y)
\end{align*}

Since $(K(\xi_1)_{**},\, K(\xi_1)_{**}K(\xi_1))$ is a Hopf algebroid, the antipode is
multiplicative, so to determine it, it suffices to show:
\begin{proposition}\label{c(e_n)}
We have:
\begin{enumerate} 
\item $c(h)=h$
\item $c(e_n)=-e_n$.
\end{enumerate}
\end{proposition}
\begin{proof}
The antipode is given by the map $\tau'_*:\Ext_P^*(k,P\cotensor_D B)\to
\Ext_P^*(k,P\cotensor_D B)$ induced by $\tau'$, defined so that $\tau'_*([x_1|\dots|x_s]m) = [\xi_1|\dots|x_s]\tau'(m)$.
Since $h = [\xi_1](1|1)\in \Ext^1_P(k,P\cotensor_D B)$, we have $c(h)=\tau'_*(h) =
h$. For (2), we need an explicit formula for the antipode in the dual Steenrod
algebra:
\begin{fact}
[{\cite[Lemma 10]{milnor-steenrod-algebra}}]\label{milnor-antipode-formula}
Let $\Part(n)$ be the set of ordered partitions of $n$, $\ell(\alpha)$ the
length of the partition $\alpha$, and $\sigma_i(\alpha) =
\sum_{j=1}^i\alpha_j$ be the partial sum. Then we have:
$$ c(\xi_n)= \sum_{\alpha\in \Part(n)} (-1)^{\ell(\alpha)}
\prod_{i=1}^{\ell(\alpha)}\xi_{\alpha_i}^{p^{\sigma_{i-1}(\alpha)}}. $$
In particular, if $n\geq 2$ then $c(\xi_n) \equiv -\xi_n +
\xi_1\xi_{n-1}^p\pmod{\bar{P}^{p^2}P}$ and $c(\xi_{n-1}^p)\equiv
-\xi_{n-1}^p\pmod{\bar{P}^{p^2}P}$.
\end{fact}
Recall (Notation \ref{xi-antipode}) that we have defined $\xi_n$ to be the
antipode of its usual definition, so here we have $\Delta(\xi_n) =
\sum_{i+j=n}\xi_i\tensor \xi_j^{p^i}$.
(Since the antipode is a ring homomorphism, the formula in Fact
\ref{milnor-antipode-formula} is the same in either case.)

Combining this antipode formula with the formula for $e_n$ in Lemma
\ref{e_n-explicit} we have:
\begin{align*}
\tau'_*(e_n) & =
\tau'_*([\xi_1](1|\xi_n)-[\xi_1^2](1|\xi_{n-1}^3)+[\xi_1](\xi_1|\xi_{n-1}^3))
\\ & = [\xi_1](1|c(\xi_n)) - [\xi_1^2](1|c(\xi_{n-1}^3)) +
[\xi_1](\xi_1|c(\xi_{n-1}^3) + 1|\xi_1c(\xi_{n-1}^3))
\\ & = [\xi_1](-1|\xi_n + 1|\xi_1\xi_{n-1}^3 + 1|A) -
[\xi_1^2](-1|\xi_{n-1}^3+1|B) 
\\ & \hspace{40pt}+ [\xi_1](-\xi_1|\xi_{n-1}^3 + \xi_1|C - 1|\xi_1\xi_{n-1}^3
+ 1|D)
\\ & = -e_n + [\xi_1](1|A + \xi_1|C+1|D) - [\xi_1^2](1|B)
\end{align*}
for $A$, $B$, $C$, and $D$ in $\bar{P}^9P = I(3)$. By Lemma \ref{B^6-lemma} these
terms are zero in $b_{10}$-local cohomology, and $c(e_n)=\tau'_*(e_n)=-e_n$.
\end{proof}
\begin{corollary}\label{eta_R}
We have $\eta_L = \eta_R$; that is, the Hopf algebroid $(K(\xi_1)_{**},\, K(\xi_1)_{**}K(\xi_1))$
is, in fact, a Hopf algebra.
\end{corollary}
\begin{proof}
One of the axioms of a Hopf algebroid is $c\circ \eta_R = \eta_L$. Since
$\eta_L$ is just the inclusion of $R$ into $b_{10}^{-1}\Ext_D^*(k,B)$, its image is
invariant under the antipode $c$.
\end{proof}

\subsection{Comultiplication}
To define the comultiplication map $$b_{10}^{-1}\Ext_P(k,B\tensor B)\to
b_{10}^{-1}\Ext_P(k,B\tensor B)^{\tensor 2},$$
first consider the maps
$$ \Ext_P(k,B\tensor B)\ttoo{\alpha_*}  \Ext_P(k,B\tensor
B\tensor B) \bttoo{\beta} \Ext_P(k,B\tensor B)\tensor
\Ext_P(k,B\tensor B)$$
where $\alpha_*$ is the map on Ext induced by $\alpha:B^{\tensor 2}\to
B^{\tensor 3}$ with $\alpha:a\tensor b\mapsto a\tensor 1\tensor b$,
and $\beta$ is defined as the map in the factorization
\begin{equation}\label{beta-factorization} \xymatrix@R=30pt@C=15pt{
\Ext_P(k,B^{\tensor 2})\tensor \Ext_P(k,B^{\tensor
2})\ar[r]^-{\text{K\"unneth}}\ar[dr] &
\Ext_P(k,B^{\tensor 2}\tensor B^{\tensor 2})\ar[r]^-{-\tensor \mu\tensor
-} & \Ext_P(k,B^{\tensor 3})
\\ & \Ext_P(k,B^{\tensor 2})\tensor_{\Ext_P(k,B)} \Ext_P(k,B^{\tensor
2})\ar@{.>}[ru]_-\beta
}\end{equation}
It follows from the shear isomorphism (Lemma \ref{shear}) and the change of
rings theorem that
$\Ext_P(k,B\tensor M)\isom \Ext_P(k,P\cotensor_D M)\isom \Ext_D(k,M)$,
and the K\"unneth isomorphism for $b_{10}$-local cohomology over $D$
(Lemma \ref{kunneth}) implies that $\beta$ is an isomorphism after inverting
$b_{10}$. We define the comultiplication map on $b_{10}^{-1}\Ext_P(k,B\tensor
B)$ by $\Delta\colonequals \beta^{-1}\circ\alpha_*$.

In particular, flatness of $K(\xi_1)_{**}K(\xi_1)$ over $K(\xi_1)_{**}$ implies that $(K(\xi_1)_{**},\,
K(\xi_1)_{**}K(\xi_1))$ is a Hopf algebroid using the definitions of
comultiplication, antipode, counit, and unit above. In a Hopf algebroid, the
comultiplication is a homomorphism, and so to determine $\Delta$ explicitly it
suffices to determine $\Delta(e_n)$. We prove this in Proposition
\ref{e_n-primitive}. Lemma \ref{e_n-explicit} gives an expression for
$e_n$ in $\Ext_P^1(k,P\cotensor_D B)$, so we prefer to calculate $\Delta:
b_{10}^{-1} \Ext_P(k,B\tensor B)\to b_{10}^{-1}\Ext_P(k,B\tensor
B)^{\tensor 2}$ after composing with the shear isomorphism; that is, there is
a commutative diagram
$$ \xymatrix{
b_{10}^{-1}\Ext_P(k,B\tensor B)\ar[r]^-{\alpha_*}\ar[d]_{(S_{B})_*}
& b_{10}^{-1}\Ext_P(k,B\tensor B\tensor B)\ar[d]^{((\Id\tensor
S_{B})\circ
S_{B\tensor B})_*}
& b_{10}^{-1}\Ext_P(k,B\tensor B)^{\tensor
2}\ar[l]_-\beta\ar[d]^{S_{B}\tensor S_{B}}
\\b_{10}^{-1}\Ext_P(k,P\cotensor_D B)\ar[r]^-{\alpha'_*} &
b_{10}^{-1}\Ext_P(k,P\cotensor_D(P\cotensor_D B)) & b_{10}^{-1}\Ext_P(k,P\cotensor_D
B)^{\tensor 2}\ar[l]_-{\beta'}
}$$
and we will show that $\alpha'_*(e_n) = \beta'(1\tensor e_n + e_n\tensor 1)$ in
$b_{10}^{-1}\Ext_P(k,P\cotensor_D (P\cotensor_D B))$.
(We have chosen to use an extra application of the shear isomorphism on the
middle term in order to apply Corollary \ref{B^6-cor}.)

\begin{lemma}\label{alpha-beta-formulas}
If $a \in \Ext_P(k,P\cotensor_D B)$ has cobar representative $[a_1|\dots|a_s](p|q)$, we have
\begin{align*}
\alpha'_*(a)  & = \sum [a_1|\dots|a_s](p|q'|q'') 
\\\beta'(1\tensor a + a\tensor 1)  & = [a_1|\dots|a_s]{\textstyle(\sum
p'|p''|q+p|q|1)} 
\end{align*}
in $\Ext_P(k,P\cotensor_D (P\cotensor_D B))$.
\end{lemma}
So to check that $a$ is primitive after inverting $b_{10}$, it suffices to check
\begin{equation}\label{alpha=beta} \sum [a_1|\dots|a_s](p|q'|q'') - [a_1|\dots|a_s]{\textstyle(\sum
p'|p''|q+p|q|1)} = 0 \end{equation}in $b_{10}^{-1}\Ext_P(k,P\cotensor_D (P\cotensor_D
B))$.

\begin{proof}
By definition, $\alpha'$ is the map induced on $\Ext$ by the composition
%$$ 
%P\cotensor_DB\ttoo{S_{B}^{-1}}   B\tensorD B\ttoo{-\tensor 1\tensor -}  B\tensorD B\tensorD
%B\ttoo{S_{B\tensor B}}
% P\cotensor_D(B\tensorD
%B)\ttoo{P\cotensor_D S_{B}}  P\cotensor_D (P\cotensor_D B)
%$$
$$ \hspace{-42pt}\xymatrix@C=35pt{
\hspace{40pt}P\cotensor_DB\ar[r]^-{S_{B}^{-1}}  & B\tensor
B\ar[r]^-{-\tensor \eta\tensor -} & B\tensor B\tensor
B\ar[r]^-{S_{B\tensor B}}
& P\cotensor_D(B\tensor
B)\ar[r]^{P\cotensor_D S_{B}} & P\cotensor_D (P\cotensor_D B).
}$$
On elements, we have:
\begin{align*}
x|y\mapstto \sum xc(y')|y''  \mapstto \sum xc(y')|1|y''
 & \mapstto \sum xc(y')(y'')'|1|(y'')''
\\ & \mapstto  \sum xc(y')(y'')'|((y'')'')'|((y'')'')''= \sum x|y'|y''
\end{align*}
where the last equality is a coassociativity argument similar to the one at the
beginning of Section \ref{subsection:antipode}.
That is, we have $\alpha'(x\tensor y)= \sum x\tensor y'\tensor y''$, which implies
$$ \alpha'_*([a_1|\dots|a_s](p|q)) = \sum [a_1|\dots|a_s](p|q'|q''). $$

The map $\beta'$ comes from the bottom composition in
$$ \hspace{-5pt}\xymatrix@C=30pt{
\Ext_P(k,B^{\tensor 2})^{\tensor 2}\ar[d]_-{(S_{B})_*\tensor
(S_{B})_*}\ar[r]^-{\text{K\"unneth}} & \Ext_P(k,B^{\tensor
2}\tensor B^{\tensor 2})\ar[r]^-{(-\tensor \mu\tensor -)_*}\ar[d]_-{(S_{B}\tensor S_{B})_*}
& \Ext_P(k,B^{\tensor 3})\ar[d]^-{(S_{B\tensor B})_*}
\\\Ext_P(k,P\cotensor_D B)^{\tensor 2}\ar[r]^-{\text{K\"unneth}} &
\Ext_P(k,(P\cotensor_D B)\tensor (P\cotensor_D B))\ar[r]^-{\gamma_*} & \Ext_P(k,P\cotensor_D (P\cotensor_D B)).
}$$
We will only give an explicit expression for $\beta'$ on elements of the form
$1\tensor a$ and $a\tensor 1$, where $1$ denotes the unit
$1\tensor 1\in \Ext_P^0(k,P\cotensor_D B)$ and
$a=[a_1|\dots|a_s](p\tensor q)\in \Ext_P^s(k,P\cotensor_D B)$.
In \cite{miller-localization}, there is a full description of the
K\"unneth map $K$ on the level of cochains, but here all we need are the maps
$K:C_P^0(k,M)\tensor C_P^s(k,N)\to
C_P^s(k,M\tensor N)$ and $K:C_P^s(k,N)\tensor C_P^0(k,M)\to C_P^s(k,M\tensor
N)$. The former sends $m\tensor [a_1|\dots|a_s]n \mapsto
[a_1|\dots|a_s](m\tensor n)$ and the latter sends $[a_1|\dots|a_s]n\tensor
m\mapsto [a_1|\dots|a_s](n\tensor m)$. In particular, $K(1\tensor a) =
[a_1|\dots|a_s](1|1|p|q)$ and $K(a\tensor 1)=[a_1|\dots|a_s](p|q|1|1)$ in
$\Ext_P^s(k,(P\cotensor_D B)\tensor (P\cotensor_D B))$.

To determine $\beta'$, it remains
to determine the map $\gamma: (P\cotensor_D B)\tensor (P\cotensor_D B)\to
P\cotensor_D (P\cotensor_D B)$ induced by $-\tensor \mu\tensor
-$. This is accomplished by calculating the effect of shear isomorphisms as
follows:
$$ \xymatrix@C=35pt{
(B\tensor B)\tensor (B\tensor B)\ar[r]^-{-\tensor \mu\tensor -} & B^{\tensor 3}\ar[d]^-{S_{B\tensor B}}
\\(P\cotensor_D B)\tensor
(P\cotensor_D B)\ar[u]^-{S_{B}^{-1}\tensor S_{B}^{-1}}
 & P\cotensor_D (B\tensor B)\ar[r]^-{P\cotensor_D S_{B}} & P\cotensor_D (P\cotensor_D B)
}$$
$$\xymatrix{
\sum xc(y')|y''\tensor zc(w')|w''\ar@{|->}[r] & \sum xc(y')|y''zc(w')|w''\ar@{|->}[d]
\\ x|y\tensor z|w\ar@{|->}[u] &
{\let\scriptstyle\textstyle\substack{\sum xc(y')(y'')'z'c(w')'(w'')'\\\hspace{10pt}\tensor
(y'')''z''c(w')''\tensor (w'')''\\=\sum xz'|yz''c(w')|w''}}\ar@{|->}[r] & \sum
xz'|yz''|w.
}$$
That is, $\gamma(x|y\tensor z|w) = \sum xz'|yz''|w$,
which implies
\begin{align*}
\beta'(1\tensor a + a\tensor 1) & = \gamma_*K(1\tensor a + a\tensor 1) 
\\ & = \gamma_*([a_1|\dots|a_s](1|1|p|q + p|q|1|1))
\\ & = [a_1|\dots|a_s]\gamma(1|1|p|q + p|q|1|1)
\\ & = [a_1|\dots|a_s]{\textstyle(\sum p'|p''|q+p|q|1)}.\qedhere
\end{align*}
\end{proof}

\begin{proposition}\label{e_n-primitive}
The element $e_n$ is primitive.
\end{proposition}
\begin{proof}
We need to check the criterion \eqref{alpha=beta} for $a = e_n$.
Recall we had the formula
$$e_n=[\xi_1](1|\xi_n)-[\xi_1^2](1|\xi_{n-1}^3)+[\xi_1](\xi_1|\xi_{n-1}^3)\in
C_P^1(P\cotensor_D B)$$
from Lemma \ref{e_n-explicit}.  It suffices to check that
$\alpha'_*(e_n)  - \beta'_*(1\tensor e_n+e_n\tensor 1)$ is zero in
$b_{10}^{-1}\Ext_P(k,P\cotensor_D (P\cotensor_D B))$.
Using Lemma \ref{alpha-beta-formulas} we have:
\begin{align*}
\alpha'_*(e_n)  - \beta'_*(1\tensor e_n+e_n\tensor 1)
 & =\big([\xi_1](1|\Delta \xi_n)-[\xi_1^2](1|\Delta \xi_{n-1}^3) +
[\xi_1](\xi_1|\Delta\xi_{n-1}^3)\big) 
\\ & \hspace{20pt}- \big([\xi_1](1|1|\xi_n + 1|\xi_n|1) -
[\xi_1^2](1|1|\xi_{n-1}^3+1|\xi_{n-1}^3|1) 
\\ & \hspace{40pt}+ [\xi_1](1|\xi_1|\xi_{n-1}^3 + \xi_1|1|\xi_{n-1}^3+\xi_1|\xi_{n-1}^3|1)\big)
\\ & = [\xi_1]\hspace{-6pt}\sum_\attop{i+j=n\\2\leq i\leq n-1}\hspace{-5pt}1|\xi_i|\xi_j^{3^i} -
[\xi_1^2]\hspace{-6pt}\sum_\attop{i+j=n-1\\ 1\leq i\leq n-2}\hspace{-5pt}1|\xi_i^3|\xi_j^{3^{i+1}} +
[\xi_1]\hspace{-6pt}\sum_\attop{i+j=n-1\\1\leq i\leq n-2}\hspace{-5pt}\xi_1|\xi_i^3|\xi_j^{3^{i+1}}
\end{align*}
But all the remaining terms in the difference are in $C_P(P\cotensor_D
(P\cotensor_D I(3)))$ so by Corollary \ref{B^6-cor} they are zero in
$b_{10}$-local cohomology.
\end{proof}
\begin{proof}[Proof of Theorem \ref{thm:section-3-main}]
The flatness assertion was proved in Corollary \ref{flatness}.
Putting together Lemma \ref{exterior-R-alg}, Proposition \ref{c(e_n)}, Corollary
\ref{eta_R}, and Proposition \ref{e_n-primitive}, we
see that the map $R\tensor E[e_2,e_3,\dots]\to b_{10}^{-1} \Ext_D^*(k,B)$ is an
isomorphism of Hopf algebras.
\end{proof}
%\begin{theorem}
%The map $R\tensor E[e_2,e_3,\dots]\to b_{10}^{-1} \Ext_D^*(k,B)$ is an isomorphism of
%Hopf algebras. That is, the Hopf algebroid $(K(\xi_1)_{**},
%K(\xi_1)_{**}K(\xi_1))$ is an exterior Hopf algebra over $R$ on the generators
%$e_2,e_3,\dots$ where $e_n$ has internal degree $2(3^n+1)$.
%\end{theorem}

\section{Computation of $d_4$} \label{section:d_4}

\subsection{Overview of the computation}
In the previous section, we've shown that the $K(\xi_1)$-based MPASS computing $b_{10}^{-1}\Ext_P(k,k)$ has the form
$$ E_2^{**} = E[h_{10}]\tensor k[b_{10}^{\pm 1},w_2,w_3,\dots]\implies b_{10}^{-1}
\Ext_P^*(k,k) $$
where $w_n$ is represented in $E_1^{1,2(3^n+1)}$ by $e_n = [\xi_1]\xi_n -
[\xi_1^2]\xi_{n-1}^3 \in b_{10}^{-1}\Ext_D^1(k,\bar{B})$.
%\begin{random}{Notation}\label{stu}
%Define the gradings $s$ (MPASS filtration), $t$ (internal homological degree),
%and $u$ (internal (topological) degree) by 
%$$ E_1^{s,t,u} = \Ext_P^{t,u}(k, \bar{B}^{\tensor s}). $$
%Furthermore, we will find it convenient to use $u' = u-6(s+t)$.
%\end{random}
Recall that $d_r$ is a map $E_r^{s,t,u} \to E_r^{s+r,t-r+1,u}$, $w_n$ has degree
$(s,t,u) = (1,1,2(3^n+1))$, $h_{10}$ has degree $(0,1,4)$, and $b_{10}$ has
degree $(0,2,12)$. Furthermore, $u'(w_n) = 2(3^n-5)$, $u'(h_{10})=-2$, and
$u'(b_{10})=0$. In Proposition \ref{d_4-d_8}, we have shown that the next
nontrivial differential is $d_4$. In this section we will completely determine
this differential. We begin by recording some $d_4$'s in low degrees.

\begin{proposition}
We have the following:
\begin{align*}
d_r(h_{10}) & =0\text{ for }r\geq 2
\\d_r(w_2) & =0\text{ for }r\geq 2
\\d_4(w_3) & = \pm b_{10}^{-4}h_{10}w_2^5
\\d_4(w_4) & = \pm b_{10}^{-4}h_{10}w_2^2w_3^3.
\end{align*}
\end{proposition}
\begin{proof}
The first two facts can be seen directly in the cobar complex $C_P(k,k)$, using
the cobar representatives $h_{10}=[\xi_1]$ and $w_2 = [\xi_1|\xi_2] -
[\xi_1^2|\xi_1^3]$, which are permanent cycles.

The differentials on $w_3$ and $w_4$ were deduced from the chart of
$\Ext_P^*(k,k)$ up to the 700 stem that appears as Figure
\ref{figure:Ext_P-chart} (generated by the software \cite{yacop}). In Proposition
\ref{b10-periodicity-line}, we show that $\Ext^*_P(k,k)$ agrees with
$b_{10}^{-1}\Ext_P^*(k,k)$ in the range of dimensions depicted in the chart.
Thus we know which classes in $E_2 = R[w_2,w_3,\dots]$ in this range of
dimensions die in the spectral sequence, and, using multiplicativity of the
spectral sequence, this forces the differentials above.
\end{proof}

The goal of this section is to prove the following:
\begin{theorem}
For $n\geq 5$, there is a differential in the MPASS
$$ d_4(w_n) = \pm b_{10}^{-4}h_{10}w_2^2w_{n-1}^3. $$
\end{theorem}
Since the spectral sequence is multiplicative, this determines $d_4$.

The main idea is to use comparison with a spectral sequence computing
$b_{10}^{-1}\Ext_{P_n}(k,k)$, where
$$P_n = k[\xi_1, \xi_2, \xi_{n-2}, \xi_{n-1}, \xi_n]/(\xi_1^9,
\xi_2^3,\xi_{n-2}^{27}, \xi_{n-1}^9,\xi_n^3).$$
(The idea is that this is the smallest algebra in which the desired differential
can be seen.)
This is a quotient Hopf algebra of $P$ by the classification of such
(see \cite[Theorem 2.1.1.(a)]{palmieri-book}).
Here's a picture:

\begin{center}
\begin{tikzpicture} [scale=0.65]
\draw (-0.5,0) -- (12,0);
\draw (0,0) rectangle (1,1);
\draw (0,1) rectangle (1,2);
\draw (1,0) rectangle (2,1);
\draw (7,0) rectangle (8,1);
\draw (7,1) rectangle (8,2);
\draw (7,2) rectangle (8,3);
\draw (8,0) rectangle (9,1);
\draw (8,1) rectangle (9,2);
\draw (9,0) rectangle (10,1);
\node[scale=0.8] at (0.5,0.5) {$\xi_1$};
\node[scale=0.8] at (0.5,1.5) {$\xi_1^3$};
\node[scale=0.8] at (1.5,0.5) {$\xi_2$};
\node[scale=0.8] at (7.5,0.5) {$\xi_{n-2}$};
\node[scale=0.8] at (7.5,1.5) {$\xi_{n-2}^3$};
\node[scale=0.8] at (7.5,2.5) {$\xi_{n-2}^9$};
\node[scale=0.8] at (8.5,0.5) {$\xi_{n-1}$};
\node[scale=0.8] at (8.5,1.5) {$\xi_{n-1}^3$};
\node[scale=0.8] at (9.5,0.5) {$\xi_n$};
\end{tikzpicture}
\end{center}
Recall $B = P\cotensor_D k$; let $B_n = P_n\cotensor_D k$. We will refer to the
spectral sequence \ref{K(xi_1)-MPASS-M} with $\Gamma = P_n$ as the
$b_{10}^{-1}B_n$-based MPASS computing $b_{10}^{-1}\Ext_{P_n}(k,k)$, and use
$E_r(k,B_n)$ to denote its $E_r$ page. For example,
\begin{equation}\label{E_1-B_n} E_1(k, B_n) = b_{10}^{-1}\Ext_D^*(k,\bar{B}_n^{\tensor *}).
\end{equation}
Let $E_r(k,B)$ denote the $b_{10}^{-1}B$-based MPASS
for $b_{10}^{-1}\Ext_{P_n}(k,k)$ we have been focusing on. Then the diagram
$$ \xymatrix{
B\ar[r]\ar[d] & P\ar[d]\ar[r] & D\ar@{=}[d]
\\B_n\ar[r] & P_n\ar[r] & D
}$$
shows there is a map of spectral sequences $E_r(k,B)\to E_r(k,B_n)$.
\begin{lemma}\label{suffices-nonzero}
It suffices to show that $d_4(w_n)\neq 0$ in $E_4(k,B)$.
\end{lemma}
\begin{proof}
Since $s(d_4(w_n)) = 4 + s(w_n) = 5$, we know that
$d_4(w_n)$ is a linear combination of terms of the form
$b_{10}^Nh_{10}w_{k_1}\dots w_{k_5}$. We have
\begin{align*}
u'(w_n) & = u'(b_{10}^Nh_{10}w_{k_1}\dots w_{k_5}) + 6
\\2(3^n-5) & = -2 + \sum_{i=1}^5 2(3^{k_i}-5) + 6
\\3^n + 18 & = \sum_{i=1}^5 3^{k_i}
\end{align*}
Note that $k_i\geq 2$. Looking at this mod 27, we see that (at least) two of the
$k_i$'s have to be $= 2$, say $k_1$ and $k_2$. Then we have $3^n = 3^{k_3} +
3^{k_4} + 3^{k_5}$. The only possibility is $n-1=k_3=k_4=k_5$. So if
$d_4(w_n)\neq 0$ then $d_4(w_n)=b_{10}^Nh_{10}w_2^2w_{n-1}^3$, and checking
internal degrees shows $N=-4$.
\end{proof}

When we discuss $E_r(k,B_n)$ it will be easy to see that there is a class
$w_n\in E_2(k,B_n)$ which is the target of $w_n\in E_2(k,B)$ along the quotient
map.
$$ \xymatrix{
E_4(k,B)\ar[r]^{d_4}\ar[d] & E_4(k,B)\ar[d]
\\E_4(k,B_n)\ar[r]^{d_4} & E_4(k,B_n)
%}
%\hspace{40pt}\xymatrix{
%w_n\ar@{|->}[r]\ar@{|->}[d] & ??\ar@{|->}[d]
%\\w_n\ar@{|->}[r] & b_{10}^{-4}h_{10}w_2^2w_{n-1}^3
}$$
Lemma \ref{suffices-nonzero} says that it suffices to show $d_4(w_n)\neq 0$ in
$E_4(k,B_n)$, but it turns out to be the same amount of work to show the
following more attractive statement.
\begin{random}{Goal}
There is a differential $d_4(w_n) = \pm b_{10}^{-4}h_{10}w_2^2w_{n-1}^3$ in
$E_r(k,B_n)$. 
\end{random}
Using the same argument as Proposition \ref{d_4-d_8}, we know that $d_2 = 0 =
d_3$ in $E_r(k,B_n)$, so $h_{10}w_2^2w_{n-1}^3$ is not the target of an earlier
differential.
We will use the following strategy to show the desired differential in $E_r(k,B_n)$:
\begin{enumerate} 
\item Calculate $E_2(k,B_n)$ in a region and identify classes $w_2, w_{n-1},
w_n$ that are the targets of their namesake classes under the quotient map
$E_2(k,B)\to E_2(k,B_n)$.
\item Show that $b_{10}^{-1}\Ext^*_{P_n}(k,k)$ is zero in the stem of
$b_{10}^{-4}h_{10}w_2^2w_{n-1}^3$. This implies that $b_{10}^{-4}h_{10}w_2^2w_{n-1}^3$ either supports
a differential or is the target of a differential.
\item Show that $b_{10}^{-4}h_{10}w_2^2w_{n-1}^3$ is a permanent cycle in the
MPASS (so it must be the target of a differential) and show that, for degree
reasons, $w_n$ is the only element that can hit it. By looking at filtrations,
we see this differential is a $d_4$.
\end{enumerate}
In order to show (2), we introduce another spectral sequence for calculating
$b_{10}^{-1}\Ext^*_{P_n}(k,k)$, the Ivanovskii spectral sequence (ISS)
\cite{ivanovskii}. This is
the ($b_{10}$-localized version of the) dual of the May spectral sequence;
that is, it is the spectral sequence obtained by filtering
the cobar complex on $P_n$ by powers of the augmentation ideal.
(For example, $[\xi_1\xi_2|\xi_{n-1}^3]$ has filtration $2 + 3 = 5$.)

In Section \ref{section:ust-gradings} we will introduce notation and record
facts about gradings. In Section \ref{section:B_n-MPASS} we will compute
$E_1(k,B_n)$ and the relevant part of $E_2(k,B_n)$, and show (1) and (3)
assuming (2). In Section \ref{section:ISS-computation} we will calculate the
relevant part of the ISS and show (2). Convergence of the localized ISS is
discussed in Section \ref{section:ISS-convergence}.

\subsection{Notation and gradings}\label{section:ust-gradings}
Since much of the work in this section consists of degree-counting arguments, we
will now record how differentials and convergence affect the various gradings at
play. We emphasize a change of coordinates on degrees that simplifies degree
arguments by putting $b_{10}$ in degree zero.

\subsubsection*{MPASS gradings}
In Section \ref{section:overview}, we introduced the gradings $(s,t,u)$.
The differential has the form
$$d_r: E_r^{s,t,u} \to E_r^{s+r,t-r+1,u}$$ and a permanent cycle in
$E_r^{s,t,u}$ converges to an element in $b_{10}^{-1}\Ext_P^{s+t,u}(k,k)$.
We also introduced $u'\colonequals u-6(s+t)$.
We prefer to track $(u',s)$ instead of $(s,t,u)$, because
$u'(b_{10})=0=s(b_{10})$, so all classes in a $b_{10}$-tower have the same
$(u',s)$-degree. The differential under the change of coordinates has the form
$$d_r:E_r^{u',s} \to E_r^{u'-6,s+r}$$ and a permanent cycle in $E_r^{u',s}$
converges to an element in $b_{10}^{-1}\Ext_P^{a,b}(k,k)$ (where $b$ is internal topological degree and $a$ is homological degree) with $b-6a = u'$.

\begin{definition}\label{def-stem}
Let \emph{stem} in $b_{10}^{-1}\Ext_P^{a,b}(k,k)$ denote the quantity $b-6a$.
Then a permanent cycle in $E_r^{u',s}$ converges to an element in the $u'$ stem.
\end{definition}

Finally, define
$$u'' \colonequals u-6t.$$
This is only useful for looking at the $E_1$ page of the
MPASS, as $d_1$ fixes $u''$.

\subsubsection*{ISS gradings}
The \emph{Ivanovskii spectral sequence}
computing $b_{10}^{-1}\Ext_{P_n}(k,k)$ is the spectral sequence obtained
by filtering the cobar complex on $P_n$ by powers of the augmentation ideal. Let $E_r^{ISS}$
denote the $E_r$ page of the Ivanovskii spectral sequence.

We use slightly different grading conventions: classes have degree $(s,t,u)$
where $s$ is ISS filtration, $t$ denotes degree in the cobar complex, and $u$
denotes internal topological degree (as in the MPASS). The differential has the
form
$$d_r^{ISS}: E_r^{s,t,u} \to E_r^{s+r,t+1,u}$$
and a permanent cycle in $E_r^{s,t,u}$ converges to an element in $b_{10}^{-1} \Ext_P^{t,u}(k,k)$.

We will use the change of coordinates
$$u'\colonequals u-6t$$ 
which is designed so that $u'(b_{10})=0$. (This has a different formula from the
MPASS change of coordinates simply because $(s,t,u)$ correspond to different
parameters here.) The differential has the form
$$d_r^{ISS}:E_r^{u',s}\to E_r^{u'-6,s+r}$$
and a permanent cycle in $E_r^{u',s}$ converges to an element in
$b_{10}^{-1}\Ext_P^{a,b}(k,k)$ with $u' = b-6a$, i.e. an element in the $u'$ stem.

Note that $u'$ has different formulas for the MPASS and ISS, but in both
spectral sequences $u'$ corresponds to stem, with the definition given above.
Now we will introduce another grading on $P_n$ (for $n\geq 5$) preserved by the
comultiplication.

\subsubsection*{Extra grading on $P_n$}
Let $P'_n =
k[\xi_1,\xi_2,\xi_{n-2}^3,\xi_{n-1},\xi_n]/(\xi_1^9,\xi_2^3,\xi_{n-2}^{27},\xi_{n-1}^9,\xi_n^3)$.
Note that every monomial in $P_n$ can be written $\xi_{n-2}^ex$ where $e\in \{
0,1,2 \}$ and $x\in P'_n$.
\begin{lemma}
For $n\geq 5$, $P'_n$ is a sub-coalgebra of $P_n$.
\end{lemma}
\begin{proof}
This is clear from the comultiplication formulas
\begin{align}
\label{Delta(xi_n)}\Delta(\xi_n) & = 1\tensor \xi_n+\xi_1\tensor \xi_{n-1}^3+\xi_2\tensor \xi_{n-2}^9
\\\notag \Delta(\xi_{n-1}) & = 1\tensor \xi_{n-1} + \xi_1\tensor
\xi_{n-2}^3+\xi_{n-1}\tensor 1
\\\notag \Delta(\xi_{n-2}^3) & = 1\tensor \xi_{n-2}^3 + \xi_{n-2}^3\tensor 1
\end{align}
and the assumption $n\geq 5$ guarantees that $\xi_1,\xi_2\neq \xi_{n-2}$.
\end{proof}
\begin{proposition}\label{def-alpha}
Let $n\geq 3$.
There is an extra grading $\alpha$ on $P_n$ that respects the comultiplication,
defined by the property that it is multiplicative on $P'_n$, and
\begin{align*}
\alpha(\xi_1)  & = \alpha(\xi_2) = \alpha(\xi_{n-2}) = 0,
\\\alpha(\xi_{n-2}^3)  & = \alpha(\xi_{n-1}) = 3,
\\\alpha(\xi_n) & =9,
%\\\alpha(\xi_{n-2}^ex)  & = e\alpha(\xi_{n-2})+\alpha(x) \text{ for }e\in \{ 0,1,2
%\} \text{ and } x\in P'_n.
\\\alpha(\xi_{n-1}^ex) & = \alpha(x) \text{ for }e\in \{ 0,1,2 \} \text{ and }
x\in P'_n.
\end{align*}
\end{proposition}
\begin{proof}
First we check that $\alpha$ respects the comultiplication when restricted to
$P'_n$. Since it is defined to be multiplicative on $P'_n$, it suffices to check
that $\alpha(y) = \alpha(\Delta y)$ for $y$ as each of the multiplicative
generators. This is clear from the comultiplication formulas
\eqref{Delta(xi_n)}.

Now suppose $y = \xi_{n-2}x$ where $x\in P'_n$. We have
\begin{align*}
\Delta(\xi_{n-2}x) & = (1\tensor \xi_{n-2} + \xi_{n-2}\tensor 1)\Delta x = \sum
\big(x'\tensor x''\xi_{n-2} + x'\xi_{n-2}\tensor x''\big)
\end{align*}
and the $\alpha$ degrees of both sides agree since $P'_n$ is a coalgebra.
Similarly, if $y=\xi_{n-2}^2x$ for $x\in P'_n$, we have 
\begin{align*}
\alpha(\Delta y)  & =
\alpha((1\tensor \xi_{n-2}^3 + 2\xi_{n-2}\tensor \xi_{n-2}+\xi_{n-2}^2\tensor
1)(\Delta x)) 
\\ & = \alpha\big(\sum x'\tensor \xi_{n-2}^2x'' + 2\xi_{n-2}x'\tensor \xi_{n-2}x'' +
\xi_{n-2}^2x'\tensor x''\big)
= \alpha(\Delta x).\qedhere
\end{align*}
\end{proof}

\subsection{The $E_2$ page of the $b_{10}^{-1}B_n$-based MPASS}
\label{section:B_n-MPASS}
The goal of this section is to prove the following:
\begin{proposition}\label{section3-result}
If $b_{10}^{-4}h_{10}w_2^2w_{n-1}^3$ is the target of a differential in the
$b_{10}^{-1}B_n$-based MPASS calculating $b_{10}^{-1}\Ext_{P_n}^*(k,k)$, that
differential must be $$ d_4(w_n) = \pm b_{10}^{-4}h_{10}w_2^2w_{n-1}^3. $$
\end{proposition}

The main task is to calculate enough of $E_2(k,B_n)$ to do a degree-counting
argument (Proposition \ref{almost-section3-result}), where
$$B_n = P_n\cotensor_D k = k[\xi_1^3,\xi_2,\xi_{n-2},
\xi_{n-1},\xi_n]/(\xi_1^9,\xi_2^3,\xi_{n-2}^{27},\xi_{n-1}^9,\xi_n^3).$$
As in the calculation of the $E_2$ page of the $b_{10}^{-1}B$-based MPASS
(Section \ref{section:B-hopf-algebroid}), the K\"unneth formula for the functor
$b_{10}^{-1}\Ext^*_D(k,-)$ (Lemma \ref{kunneth}) guarantees flatness of
$(b_{10}^{-1}B_n)_{**} (b_{10}^{-1} B_n)$ over $(b_{10}^{-1}B_n)_{**}$. So we
can use the formula
\begin{equation}\label{E_2-B_n} E_2 \isom \Ext^*_{(b_{10}^{-1} B_n)_{**}b_{10}^{-1} B_n}((b_{10}^{-1}B_n)_{**},
(b_{10}^{-1}B_n)_{**}) \end{equation}
where $(b_{10}^{-1} B_n)_{**} = b_{10}^{-1}\Ext_{P_n}^*(k,B_n) =R$
and $(b_{10}^{-1}B_n)_{**}(b_{10}^{-1} B_n) = b_{10}^{-1}\Ext^*_{P_n}(R,
B_n^{\tensor 2}) \isom b_{10}^{-1}\Ext^*_D(k,B_n)$ by the change of
rings theorem. We will
simultaneously determine the vector space structure and the comultiplication on
$b_{10}^{-1} \Ext_D(k,B_n)$.
\begin{remark}\label{MPASS-B_n-filtration-coincide}
By \eqref{E_1-B_n} and the K\"unneth formula mentioned above, we have 
$$E_1^{s,*}(k,B_n) \isom b_{10}^{-1}\Ext_D(k,\bar{B}_n)^{\tensor s}$$
and so the coproduct on $b_{10}^{-1}\Ext_D^*(k,B_n)$ coincides with $d_1$ on
$E_1^{1,*}$.
\end{remark}

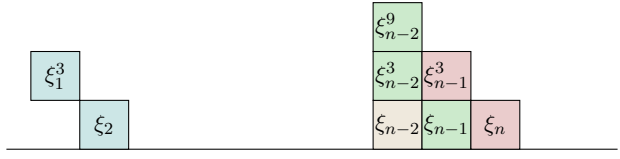
\begin{figure}[H] 
\begin{center}
\begin{tikzpicture} [scale=0.65]
\draw (-0.5,0) -- (12,0);
\draw[fill=teal!20] (0,1) rectangle (1,2);
\draw[fill=teal!20] (1,0) rectangle (2,1);
\draw[fill=gold!20] (7,0) rectangle (8,1);
\draw[fill=mgreen!20] (7,1) rectangle (8,2);
\draw[fill=mgreen!20] (7,2) rectangle (8,3);
\draw[fill=mgreen!20] (8,0) rectangle (9,1);
\draw[fill=maroon!20] (8,1) rectangle (9,2);
\draw[fill=maroon!20] (9,0) rectangle (10,1);
\node[scale=0.8] at (0.5,1.5) {$\xi_1^3$};
\node[scale=0.8] at (1.5,0.5) {$\xi_2$};
\node[scale=0.8] at (7.5,0.5) {$\xi_{n-2}$};
\node[scale=0.8] at (7.5,1.5) {$\xi_{n-2}^3$};
\node[scale=0.8] at (7.5,2.5) {$\xi_{n-2}^9$};
\node[scale=0.8] at (8.5,0.5) {$\xi_{n-1}$};
\node[scale=0.8] at (8.5,1.5) {$\xi_{n-1}^3$};
\node[scale=0.8] at (9.5,0.5) {$\xi_n$};
\end{tikzpicture}
\end{center}
\caption{Illustration of the decomposition of $B_n$ into tensor factors}
\label{figure:colored-boxes}
\end{figure}

We can write $B_n$ as a tensor product
$$ B_n = k[\xi_2,\xi_1^3]/(\xi_2^3,\xi_1^9)\tensor
k[\xi_{n-2}]/\xi_{n-2}^3\tensor k[\xi_{n-1},\xi_{n-2}^3]/(\xi_{n-1}^3, \xi_{n-2}^{27})\tensor
k[\xi_n,\xi_{n-1}^3]/(\xi_n^3,\xi_{n-1}^9) $$
illustrated in Figure \ref{figure:colored-boxes}.

Since we have a K\"unneth formula for $b_{10}^{-1}\Ext^*_D(k,-)$, it suffices to
apply this functor to each of the four factors of $B_n$ above.

%Start with the first piece.
\subsubsection*{Factor 1: $\color{mblue}k[\xi_2,\xi_1^3]/(\xi_2^3,\xi_1^9)$}
As a $D$-comodule, this decomposes as:
\begin{equation}\label{D-comodule_2} k[\xi_2,\xi_1^3]/(\xi_2^3,\xi_1^9) \isom 
\u{k\{ 1 \}}_{\isom k} \dsum \u{k\{ \xi_2,\xi_1^3 \}}_{\isom M(1)}\dsum \u{k\{
\xi_2^2,\xi_1^3\xi_2,\xi_1^6 \}}_{\isom D}\dsum \u{k\{ \xi_1^3\xi_2^2,
\xi_1^6\xi_2 \}}_{\isom M(1)}\dsum \u{k\{ \xi_1^6\xi_2^2 \}}_{\isom
k}.\end{equation} (Recall $M(1)$ was defined to be the $D$-comodule
$k[\xi_1]/\xi_1^2$, and every $D$-comodule is a sum of copies of $k$, $M(1)$,
and $D$.) As a module over $R\colonequals E[h_{10}]\tensor k[b_{10}^{\pm 1}]$,
this is generated by a class $e_2 = e(\xi_2)$ in $b_{10}^{-1}\Ext^1_D(k,k\{
\xi_2,\xi_1^3 \})$, a class $f_{20} = e(\xi_1^3\xi_2^2)$ in $b_{10}^{-1}
\Ext^1_D(k,k\{ \xi_1^3\xi_2^2,\xi_1^6\xi_2 \})$, and a class $c_2$ in
$b_{10}^{-1} \Ext^0_D(k,k\{ \xi_1^6\xi_2^2 \})$. As $b_{10}^{-1}\Ext_D^*(k,D) =0$,
we may ignore the free summands.

Using Lemma \ref{e(x)e(y)}, we can give explicit representatives for the classes
in $b_{10}^{-1}\Ext_D^*(k,k[\xi_2,\xi_1^3]/(\xi_2^3,\xi_1^9))$ coming from the decomposition \eqref{D-comodule_2}:
\begin{align*}
e_2 \colonequals e(\xi_2)  & = [\xi_1]\xi_2 - [\xi_1^2]\xi_1^3 \in \Ext_D^*(k,
k[\xi_2,\xi_1^3]/(\xi_2^3,\xi_1^9))
\\f_{20} \colonequals e(\xi_1^3\xi_2^2)  & = [\xi_1]\xi_1^3\xi_2^2 + [\xi_1^2]\xi_1^6\xi_2
\\c_2 & = \xi_1^6\xi_2^2
\end{align*}
satisfying relations $e_2^2 = 0=f_{20}^2$ and $b_{10}c_2 = e_2f_{20}$.

\begin{lemma}\label{e_2-f_20}
The classes $e_2$ and $f_{20}$ are primitive in
the coalgebra $b_{10}^{-1}\Ext_D^*(k,B_n)$.
\end{lemma}
\begin{proof}
As described in Section \ref{section:main-tool}, we can interpret the
MPASS as a filtration spectral sequence on the cobar complex $C_{P_n}(k,k)$, where
$[a_1|\dots|a_s]$ is in filtration $n$ if $\geq n$ $a_i$'s are in
$\bar{B}_nP_n$. The elements $e_2$ and $f_{20}$ correspond to elements in
$F^1/F^2C_{P_n}^2(k,k)$ with the same formulas, and by Remark \ref{MPASS-B_n-filtration-coincide} it suffices to show that
$d_1(e_2)=0=d_1(f_{20})$ in the filtration spectral sequence.
One checks explicitly that
$d_{\text{cobar}}(e_2)=0$, so it is a permanent cycle. This is not true of
$f_{20}$, but we can write down explicit correcting terms in higher filtration:
$$ f_{20}\equiv \til{f}_{20}\colonequals [\xi_2|\xi_2^2] + [\xi_2^2|\xi_2] -
[\xi_1\xi_2|\xi_2\xi_1^3]+ [\xi_1\xi_2^2|\xi_1^3]+[\xi_1^2\xi_2|\xi_1^6] +
[\xi_1^2|\xi_2\xi_1^6] + [\xi_1|\xi_2^2\xi_1^3]$$
and then check that $d_{\text{cobar}}(\til{f}_{20})=
[\xi_1^3|\xi_1^6|\xi_1^3]+[\xi_1^3|\xi_1^3|\xi_1^6]$. This has filtration 3, and
so $d_1(f_{20})=0$.
\end{proof}
So we've proved:
\begin{proposition}\label{MPASS-E_1-P(1)}
There is an isomorphism of Hopf algebras
$$b_{10}^{-1}\Ext_D(k,k[\xi_2,\xi_1^3]/(\xi_2^3,\xi_1^9))\isom R\tensor E[e_2,f_{20}]$$
where $e_2$ and $f_{20}$ are primitive.
\end{proposition}
We can summarize the degree information as follows:
\begin{center}\renewcommand{\arraystretch}{1.1}
\begin{tabular}{|c|c|c|c|c|c|}
\hline element & $s$ & $t$ & $u$ & $u''=u-6t$ & $\alpha$
\\\hline\hline 1 & 0 & 0 & 0 & 0 & 0
\\\hline $h_{10}$ & 0 & 1 & 4 & $-2$ & 0
\\\hline $b_{10}$ & 0 & 2 & 12 & 0 & 0
\\\hline\hline $e_2 = [\xi_1]\xi_2 - [\xi_1^2]\xi_1^3$ & 1 & 1 & 20 & 14 &
0
\\\hline $f_{20} =  [\xi_1]\xi_1^3\xi_2^2 + [\xi_1^2]\xi_1^6\xi_2$
& 1 & 1 & 48 & 42 & 0
\\\hline $c_2= \xi_1^6\xi_2^2$ & 1 & 0 & 56 & 56 & 0
\\\hline
\end{tabular}
\end{center}

%Next we record the decomposition of the second factor of $B_n$, namely
%$k[\xi_{n-2}]/\xi_{n-2}^3$.
\subsubsection*{Factor 2: $\color{brown}k[\xi_{n-2}]/\xi_{n-2}^3$}
This decomposes as $k\{ 1 \}\dsum k\{ \xi_{n-2} \}\dsum k\{ \xi_{n-2}^2 \}$ so we
have three $R$-module generators:
\begin{center}\renewcommand{\arraystretch}{1.1}
\begin{tabular}{|c|c|c|c|c|c|}
\hline element & $s$ & $t$ & $u$ & $u''=u-6t$ & $\alpha$
\\\hline\hline 1 & 0 & 0 & 0 & 0 & 0
\\\hline $\xi_{n-2}$& 1 & 0 & $2(3^{n-2}-1)$ & $2(3^{n-2}-1)$ & 0
\\\hline $\xi_{n-2}^2$& 1 & 0 & $2\cdot 2(3^{n-2}-1)$ & $2\cdot 2(3^{n-2}-1)$ & 0
\\\hline
\end{tabular}
\end{center}
As a Hopf algebra we have
$$ b_{10}^{-1}\Ext_D(k,k[\xi_{n-2}]/\xi_{n-2}^3) \isom R\tensor
D[\xi_{n-2}].  $$

\subsubsection*{Factor 3:
$\color{mgreen}k[\xi_{n-1},\xi_{n-2}^3]/(\xi_{n-1}^3,\xi_{n-2}^{27})$}
Similarly to \eqref{D-comodule_2}, for the third factor of $B_n$ we have a $D$-comodule decomposition
$$ k[\xi_{n-1},\xi_{n-2}^3]/(\xi_{n-1}^3,\xi_{n-2}^{27}) \isom \u{k\{ 1
\}}_{\isom k}\dsum \u{k\{ \xi_{n-1},\xi_{n-2}^3 \}}_{\isom M(1)}\dsum \u{k\{
\xi_{n-1}^2\xi_{n-2}^{21},\xi_{n-1}\xi_{n-2}^{24} \}}_{\isom M(1)}\dsum \u{k\{
\xi_{n-1}^2\xi_{n-2}^{24} \}}_{\isom k}\dsum F$$
where $F$ is a free $D$-comodule,
which gives the following $R$-module generators of
\\$b_{10}^{-1}\Ext_D^*(k,k[\xi_{n-1},\xi_{n-2}^3]/(\xi_{n-1}^3,\xi_{n-2}^{27})):$

\begin{center}\renewcommand{\arraystretch}{1.1}
\begin{tabular}{|c|c|c|c|c|c|}
\hline element & $s$ & $t$ & $u$ & $u''=u-6t$ & $\alpha$
\\\hline\hline 1 & 0 & 0 & 0 & 0 & 0
\\\hline $e_{n-1}\colonequals [\xi_1]\xi_{n-1} - [\xi_1^2]\xi_{n-2}^3$ & 1 &
1 & $2(3^{n-1}+1)$ & $2(3^{n-1}-2)$ & 3
\\\hline $y_{n-1}  \colonequals [\xi_1]\xi_{n-1}^2\xi_{n-2}^{21} +[\xi_1^2]\xi_{n-1}\xi_{n-2}^{24}$  & 1 & 1 & $2(3^{n+1}-21)$ & $2(3^{n+1}-24)$
& 27
\\\hline $z_{n-1}\colonequals \xi_{n-1}^2\xi_{n-2}^{24}$ & 1 & 0 & $2(3^{n+1}+3^{n-1}-26)$
& $2(3^{n+1}+3^{n-1}-26)$ & 30
\\\hline
\end{tabular}
\end{center}

\begin{lemma}\label{d_1(e_{n-1})}
$e_{n-1}$ is a permanent cycle in $E_r(k,B_n)$. In particular, $d_1(e_{n-1})=0$.
\end{lemma}
\begin{proof}
Use the filtration spectral sequence
interpretation of the MPASS described in the proof of Lemma \ref{e_2-f_20},
where $e_{n-1}$ has representative
$$ [\xi_1|\xi_{n-1}] - [\xi_1^2|\xi_{n-2}^3] $$
in $C_{P_n}(k,k)$. It is clear that this is a cycle in $C_{P_n}(k,k)$, hence a
permanent cycle in the spectral sequence.
\end{proof}

\subsubsection*{Factor 4: $\color{red}k[\xi_n,\xi_{n-1}^3]/(\xi_n^3,\xi_{n-1}^9)$}
There is a $D$-comodule decomposition
\begin{align*} k[\xi_n,\xi_{n-1}^3] & /(\xi_n^3,\xi_{n-1}^9)
\\ & \isom 
\u{k\{ 1 \}}_{\isom k} \dsum \u{k\{ \xi_n,\xi_{n-1}^3 \}}_{\isom M(1)}\dsum \u{k\{ \xi_n^2,\xi_{n-1}^3\xi_n,\xi_{n-1}^6
\}}_{\isom D}\dsum \u{k\{ \xi_{n-1}^3\xi_n^2, \xi_{n-1}^6\xi_n \}}_{\isom M(1)}\dsum
\u{k\{ \xi_{n-1}^6\xi_n^2 \}}_{\isom k}.\end{align*}
The non-free summands lead to $R$-module generators of
$b_{10}^{-1}\Ext_D^*(k,k[\xi_n,\xi_{n-1}^3]/(\xi_n^3,\xi_{n-1}^9))$ which have
representatives (in order):

\begin{center}\renewcommand{\arraystretch}{1.1}
\begin{tabular}{|c|c|c|c|c|c|}
\hline element & $s$ & $t$ & $u$ & $u''=u-6t$ & $\alpha$
\\\hline\hline 1 & 0 & 0 & 0 & 0 & 0
\\\hline $e_n\colonequals [\xi_1]\xi_n - [\xi_1^2]\xi_{n-1}^3$ & 1 & 1 &
$2(3^n+1)$ & $2(3^n-2)$ & 9
\\\hline $f_{n0} \colonequals [\xi_1]\xi_{n-1}^3\xi_n^2 - [\xi_1^2]\xi_{n-1}^6\xi_n$ & 1 & 1 & $2(3^{n+1}-3)$ & $2(3^{n+1}-6)$ & 27
\\\hline $c_n \colonequals \xi_{n-1}^6\xi_n^2$ & 1 & 0 & $2(3^{n+1}+3^n-8)$ &
$2(3^{n+1}+3^n-8)$ & 36
\\\hline
\end{tabular}
\end{center}
\begin{corollary}
There is an isomorphism of $R$-modules
\begin{align*}
b_{10}^{-1}\Ext_{D}(k,B_n) & \isom R\{ 1, e_2,f_{20},c_2 \}
\tensor R\{ 1, \xi_{n-2},\xi_{n-2}^2 \}
\\ & \hspace{30pt}\tensor R\{ 1, e_{n-1},y_{n-1},z_{n-1}
\}\tensor R\{ 1, e_n,f_{n,0},c_n \}.
\end{align*}
\end{corollary}
We have already computed part of the Hopf algebra structure on
$b_{10}^{-1}\Ext_D(k,B_n)=E_1^{1,*}(k,B_n)$ but
do not need to finish this; we just need one more piece of information.
\begin{lemma}\label{d_1(e_n)}
$e_n$ is primitive in $b_{10}^{-1}\Ext_D(k,B_n)$
\end{lemma}
\begin{proof}
Write $\psi(e_n)= \sum_i c[x_i|y_i]$,
where $c\in R$ and $x_i,y_i\in b_{10}^{-1}\Ext_D(k,B_n)$. As the cobar
differential preserves the grading $\alpha$ (see Proposition \ref{def-alpha}) and $\psi$ can be given in terms of the cobar
differential (see e.g. Remark \ref{MPASS-B_n-filtration-coincide}), $\psi$ also
preserves $\alpha$.
Since $\alpha(e_n)=9$, in order for
$d_1(e_n)$ to have $\alpha = 9$, we need $\alpha(x_i) + \alpha(y_i) = 9$.
Looking at $\alpha$ degrees in the above charts of $R$-module generators in
$b_{10}^{-1}\Ext_D(k,B_n)$, the only options are for $e_n\mid x_i$ or $y_i$, or
for $e_{n-1}^2 \mid x_i$ or $y_i$. But $e_{n-1}^2=0$ by Lemma \ref{e(x)e(y)},
and so the only option is for $e_n$ to be primitive.
\end{proof}
Combining Lemmas \ref{e_2-f_20}, \ref{d_1(e_{n-1})}, and \ref{d_1(e_n)} we have:
\begin{corollary}
In $b_{10}^{-1}\Ext_D(k,B_n)$, the elements $e_2$, $f_{20}$, $e_{n-1}$, and
$e_n$ are exterior generators in the Hopf algebra sense---they are primitive and
square to zero.
\end{corollary}

Now we have computed enough of $E_2(k,B_n)$ to show Proposition
\ref{section3-result}. If $b_{10}^{-4}h_{10}w_2^2w_{n-1}^3$ (which is in degree
$\alpha = 9$, $u'=2(3^n-8)$, and $u=2(3^n+1)$) is the target of a differential,
it must be a $d_r$ for $r\leq 4$ (since the target is in filtration 5), and the
source of that differential must have degree $\alpha = 9$, $u'=2(3^n-5)$, and
$u=2(3^n+1)$. Thus it suffices to prove Proposition
\ref{almost-section3-result}.

\begin{proposition}\label{almost-section3-result}
The only element in $E_2(k,B_n)$ with $s\leq 4$, $\alpha = 9$, $u' =
2(3^n-5)$, and $u=2(3^n+1)$ is $\pm w_n$.
\end{proposition}
\begin{proof}
There is a map $R\tensor E[e_2,f_{20},
e_{n-1},e_n]\tensor D[\xi_{n-2}]\to b_{10}^{-1}\Ext_D(k,B_n)$ that is an
isomorphism on degree $u''< 2(3^{n+1}-24)$ and induces a map on cobar complexes
$$ C^s_{R\tensor E[e_2,f_{20},e_{n-1},e_n]\tensor D[\xi_{n-2}]}(R,R) \to
C^s_{b_{10}^{-1}\Ext_D(k,B_n)}(R,R).$$
We claim the map of cobar complexes is an isomorphism in degree $u'' < -2 +
2(3^{n+1}-24) + 14(s-1)$. One can see this by noting that a minimal-degree
element in $C^s_{b_{10}^{-1}\Ext_D(k,B_n)}(R,R)$ not in the image is
$h_{10}[y_{n-1}|e_2|\dots|e_2]$, in degree $-2 + 2(3^{n+1}-24) + 14(s-1)$. (We use
$u''$ degree here because it is additive with respect to multiplication within
$b_{10}^{-1}\Ext_D(k,B_n)=E_1^{1,*}$, whereas $u'$ degree is additive with
respect to multiplication of cohomology classes in $H^*E_1 = E_2$.) Note that
the desired degrees $u'' = u'+6s = 2(3^n-5) + 6s$ fall into the region described
here for every $s$.

Now we look at the map induced on $\Ext$ in this region.
Since $d_r$ differentials increase $u''$ degree by $6(r-1)$ (they preserve $u$
and decrease $t$ by $r-1$) and increase $s$ by $r$, differentials originating in
the region $u''< -2+2(3^{n+1}-24)+14(s-1)$ stay in the region, but
there might be differentials originating outside the region hitting elements in
the region. Instead of showing that the map on $\Ext$ is an isomorphism in a
smaller region, note that this is already enough for our purposes: we want to
check that $\Ext_{b_{10}^{-1}\Ext_D(k,B_n)}(R,R)$ is zero in particular dimensions, and it
suffices to check that in $\Ext_{R\tensor E[e_2,f_{20},e_{n-1},e_n]\tensor
D[\xi_{n-2}]}(R,R)$.

We have $$\Ext_{R\tensor E[e_2,f_{20},e_{n-1},e_n]\tensor D[\xi_{n-2}]}(R,R)
\isom R[w_2,b_{20},b_{n-2,0},w_{n-1},w_n]\tensor E[h_{n-2,0}]$$
where $w_i=[e_i]$, $b_{20}=[f_{20}]$, and $\Ext_{D[\xi_{n-2}]}(R,R)= R\tensor
E[h_{n-2,0}]\tensor k[b_{n-2,0}]$. Degree information is as follows:
\begin{center}\renewcommand{\arraystretch}{1.1}
\begin{tabular}{|c|c|c|c|c|c|}
\hline element & $s$ & $t$ & $u$ & $u'$ & $\alpha$
\\\hline\hline $w_2$ & 1 & 1 & 20 & 8 & 0
\\\hline $b_{20}$ & 1 & 1 & 48 & 36 & 0
\\\hline $h_{n-2,0}$ & 1 & 0 & $2(3^{n-2}-1)$ & $2(3^{n-2}-1)$ & 0
\\\hline $b_{n-2,0}$ & 2 & 0 & $2(3^{n-1}-3)$ & $2(3^{n-1}-3)$ & 0
\\\hline $w_{n-1}$ & 1 & 1 & $2(3^{n-1}+1)$ & $2(3^{n-1}-5)$ & 3
\\\hline $w_n$ & 1 & 1 & $2(3^n+1)$ & $2(3^n-5)$ & 9
\\\hline\hline $h_{10}$ & 0 & 1 & 4 & $-2$ & 0
\\\hline $b_{10}$ & 0 & 2 & 12 & 0 & 0
\\\hline
\end{tabular}
\end{center}
Of course, $w_n$ has the right degree. Any other monomial with the right degree must be in
$R[w_2,b_{20},b_{n-2,0},w_{n-1}]\tensor E[h_{n-2,0}]$, and it is clear
from looking at $\alpha$ degree above
that it must have the form $w_{n-1}^3x$ (where $x\in 
R[w_2,b_{20},b_{n-2,0}]\tensor E[h_{n-2,0}]$).
Since $u'(w_{n-1}^3) = 2(3^n-15)$, we need $u'(x)=20$, which is not possible
using $w_2$ in degree 8, $b_{20}$ in degree 36, $h_{10}$ in degree $-2$ (where
$h_{10}^2=0$), and $h_{n-2,0}$ and $b_{n-2,0}$ in higher degree.

So the element must be $\pm b_{10}^Nw_n$, and by checking $u$ degree we see that
the power $N$ has to be zero.
\end{proof}

\subsection{Degree-counting in the ISS} \label{section:ISS-computation}
Recall that $b_{10}^{-4}h_{10}w_2^2w_{n-1}^3$ has $\alpha = 9$ and 
$u'=2(3^n-8)$; if it were a permanent cycle, it would converge to an element of
$b_{10}^{-1}\Ext_{P_n}^{a,b}(k,k)$ with stem $b-6a=2(3^n-8)$ (see Definition
\ref{def-stem}) and $\alpha=9$. The goal of this section is to prove:
\begin{proposition}\label{hw_2^2w_{n-1}^3_stem}
The sub-vector space of $b_{10}^{-1}\Ext^*_{P_n}(k,k)$ consisting of elements in
stem $2(3^n-8)$ and $\alpha=9$ is zero.
\end{proposition}
We will prove this using a (localized) Ivanovskii spectral sequence
(ISS) computing $b_{10}^{-1}\Ext_{P_n}(k,k)$. In our case, the ISS is
constructed by filtering the cobar complex for $P_n$ by powers of the
augmentation ideal. For example, $[\xi_n]$ is
in filtration 1, and in the Milnor diagonal
$$ d_{\text{cobar}}([\xi_n]) = [\xi_1| \xi_{n-1}^3] +
[\xi_2| \xi_{n-2}^9], $$
$[\xi_1|\xi_{n-1}^3]$ is in filtration 4 (since
$[\xi_1]$ is in filtration 1 and $[\xi_{n-1}^3]$ is in filtration 3), and
$[\xi_2|\xi_{n-2}^9]$ is in filtration 10.
In general, all of the multiplicative generators
$\xi_1,\xi_2,\xi_{n-2},\xi_{n-1},\xi_n$ are primitive in the associated graded,
i.e. they are in $\ker d_0$. To form the $b_{10}$-localized spectral sequence,
take the colimit of multiplication by $b_{10}$. In Section
\ref{section:ISS-convergence} we show that the (localized and un-localized) ISS
converges in our case.

So we have $E_0 \isom D[\xi_1,\xi_1^3, \xi_2,
\xi_{n-2},\xi_{n-2}^3,\xi_{n-2}^9,\xi_{n-1},\xi_{n-1}^3,\xi_n]$ and
$$ E_1^{ISS} = E[h_{1i},h_{20}, h_{n-2,j}, h_{n-1,i}, h_{n0}]_{i\in \{ 0,1 \}\atop
j\in \{ 0,1,2 \}}\tensor k[b_{10}^{\pm1}, b_{11}, b_{20},
b_{n-2,j},b_{n-1,i},b_{n,0}]_{i\in \{ 0,1 \}\atop j\in \{ 0,1,2 \}}. $$
Here $h_{ij} = [\xi_i^{3^j}]$ has filtration $3^j$ and $b_{ij}$ has filtration
$3^{j+1}$.
To help with the degree-counting argument in Proposition
\ref{hw_2^2w_{n-1}^3_stem}, here is a table of the degrees of the multiplicative
generators of the $E_1$ page.
\begin{center}
\renewcommand{\arraystretch}{1.1}
\begin{longtable}{|c|c|c|c|c|c|}
\hline element & $s$ & $t$ & $u$ & $u'=u-6t$ & $\alpha$
\\\hline\hline $h_{10}$ & 1 & 1 & 4 & $-2$ & 0
\\\hline $b_{10}$ & 3 & 2 & 12 & 0 & 0
%\\\hline $w_2$ & 2 & 2 & 20 & 8 & 0
\\\hline $h_{11}$ & 3 & 1 & 12 & 6 & 0
\\\hline $b_{11}$ & 9 & 2 & 36 & 24 & 0
\\\hline $h_{20}$ & 1 & 1 & 16 & 10 & 0
\\\hline $b_{20}$ & 3 & 2 & 48 & 36 & 0
\\\hline $h_{n-2,0}$ & 1 & 1 &
$2(3^{n-2}-1)$ & $2(3^{n-2}-4)$ & 0
\\\hline $b_{n-2,0}$ & 3 & 2 &
$2(3^{n-1}-3)$ & $2(3^{n-1}-9)$ & 0
\\\hline $h_{n-2,1}$ & 3 & 1 & $2(3^{n-1}-3)$ & $2(3^{n-1}-6)$ & 3
\\\hline $b_{n-2,1}$ & 9 & 2 & $2(3^n-9)$ & $2(3^n-15)$ & 9
\\\hline $h_{n-2,2}$ & 9 & 1 & $2(3^n-9)$ & $2(3^n-12)$ & 9
\\\hline $b_{n-2,2}$ & 27 & 2 & $2(3^{n+1}-27)$ & $2(3^{n+1}-33)$ & 27
\\\hline $h_{n-1,0}$ & 1 & 1 & $2(3^{n-1}-1)$ & $2(3^{n-1}-4)$ & 3
\\\hline $b_{n-1,0}$ & 3 & 2 & $2(3^n-3)$ & $2(3^n-9)$ & 9
\\\hline $h_{n-1,1}$ & 3 & 1 & $2(3^n-3)$ & $2(3^n-6)$ & 9
\\\hline $b_{n-1,1}$ & 9 & 2 & $2(3^{n+1}-9)$ & $2(3^{n+1}-15)$ & 27
\\\hline $h_{n,0}$ & 1 & 1 & $2(3^n-1)$ & $2(3^n-4)$ & 9
\\\hline $b_{n,0}$ & 3 & 2 & $2(3^{n+1}-3)$ & $2(3^{n+1}-9)$ & 27
\\\hline
\end{longtable}
\end{center}
\begin{proof}
[Proof of Proposition \ref{hw_2^2w_{n-1}^3_stem}]
The argument has two parts:
\begin{enumerate} 
\item show that (up to powers of $b_{10}$) the only generators in $E_1^{ISS}$
in degree $(u'=2(3^n-8),\alpha=9)$ are $h_{10}h_{20}h_{n-2,2}$ and
$h_{10}h_{11}h_{20}b_{n-2,1}$; 
\item show that those elements are targets of higher differentials in the
$b_{10}$-local ISS.
\end{enumerate}

From looking $\alpha$ degrees we see that no monomial in $E_1$ in degree
$(u'=2(3^n-8),\alpha=9)$
can be divisible by $b_{n-2,2}$, $b_{n-1,1}$, or $b_{n,0}$, and moreover by
looking at $u'$ degree we see it is not possible for $b_{n-1,0}$, $h_{n-1,1}$,
or $h_{n,0}$ to be a factor of such a monomial.
The only monomial of the right degree divisible by $h_{n-2,2}$ is
$b_{10}^Nh_{10}h_{20}h_{n-2,2}$. Any remaining elements of the right degree are in
$$E[h_{10},h_{11},h_{20}, h_{n-2,0},h_{n-2,1},h_{n-1,0}]\tensor k[b_{10}^{\pm
1}, b_{11}, b_{20}, b_{n-2,0},b_{n-2,1}].$$
Of these generators, only $h_{n-2,1}$, $h_{n-1,0}$, and $b_{n-2,1}$ have $\alpha>0$. Since
$h_{n-2,1}^2 =0 = h_{n-1,0}^2$, a monomial with $\alpha=9$ needs to be divisible
by $b_{n-2,1}$. If $u'(b_{n-2,1}x) = 2(3^n-8)$ then $u'(x) = 14$, and the
only possibility is $x = b_{10}^Nh_{10}h_{11}h_{20}$. (Here we are using the
assumption $n\geq 5$ to determine that $u'(h_{n-2,0})= 2(3^{n-2}-4)\geq 46$, and
the elements following it in the chart have greater degree).

This concludes part (1) of the argument; for (2) it suffices to show
\begin{align}
\label{ISS-d_9-1} d_9(h_{10}h_{20}b_{n-1,0}) & = h_{10}h_{20}h_{11}b_{n-2,1} - b_{10}h_{10}h_{20}h_{n-2,2}
\\\label{ISS-d_9-2} d_9(b_{10}h_{10}h_{n0}) & = -b_{10}h_{10}h_{20}h_{n-2,2}.
\end{align}
First, we claim that $h_{10}h_{20}$ is a permanent cycle; it is represented by
$[\xi_1|\xi_2] - [\xi_1^2|\xi_1^3] = w_2$, which we've seen is a permanent cycle
in the cobar complex.
The class $b_{n-1,0}$ has
cobar representative $[\xi_{n-1}|\xi_{n-1}^2] + [\xi_{n-1}^2|\xi_{n-1}]$ and
\begin{align*}
b_{n-1,0} & \equiv [\xi_{n-1}|\xi_{n-1}^2] + [\xi_{n-1}^2|\xi_{n-1}] -
[\xi_1\xi_{n-1}|\xi_{n-1}\xi_{n-2}^3]+ [\xi_1\xi_{n-1}^2|\xi_{n-2}^3]
\\ & \hspace{20pt}+[\xi_1^2\xi_{n-1}|\xi_{n-2}^6] +
[\xi_1^2|\xi_{n-1}\xi_{n-2}^6] + [\xi_1|\xi_{n-1}^2\xi_{n-2}^3] \in
(F^3/F^4)C_{P_n}^2(k,k).
\end{align*}
Computing the cobar differential on this class (and remembering that
$\xi_{n-3}^9=0$ in $P_n$), we see that $d_9(b_{n-1,0})=h_{11}b_{n-2,1} -
b_{10}h_{n-2,2}$. So
$$d_9(h_{10}h_{20}b_{n-1,0}) = h_{10}h_{20}\,d_9(b_{n-1,0})=
h_{10}h_{20}(h_{11}b_{n-1,1}- b_{10}h_{n-2,2}).$$
We have $h_{10}h_{n0}\equiv [\xi_1|\xi_n] -
[\xi_1^2|\xi_{n-1}^3] = w_n \in F^2/F^3$ and there is a cobar differential
$$ d_{\text{cobar}}([\xi_1|\xi_n]-[\xi_1^2|\xi_{n-1}^3]) = -
[\xi_1|\xi_2|\xi_{n-2}^9] + [\xi_1^2|\xi_1^3|\xi_{n-2}^9]. $$
This implies \eqref{ISS-d_9-2}.
(We did not check that $h_{10}h_{20}h_{11}b_{n-2,1}$ and
$h_{10}h_{20}b_{10}h_{n-2,2}$ survive to the $E_9$ page, because that is not
necessary: we only have to check that these elements die somehow in the spectral
sequence, and if they have already died before the $E_9$ page, then that is good
enough for this argument.)
\end{proof}

\section{Some results on higher differentials} \label{section:d_8}
In the case $r=4$, the following proposition gives an explicit way to compute
$d_8$ on any class, given our knowledge of $d_4$ from the previous section.
\begin{proposition}\label{prop:d_{r+4}}
Suppose $\bar{x}\in E_2$ satisfies $d_{r'}(\bar{x})=0$ for $r'<r$ and
$d_r(\bar{x}) = h_{10}\til{y}\in E_r$. Also suppose $\bar{y}$ is an $E_2$
representative for $\til{y}$ and $d_4(\bar{y}) = h_{10}\til{z}$. Then
$d_{r+4}(h_{10}\bar{x}) = b_{10}\til{z}$.
\end{proposition}
Note that the choice $\bar{y}$ does not matter, as two such choices differ (up
to $E_2$ class) by a boundary.

One is tempted to use Massey product arguments, e.g. try to apply the Massey
product differential and extension theorem \cite[4.5, 4.6]{may-matric} to
$\an{h_{10},h_{10},\til{z},h_{10}}$, but the following explicit argument avoids
Massey product technicalities.
\begin{lemma}\label{dx-form}
Suppose $0\neq \bar{x}\in E_2^{u'(x),s(x)}$ is not $h_{10}$-divisible, and
define $\bar{y}\in E_r^{u'(x)-4,s(x)+r}$ such that
$d_r(\bar{x})=h_{10}\bar{y}$ and $d_{r'}(\bar{x}) = 0$ for
$r'<r$. Furthermore, suppose $d_4(\bar{y})=h_{10}\bar{z}$. Then there is a cobar representative $x\in F^{s(x)}$ of
$b_{10}^N\bar{x}$ for some $N$, a cobar representative $y\in F^{s(x)+r}$ of
$b_{10}^N\bar{y}$, and a cobar representative $z\in F^{s(x)+r+4}$ of
$b_{10}^N\bar{z}$ such that
\begin{equation}\label{d(x)-goal} d(x)=[\xi_1|y]-[\xi_1^2|z]. \end{equation}
%Moreover, either $d(y)=0$ or $d(y) = [\xi_1|z]$ and $d(z) = 0$.
\end{lemma}
\begin{proof}
We prove this by induction on $u'$. The statement is trivially true for $u'<
-2$, since there are no elements of $E_2$ in those degrees.
So let $\bar{x}\in E_2$ with $u'(\bar{x})\geq-2$, and assume the inductive hypothesis.

By Proposition \ref{d_4-d_8}, $d_r(\bar{x})$ has the form $h_{10}\bar{y}$.
If $\bar{y}$ is not a permanent cycle, we abuse notation by letting $\bar{y}$
denote an $E_2$ representative.
By Proposition \ref{d_4-d_8}, there is a nontrivial differential
$d_R(\bar{y})=h_{10}\bar{z}$ for some $R\geq 4$ such that $d_{r'}(\bar{y})=0$ for
$r'<R$.
Since $u'(\bar{y}) = u'(\bar{x})-4$,
we may apply the inductive hypothesis to $\bar{y}$, obtaining a
cobar representative $y$ of $b_{10}^N\bar{y}$ for some $N$, a cobar
representative $z\in F^{s(x)+r+R}$ of $b_{10}^N\bar{z}$, and a cobar element
$w\in F^{s(x)+r+R+4}$ such
that \begin{equation}\label{d(y)} d(y) = [\xi_1|z] - [\xi_1^2|w]. \end{equation}
If $\bar{y}$ is a permanent cycle, \eqref{d(y)} holds with $z=0=w$.

Since $d_r(b_{10}^N\bar{x})=b_{10}^Nh_{10}\bar{y}$, 
there exists a cobar representative $x\in F^{s(x)}$ for $b_{10}^N\bar{x}\in
E_2^{s=s(x)}$ such that
$d(x) \equiv [\xi_1|y] \pmod {F^{s(x)+r+1}}$. In particular, we may write
\begin{equation}\label{x' def} d(x) = [\xi_1|y] - [\xi_1^2|z] + x' \end{equation}
where $x'\in F^{s(x)+r+1}$.
(Note that $[\xi_1^2|z]$ is also in higher filtration
than $y$, and this term is added because it simplifies the next calculation.)

\begin{claim}
We may choose $x$ and $x'$ such that $x'\in F^{s(x)+r+5}$.
\end{claim}
\begin{proof}[Proof of claim]
Applying $d$ to \eqref{d(y)}, we have 
$$ 0 = -[\xi_1|d(z)] + [\xi_1|\xi_1|w] + [\xi_1^2|d(w)]. $$
Equating terms starting with $\xi_1$, we obtain $d(z)=[\xi_1|w]$; equating terms
starting with $\xi_1^2$, we obtain $d(w)=0$.
Applying $d$ to \eqref{x' def}, we have
\begin{align*} 0  & = -[\xi_1|d(y)] + [\xi_1|\xi_1|z] + [\xi_1^2|d(z)] + d(x')
\\ & = [\xi_1|\xi_1^2|w] + [\xi_1^2|\xi_1|w] + d(x')
\end{align*}
so $d(x') = -b_{10}w\in F^{s(x)+r+R+4} \subseteq F^{s(x)+8}$.
So $x'$ represents an element of $E_2^{s=s(x)+r+1}$. Since $u'(x')=
u'(h_{10}y)$,
Lemma \ref{s-possibilities} implies that if $x'$ were nonzero in
$E_2^{s(x')}$, then $s(x')\equiv s(h_{10}y) = s(x)+r\pmod
9$. In particular, $x'$ is zero as an element of $E_2^{s(x)+r+1}$, so it must
have a representative in higher filtration. Repeating
this argument, we find $x'$ is zero as an element of $E_2^{s(x)+r+i}$ for
for $1\leq i\leq 5$. So we may write $x' + d(x_1)\in F^{s(x)+r+5}$, where
$x_1\in F^{s(x)+r}$. Thus by adjusting the representative $x$ by $x_1$, we may
assume $x'\in F^{s(x)+r+5}$.
\end{proof}

Then
\begin{align*}
d(b_{10}x) & = [\xi_1|\xi_1^2|\xi_1|y] + [\xi_1^2|\xi_1|\xi_1|y] -
[\xi_1|\xi_1^2|\xi_1^2|z] - [\xi_1^2|\xi_1|\xi_1^2|z] + b_{10}x'
\\d(b_{10}x - [\xi_1^2|\xi_1^2|y]) & = [\xi_1|\xi_1^2|\xi_1|y] -
[\xi_1|\xi_1^2|\xi_1^2|z] - [\xi_1^2|\xi_1|\xi_1^2|z] + 
b_{10}x' 
\\ & \hspace{30pt}+ [\xi_1|\xi_1|\xi_1^2|y]- [\xi_1^2|\xi_1^2|\xi_1|z] + [\xi_1^2|\xi_1^2|\xi_1^2|w]
\\ & =: [\xi_1|\til{y}] - [\xi_1^2|\til{z}]
\end{align*}
where
\begin{align*}
\til{y} & := b_{10}y - [\xi_1^2|\xi_1^2|z] + [\xi_1^2|x']
\\\til{z} & := b_{10}z - [\xi_1^2|\xi_1^2|w] - [\xi_1|x'].
\end{align*}
By our assumptions on the filtrations of all the elements involved,
$\til{y}\equiv b_{10}\pmod {F^{s(x)+r}}$ and $\til{z}\equiv b_{10}z\pmod
{F^{s(x)+r+5}}$, so $\til{y}$ is a representative of $b_{10}^{N+1}\bar{y}$ and
$\til{z}$ is a representative of $b_{10}^{N+1}\bar{z}$.
\end{proof}
%\fixme{Issue: $\til{y}$ could be a boundary. Maybe OK since the image of
%$\til{y}$ in $E_2$ = the image of $y$ in $E_2$? Can't explicitly rule out
%$\bar{y}=0$ since maybe $\bar{x}=0$. The condition ``if $y\in F^{s_y}$ then
%$z\in F^{s_y+1}$'' isn't capturing the fact that we don't want $y+$ a correction
%term to have higher filtration than $z$. Maybe there should be two cases: either
%$\bar{y}=0$ and $\bar{x}$ is a permanent cycle, or $d_r(\bar{x})=h_{10}\bar{y}$
%for $0\neq \bar{y}\in E_2^{s=s_y}$, and $y$ is a representative of $\bar{y}$.}

\begin{proof}[Proof of Proposition \ref{prop:d_{r+4}}]
Use Lemma \ref{dx-form} to write
\begin{equation}\label{d(x)} d(x)=[\xi_1|y]-[\xi_1^2|z] \end{equation}
where $x$ is a cobar representative for $b_{10}^N\bar{x}$, $y$ is a cobar
representative for $b_{10}^N\bar{y}$, and $z$ is a cobar representative for
$b_{10}^N\til{z}$.
Applying $d$ to \eqref{d(x)},
$$ 0 = -[\xi_1|d(y)] + [\xi_1|\xi_1|z] - [\xi_1^2|d(z)] .$$
Equating terms whose first component is $\xi_1$, we have $d(y)=[\xi_1|z]$; equating
terms whose first component is $\xi_1^2$, we have $d(z) = 0$.
Then $[\xi_1|x]-[\xi_1^2|y]$ is a
representative for $h_{10}\bar{x}$, and we have
\begin{align*}
d([\xi_1|x]-[\xi_1^2|y]) & = [\xi_1|\xi_1^2|z] + [\xi_1^2|\xi_1|z] =
b_{10}z.\qedhere
\end{align*}
Thus, in the $b_{10}$-localized spectral sequence,
$d_{r+4}(b_{10}^Nh_{10}\bar{x}) = b_{10}^N \til{z}$ implies
$d_{r+4}(h_{10}\bar{x})=b_{10}\til{z}$.
\end{proof}

\begin{random}{Conjecture}\label{conj:collapse}
The $K(\xi_1)$-based MPASS collapses at $E_9$.
\end{random}
Using computer calculations, we verified the conjecture for stems $\leq 600$.
However, it is not possible to rule out higher differentials based only on
degree.

\begin{proposition}
Assuming Conjecture \ref{conj:collapse}, we have
$$ b_{10}^{-1}\Ext_P(k,k) \isom b_{10}^{-1}\Ext_D(k,k[\til{w}_2,\til{w}_3,\dots]) $$
where $\til{w}_n = b_{10}^{-1} w_n$ and the $D$-coaction on the $E_2$ page is
given by $\psi(\til{w}_n) =1\tensor \til{w}_n+
\xi_1\tensor h_{10} \til{w}_2^2\til{w}_{n-1}^3$ for $n\geq 3$.
\end{proposition}

\begin{proof}
Let $\til{W} = k[\til{w}_2,\til{w}_3,\dots]$.
We have $d_4(\til{w}_n) = h_{10}\til{w}_2^2\til{w}_{n-1}^3$.
By Proposition \ref{d_4-d_8} and Conjecture \ref{conj:collapse}, the $E_\iy$
page of the MPASS is obtained by taking the cohomology of $E_2$ by $d_4$ and
$d_8$; more precisely, we have
\begin{align*}
E_\iy & \isom \ker(d_4|_{W_+})/\im(d_8|_{W_+}) \dsum
\ker(d_8|_{W_-})/\im(d_4|_{W_-}).
\end{align*}
If we let $\partial(x) = {1\over h_{10}}d_4(x)$, then Proposition
\ref{prop:d_{r+4}} says that $b_{10}\partial^2(x) = d_8(h_{10}x)$.
%\begin{align*}
%E_\iy & \isom k[b_{10}^{\pm}]\tensor
%\big(\ker(\partial|_{\til{W}})/\im(\partial^2|_{\til{W}}) \dsum h_{10}
%\ker(\partial^2|_{\til{W}})/\im(\partial|_{\til{W}})\big).
%\end{align*}
Thus we may write down an isomorphism $f$ of chain complexes
$$ \xymatrix{
\dots\ar[r] & \til{W}\ar[r]^\partial\ar[d]^-{f^{2n}} &
\til{W}\ar[r]^-{\partial^2}\ar[d]^-{f^{2n+1}} &
\til{W}\ar[r]^-\partial\ar[d]^-{f^{2n+2}} & \dots
\\\dots\ar[r] & \til{W}\{ b_{10}^n \}\ar[r]^-{d_4} & \til{W}\{ h_{10}b_{10}^n
\}\ar[r]^-{d_8} & \til{W}\{ b_{10}^{n+1} \} & 
}$$
By Lemma \ref{M-resolution}, the cohomology of the top complex is
$b_{10}^{-1}\Ext_D(k,\til{W})$, and we have argued below that the cohomology of
the bottom complex is $E_\iy$. Thus we have an isomorphism of vector spaces
$b_{10}^{-1}\Ext_P(k,k)\isom b_{10}^{-1}\Ext_D(k,\til{W})$.

It remains to show that this is an isomorphism of $R$-modules. We will just
check that the induced map $f_*$ on cohomology respects $h_{10}$-multiplication.
If $\omega = [x]\in \til{W}^{2n}$ is a cycle, then $h_{10}\omega$ is represented
by $[x]\in \til{W}^{2n+1}$. If $\nu=[y]\in \til{W}^{2n+1}$ is a cycle, then
$h_{10}\nu$ is represented by $[\partial y]\in \til{W}^{2n}$. So
$f_*^{2n+1}(h_{10}\omega) = [h_{10}b_{10}^nx] = h_{10}[b_{10}^nx] =
h_{10}f_*^{2n}(\omega)$. For the other case, we need to show that
$f_*^{2n+2}(h_{10}\nu)=[b_{10}^{n+1}(\partial y)]$ can be represented as
$h_{10}\cdot [h_{10}b_{10}^n y]=h_{10}f_*^{2n+1}(\nu)$. This corresponds to a
hidden multiplication in the MPASS. From the commutativity of the diagram we
have $d_4([b_{10}^ny])=[h_{10}b_{10}^n\partial y]=h_{10}[b_{10}^n\partial y]$.
The desired relation $h_{10}[h_{10}b_{10}^n y] = [b_{10}^{n+1}(\partial y)]$
follows from Lemma \ref{hidden-mult}.
\end{proof}

\begin{lemma} \label{hidden-mult}
Suppose $d_r(\bar{x})=h_{10}\bar{y}$ where $\bar{x}\in W_+$ and
$d_{r'}(\bar{x})=0$ for $r'<r$. Then there is a hidden multiplication
$h_{10}\cdot (h_{10}\bar{x}) = -b_{10}\bar{y}$.
\end{lemma}
This is closely related to the Massey product shuffle $h_{10} (h_{10}\bar{x}) =
h_{10}\an{h_{10}, h_{10}, \bar{y}} = \an{h_{10},h_{10},h_{10}}\bar{y}$, though
the following explicit argument avoids Massey product technicalities.
\begin{proof}
Use Lemma \ref{dx-form} to find a representative $x$ such that
$d(x)=[\xi_1|y]-[\xi_1^2|z]$ where $y$ is a representative for $\bar{y}$ and $z$
is a representative for $\bar{z}$ such that $d_4(\bar{y})=h_{10}\bar{z}$.
We use $[\xi_1|x]-[\xi_1^2|y]$ as a representative for $h_{10}\bar{x}$. Then
$h_{10}\cdot (h_{10}\bar{x})$ is represented by
$[\xi_1|\xi_1|x]-[\xi_1|\xi_1^2|y]$. Since $d([\xi_1^2|x]) = -[\xi_1|\xi_1|x] -
[\xi_1^2|\xi_1|y] + [\xi_1^2|\xi_1^2|z]$, we have
\begin{align*}
[\xi_1|\xi_1|x]-[\xi_1|\xi_1^2|y] + d(\xi_1^2|x)  & = -b_{10}y +
[\xi_1^2|\xi_1^2|z]\equiv -b_{10}\bar{y}. \qedhere
\end{align*}
\end{proof}

\section{Localized cohomology of a large quotient of $P$}\label{section:D_{1,iy}}
In this section we will prove Theorem \ref{D_iy}, a complete calculation of
$b_{10}$-local cohomology of a small $P$-comodule. Using the change of rings
theorem, this is equivalent to the following.
\begin{theorem}\label{D_{1,iy}}
Let $D_{1,\iy} = k[\xi_1,\xi_2,\dots]/(\xi_1^3)$. Then
$$b_{10}^{-1}\Ext^*_{D_{1,\iy}}(k,k) \isom E[h_{10},h_{20}]\tensor P[b_{10}^{\pm
1},b_{20},w_3,w_4,\dots]. $$
In particular, one can write
$$ b_{10}^{-1}\Ext^*_{D_{1,\iy}} (k,k) \isom b_{10}^{-1}\Ext^*_D(k,k[h_{20},
b_{20},w_3,w_4,\dots]/(h_{20}^2) $$
where all the generators $h_{20},b_{20},w_n$ are $D$-primitive.
\end{theorem}
Though $D_{1,\iy}$ seems reasonably close to $P$ in size,
the computation of its $b_{10}$-local cohomology is much simpler. In
particular, attempting to apply the methods in this section (especially the
explicit construction in Lemma \ref{theta-on-C}) to computing
$b_{10}^{-1}\Ext_P^*(k,k)$ quickly becomes intractable.

The strategy is to explicitly construct a map from the cobar complex
$C_{D_{1,\iy}}(k,k)$ to another complex which is designed to have the right
cohomology, and then show the map is a quasi-isomorphism. Note that the cobar
complex is a dga under the concatenation product, so every element is a product
of elements in degree 1. Thus if our target complex is a dga, it suffices to
construct a map out of $C^1_{D_{1,\iy}}(k,k) = \bar{D_{1,\iy}}$, and then extend the
map to all of $C_{D_{1,\iy}}^*(k,k)$ by multiplicativity. In order to ensure the
resulting map is a map of complexes, there is a criterion that the map on degree
1 needs to satisfy:
\begin{proposition}\label{twisting-cobar-map}
Let $\Gamma$ be a Hopf algebra over $k$, $Q^*$ be a dga with augmentation $k\to
Q^*$, and $\theta:\bar{\Gamma}\to Q^1$ be a $k$-linear map such that
\begin{equation}\label{twisting-condition}
d_Q(\theta(x)) = \sum \theta(x')\theta(x'')  \end{equation} for all $x\in
\bar{\Gamma}$, where $\sum x'\tensor x''$ is the reduced diagonal
$\bar{\Delta}(x)$. Then there is a map of dga's $f: C_\Gamma^*(k,k)\to Q^*$ sending
$[a_1|\dots|a_n]$ to $\prod \theta(a_i)$.
\end{proposition}
\begin{proof}
We just need to check that $f$ commutes with the differential; that is, we have
to check the following diagram commutes:
$$ \xymatrix{
C_\Gamma^n(k,k)\ar[r]^-f\ar[d]_-{d_{cobar}} & Q^n\ar[d]^-{d_Q}
\\C_\Gamma^{n+1}(k,k)\ar[r]^-f & Q^{n+1}
}$$
For $n=1$, this is precisely what the condition \eqref{twisting-condition}
guarantees. Commutativity for $n>1$ follows from the Leibniz rule. The map on
$n=0$ is the augmentation.
\end{proof}
\begin{remark}
This is an example of the more general construction of \emph{twisting cochains};
see \cite[\S II.1]{husemoller-moore-stasheff}. A morphism $\theta$ satisfying
\eqref{twisting-condition} will be called a \emph{twisting morphism}.
\end{remark}
The target of our desired twisting morphism will be the complex
$b_{10}^{-1}\til{U}^*\tensor W'$, where
\begin{itemize} 
\item $W' = k[w_3,w_4,\dots]$, with $u(w_n)=2(3^n-1)$, 
is in homological degree zero with zero differential, and
\item $\til{U}^* := UL^*(\xi_1)\tensor UL^*(\xi_2) \subset C_{D[\xi_1,\xi_2]}^*(k,k)$
where the sub-dga $UL^*(x)\subset C_{D[x]}^*(k,k)$ is defined below.
\end{itemize}
\begin{definition}
Given a height-3 truncated polynomial algebra $D[x]$,
let $UL^*(x)$ be the sub-dga of $C_{D[x]}^*(k,k)$ multiplicatively generated by
the elements $\alpha = [x]$,  $\beta= [x^2]$, and $\gamma= [x|x^2]+ [x^2|x]$.
This inherits from $C_{D[x]}^*(k,k)$ the differentials $d(\alpha)=0$,
$d(\beta)=-\alpha^2$, and $d(\gamma)=0$, along with the relations
$\alpha\beta+\beta\alpha=\gamma$, $\alpha^3=0$, and $\beta^2=0$.
\end{definition}
\begin{remark}
This is (up to signs) the $p=3$ case of a construction due to Moore: let $UL^*$
be the dga which has multiplicative generators
$a_1,\dots,a_{p-1}$ in degree 1 and $t_2,\dots,t_p$ in degree 2 with
$d(a_i)=t_i$, subject to
\begin{align*}
a_1^2 & =t_2 & a_i^2 & = 0\text{ for }i\neq 1 &  a_1^p & =0 & a_ia_j & =-a_ja_i\text{ for $i,j\neq 1$} 
\\a_ja_1 & = -a_1a_j + t_{j+1} & a_it_j & = t_ja_i & t_it_j & =t_jt_i.
\end{align*}
This is a dga quasi-isomorphic to, and much smaller than, $C_{k[x]/x^p}(k,k)$.
It also has the nice property that $t_p$ (which, in the case $x = \xi_1$,
represents $b_{10}$) is central.
\end{remark}

\begin{random}{Notation}
Denote the generators of $UL^*(\xi_1)$ by $a_1=[\xi_1]$, $a_2=[\xi_1^2]$, and
$b_{10}= [\xi_1|\xi_1^2]+[\xi_1^2|\xi_1]$, and the generators of $UL^*(\xi_2)$
by $q_1=[\xi_2]$, $q_2=[\xi_2^2]$, and $b_{20} =
[\xi_2|\xi_2^2]+[\xi_2^2|\xi_2]$.
(This definition of $b_{10}$ and $b_{20}$ does, of course, match up with the
image of $b_{10}$ and $b_{20}$ along $\Ext^*_P(k,k)\to
\Ext^*_{D[\xi_1,\xi_2](k,k)}$, and even $\Ext^*_P(k,k)\to
\Ext^*_{D_{1,\iy}}(k,k)$.) Note that
$$H^*(\til{U}) = H^*(C_{D[\xi_1,\xi_2]}(k,k))= E[h_{10},h_{20}]\tensor
P[b_{10},b_{20}].$$ 
So our target complex $b_{10}^{-1}\til{U}\tensor W'$ has
cohomology \begin{equation}\label{coh-UW'} H^*(b_{10}^{-1}\til{U}\tensor W')=H^*(b_{10}^{-1}\til{U})\tensor W' = E[h_{10},h_{20}]\tensor
P[b_{10}^{\pm 1},b_{20}]\tensor W'.\end{equation}
\end{random}

\subsection{Defining $\theta: \bar{D_{1,\iy}}\to b_{10}^{-1}\til{U}\tensor
W'$}\label{section:def-theta}
The definition of the map $\theta: \bar{D_{1,\iy}}\to b_{10}^{-1}\til{U}^*\tensor W'$ is
quite ad hoc, and will be done in several stages. The map will arise as a
composition $D_{1,\iy}\to D' \to \til{U}^*\tensor W'\to
b_{10}^{-1}\til{U}^*\tensor W'$, where the first map is the natural surjection to
$$D'\colonequals k[\xi_1,\xi_2,\dots]/(\xi_1^3,\xi_2^9,\xi_3^9,\dots)$$ and
the last map is the natural localization map; the main goal is to construct a
map $D'\to \til{U}^*\tensor W'$ satisfying the twisting morphism condition, and
we begin by constructing a map out of a slightly smaller coalgebra.

%The first step (and most of the work) is to construct a map out of
%$$C = k[\xi_1,\xi_2^3,\xi_3,\xi_4,\dots]/(\xi_1^3,\xi_2^9,\xi_3^9,\dots).$$ 
%This is a sub-coalgebra (though not sub-Hopf algebra) of $D'$.
%Since $\xi_2$ is primitive, we have an isomorphism of coalgebras $D' \isom
%C\tensor D[\xi_2]$. We will first define a twisting morphism $\theta$ out of $C$,
%and then in Lemma \ref{theta-on-D'} we will extend it to $D'$, eventually
%pre-composing with the quotient map $D_{1,\iy}\to D'$ to get a twisting morphism
%out of $D_{1,\iy}$.
\begin{lemma}\label{theta-on-C}
Let 
$$C = k[\xi_1,\xi_2^3,\xi_3,\xi_4,\dots]/(\xi_1^3,\xi_2^9,\xi_3^9,\dots).$$
There is a twisting morphism $\theta: \bar{C}\to UL^1(\xi_1)\tensor W'$.
\end{lemma}
\begin{proof}
For $n,m,k\geq 3$, make the following definitions:
\begin{align*}
\theta(\xi_1) & = a_1
\\\theta(\xi_1^2) &= a_2
\\\theta(\xi_{n-1}^3) & = -a_1w_n
\\\theta(\xi_n) & = a_2w_n
\\\theta(\xi_1\xi_{n-1}^3) & = -a_2w_n
\\\theta(\xi_1\xi_n) & = 0
\\\theta(\xi_{n-1}^3\xi_{m-1}^3) & = a_2w_nw_m
\\\theta(\xi_n\xi_{m-1}^3) & = 0
\\\theta(\xi_n\xi_m) & = 0
\\\theta(\xi_1^2\xi_{n-1}^3) & = 0
\\\theta(\xi_1\xi_{n-1}^3\xi_{m-1}^3) & = 0
\\\theta(\xi_{n-1}^3\xi_{m-1}^3\xi_{k-1}^3) & = 0
\end{align*}
It is a straightforward computation with the cobar differential to check that each of these does not violate the twisting morphism condition
\begin{equation} d(\theta(x))=\sum \theta(x')\cdot \theta(x'')\end{equation}
where $\bar{\Delta}(x)=\sum x'\tensor x''$.
(Note that, in $C$, we have $\bar{\Delta}(\xi_{n-1}^3)=0$ and
$\bar{\Delta}(\xi_n)=\xi_1|\xi_{n-1}^3$.)

Now it suffices to prove the following.
\begin{claim}
Defining $\theta(X)=0$ for all monomials $X$ except the ones listed above
defines a twisting morphism.
\end{claim}
Define a (non-multiplicative) grading $\rho$ on $C$ where
\begin{align*}
\rho(1) & =0 & \rho(\xi_1) & =1 & \rho(\xi_1^2) & =2 &
\rho(\xi_{n-1}^3) & = 1 & \rho(\xi_{n-1}^6) & =2 & \rho(\xi_n) & = 2 
& \rho(\xi_n^2) & = 4
\end{align*}
for $n\geq 3$, and $\rho(\prod_i \xi_i^{a_i + 3b_i}) = \sum
\rho(\xi_i^{a_i}) + \rho(\xi_i^{3b_i})$
(where $a_i,b_i\in \{ 0,1,2 \}$).
The reason for considering this grading is the following:
\begin{claim}\label{alpha-grading-claim}
Writing $\Delta(x)=\sum x'\tensor x''$, we have $\rho(x')+\rho(x'')\leq
\rho(x)$.
\end{claim}
\begin{proof}[Proof of Claim \ref{alpha-grading-claim}]
If $X = \prod \xi_i^{a_i+3b_i}$ for $a_i,b_i\in \{ 0,1,2 \}$, consider the
collection $\mathscr{T}_X = \{ \xi_i^{a_i} \st a_i\neq 0 \}\union \{ \xi_i^{3b_i} \st
b_i\neq 0 \}$. Use induction on $n:=\#\mathscr{T}_X$. If $n=1$, then it suffices to check
explicitly the Milnor diagonal of each of the terms $\{
\xi_1,\xi_1^2,\xi_{i-1}^3,\xi_{i-1}^6,\xi_i,\xi_i^2 \}$. (In fact, we find
$\rho(x)=\rho(x') + \rho(x'')$ for each of these terms.)

For general monomials $a,b$, we have \begin{equation}\label{alpha(ab)}\rho(ab)\leq
\rho(a)+\rho(b). \end{equation}
By definition, if $x$ and $y$ are products of non-overlapping subsets of
$\mathscr{T}_X$, then
\begin{equation}\label{alpha(xy)}\rho(xy)=\rho(x)+\rho(y). \end{equation}
Write $X = xy$ where $x\in \mathscr{T}_X$ and $y$ is a product of terms in
$\mathscr{T}_X$
(different from $x$). Since $\Delta(xy) = \sum x'y'|x''y''$ it suffices to prove
$\rho(x'y') + \rho(x''y'') \leq \rho(xy)$. We have
\begin{align*}
\rho(x'y') + \rho(x''y'') & \leq \rho(x') + \rho(y') + \rho(x'') +
\rho(y'')
\\ & \leq \rho(x) + \rho(y)
\\ & = \rho(xy)
\end{align*}
where the first inequality is by \eqref{alpha(ab)}, the second inequality is by
the inductive hypothesis, and the last equality is by \eqref{alpha(xy)}.
\end{proof}
So the monomials in $C$ with degree 1 are $\xi_1$
and $\xi_{n-1}^3$ for $n\geq 3$, the monomials with $\rho$-degree 2 are 
$\xi_1^2$, $\xi_n$, $\xi_{n-1}^3\xi_{m-1}^3$, and $\xi_1\xi_{n-1}^3$ for
$n,m\geq 3$, and the monomials with degree 3 are $\xi_1^2\xi_{n-1}^3$,
$\xi_1\xi_{n-1}^3\xi_{m-1}^3$, $\xi_{n-1}^3\xi_{m-1}^3\xi_{k-1}^3$,
$\xi_1\xi_n$, and $\xi_{n-1}^3\xi_m$ for $n,m\geq 3$. Notice that $\theta$ has
already been defined for these monomials above. So it remains to show that
$\theta$ can be defined consistently for monomials with $\rho\geq 4$. In
particular, we will show using induction on $\rho$ degree that we can define
$\theta(x)=0$ if $\rho(x)\geq 3$ while preserving the twisting morphism
condition \eqref{twisting-condition}.

Since we have already checked above that we can define $\theta(x)=0$ on the
monomials $x$ with $\rho(x)=3$, let $\rho(x)=n>3$ and assume inductively
that we have already defined $\theta(y)=0$ if $3\leq \rho(y)\leq n-1$. Any
monomial $y$ with $\rho(y)=0$ is in $k$ (and hence $\theta(y)=0$), so we can
assume that $\rho(x')< \rho(x)$ and $\rho(x'') < \rho(x)$. So by the
inductive hypothesis we have $\sum \theta(x')\cdot \theta(x'') = 0$, and so we
can set $\theta(x)=0$ without violating \eqref{twisting-condition}.
\end{proof}

\begin{lemma} \label{theta-on-D'}
One may extend $\theta$ constructed in Lemma \ref{theta-on-C} to a twisting
morphism $\bar{D'}\to \til{U}^1\tensor W'$ by defining:
\begin{align*}
\theta(\xi_2) & = q_1
\\\theta(\xi_2^2) & = q_2
\\\theta(\xi_2x) & = 0 \hspace{20pt}\text{ for $x\in \bar{C}$}
\\\theta(\xi_2^2x) & = 0 \hspace{20pt}\text{ for $x\in \bar{C}$}
\end{align*}
where $\bar{C}$ is the cokernel of the unit map $k\to C$.
\end{lemma}
\begin{proof}
Note that $\xi_2$ is primitive in $D'$, and $C$ is a sub-coalgebra of $D'$, so
we need to define $\theta$ on $\xi_2C$ and $\xi_2^2C$. It is straightforward to
check that $\theta(\xi_2)=q_1$ and $\theta(\xi_2^2)=q_2$ is consistent with
\eqref{twisting-condition}.

If $x = \xi_2y$ for $y\in \bar{C}$ then every $y',y''$ in $\Delta y$ is in $C$,
and
\begin{align*}
\sum\theta(x')\cdot \theta(x'')  & = \sum \big(\theta(\xi_2 y')\cdot \theta(y'') +
\theta(y')\cdot \theta(\xi_2 y'')\big)
\\ & = \theta(\xi_2)\theta(y) + \theta(y)\theta(\xi_2)
+ \sum_{y',y''\notin k} \big(\theta(\xi_2 y')\cdot \theta(y'') +
\theta(y')\cdot \theta(\xi_2 y'')\big)
\\ & = q_1\theta(y) + \theta(y)q_1
+ \sum_{y',y''\notin k} \big(\theta(\xi_2 y')\cdot \theta(y'') +
\theta(y')\cdot \theta(\xi_2 y'')\big).
\end{align*}
Since $\theta(y)\in UL^1(\xi_1)\tensor W'$ and $q_1$ anti-commutes with the
generators $a_1$ and $a_2$ of $UL^1(\xi_1)$, we have
$q_1\theta(y)+\theta(y)q_1=0$. Thus defining $\theta(\xi_2y)=0$ does not violate
\eqref{twisting-condition}.

Similarly, if $x = \xi_2^2y$ for $y\in \bar{C}$, then
\begin{align*}
\sum \theta(x')\cdot \theta(x'') & = \sum \big(\theta(\xi_2^2y')\cdot \theta(y'') +
2\theta(\xi_2 y')\cdot \theta(\xi_2y'') + \theta(y')\cdot
\theta(\xi_2^2y'')\big)
\\ & = \theta(\xi_2^2)\theta(y) +2\theta(\xi_2)\theta(\xi_2y) +
2\theta(\xi_2y)\theta(\xi_2) + \theta(y)\theta(\xi_2^2) 
\\ & \hspace{40pt}+ \sum_{y',y''\notin
k} \big(\theta(\xi_2^2y')\cdot \theta(y'') +
2\theta(\xi_2 y')\cdot \theta(\xi_2y'') + \theta(y')\cdot
\theta(\xi_2^2y'')\big)
\\ & = \theta(\xi_2^2)\theta(y) + \theta(y)\theta(\xi_2^2) +\sum_{y',y''\notin k}
\big(\theta(\xi_2^2y')\theta(y'') + \theta(y')\theta(\xi_2^2y'')\big)
\end{align*}
where in the third equality we use the fact that
$0=\theta(\xi_2y)=\theta(\xi_2y')=\theta(\xi_2y'')$ (for $y',y''\notin k$).
Again, $\theta(\xi_2^2)\theta(y) + \theta(y)\theta(\xi_2^2)=q_2\theta(y)+ 
\theta(y)q_2$ which is zero since $\theta(y)$ is in $UL^1(\xi_1)\tensor W'$ and $q_2$
anti-commutes with the generators $a_1$ and $a_2$ of $UL^1(\xi_1)$. So it is consistent
with \eqref{twisting-condition} to define $\theta(\xi_2^2y)=0$.
\end{proof}

Now precompose with the surjection $q:D_{1,\iy}\to D'$ to obtain a twisting
morphism $$\theta:D_{1,\iy}\to D'\to \til{U}^1\tensor W'.$$ This remains a
twisting morphism because it is a coalgebra map---in particular, $q$ commutes
with the coproduct---and so $d(\theta(q(x))) = \sum \theta(q(x)')\theta(q(x)'')
= \sum \theta(q(x'))\theta(q(x''))$. So by Proposition \ref{twisting-cobar-map}
we get an induced map \begin{equation}\label{theta-prime}\theta': C_{D_{1,\iy}}^*(k,k)\to
\til{U}^*\tensor W'\end{equation}
by extending $\theta$ multiplicatively using the concatenation product on the
cobar complex.

\subsection{Showing $\theta$ is a quasi-isomorphism via spectral sequence
comparison}\label{section:theta-qis}
Our goal is to show that
the map $$\theta':C_{D_{1,\iy}}^*(k,k)\to \til{U}^*\tensor W'$$ induces an
isomorphism in cohomology after inverting $b_{10}$.

To prove this, we define filtrations on $C^*_{D_{1,\iy}}(k,k)$ and
on $\til{U}^*\tensor W'$ in a way that makes $\theta'$
a filtration-preserving map; this induces a map of filtration spectral
sequences. We compute the $E_2$ pages of both sides and show that $\theta'$
induces an isomorphism of $E_2$ pages, hence an isomorphism of $E_\iy$ pages.

Let $B_{1,\iy}:= k[\xi_2,\xi_3,\dots] = D_{1,\iy}\cotensor_D k$. Define a
decreasing filtration on $C_{D_{1,\iy}}^*(k,k)$ where $[a_1|\dots|a_n]$ is in
$F^sC_{D_{1,\iy}}^*(k,k)$ if at least $s$ of the $a_i$'s are in
$\ker(D_{1,\iy}\to D)=\bar{B}_{1,\iy}D_{1,\iy}$. Define a decreasing filtration
on $\til{U}^*\tensor W'$ by the following multiplicative grading:
\begin{itemize} 
\item $|a_1|= |a_2| = |b_{10}|=0$
\item $|q_1| = |q_2|=1$
\item $|b_{20}| = 2$
\item $|w_n| = 1$.
\end{itemize}
Looking at the definition of $\theta$ in Lemma \ref{theta-on-C} and Lemma
\ref{theta-on-D'}, it is clear that $\theta$ is filtration-preserving, and hence
so is $\theta'$.

For the same reasons that the $b_{10}^{-1}B$-based MPASS coincides at $E_1$ with
the filtration spectral sequence mentioned in Section \ref{section:main-tool},
the $b_{10}^{-1}B_{1,\iy}$-based MPASS for computing
$b_{10}^{-1}\Ext^*_{D_{1,\iy}}(k,k)$ coincides with the $b_{10}$-localized
version of the filtration spectral sequence on $C_{D_{1,\iy}}^*(k,k)$ defined
above. Our next goal is to calculate the $E_2$ page of (the $b_{10}$-localized
version of) the filtration spectral sequence on $C_{D_{1,\iy}}^*(k,k)$, and
using this correspondence we may instead calculate the MPASS $E_2$ term
\begin{equation}\label{E_2_MPASS} E_2^{s,*} =
b_{10}^{-1}\Ext^s_{b_{10}^{-1}\Ext^*_D(k,B_{1,\iy})}(b_{10}^{-1}\Ext^*_D(k,k),\
b_{10}^{-1} \Ext^*_D(k,k)).\end{equation}
So we need to compute $b_{10}^{-1}\Ext^*_D(k,B_{1,\iy})$ and its coalgebra
structure. The correspondence of spectral sequences further gives that
\begin{equation}\label{MPASS-equals-filtration-sseq}E_1^{1,*}=b_{10}^{-1}\Ext_D^*(k,\bar{B}_{1,\iy}) \isom b_{10}^{-1}H^*(F^1/F^2 C_{D_{1,\iy}}^*(k,k))\end{equation}
and the reduced diagonal on $b_{10}^{-1}\Ext_D^*(k,B_{1,\iy})$ coincides with
$d_1$ in the filtration spectral sequence.

\begin{proposition}
As coalgebras, we have
$$ b_{10}^{-1}\Ext^*_D(k,B_{1,\iy}) \isom b_{10}^{-1}E[e_3,e_4,\dots] \tensor
D[\xi_2]$$
i.e. $e_n$ and $\xi_2$ are primitive and $\bar{\Delta}(\xi_2^2)=2\xi_2\tensor
\xi_2$.
\end{proposition}
\begin{proof}
The first task is to determine the $D$-comodule structure on $B_{1,\iy}$. Let
$\psi$ denote the $D$-coaction induced by the $D$-coaction on $P$, and
$\partial:B_{1,\iy}\to B_{1,\iy}$ denote the operator defined by $\psi(x) =
1\tensor x + \xi_1\tensor \partial x - \xi_1^2\tensor \partial^2 x$ (see
Definition \ref{def-partial}). For example, $\partial(\xi_n)=\xi_{n-1}^3$,
$\partial(\xi_{n-1}^3)=0$, and $\partial$ satisfies the Leibniz rule.

We have a coalgebra isomorphism $B_{1,\iy} \isom D[\xi_2]\tensor
k[\xi_2^3,\xi_3,\xi_4,\dots]$. Since 1, $\xi_2$, and $\xi_2^2$ are all
primitive, $D[\xi_2]$ splits as $D$-comodule into three trivial $D$-comodules,
generated by 1, $\xi_2$, and $\xi_2^2$ respectively. So it suffices to determine
the $D$-comodule structure of $k[\xi_2^3,\xi_3,\xi_4,\dots]$.

As part of the determination of the structure of $b_{10}^{-1}\Ext^*_D(k,B)$ in
Section \ref{section:B}, we showed that there is a $D$-comodule decomposition
$$ B \isom
\dsums_\attop{\let\scriptstyle\scriptstyle\xi_{n_1}\dots
\xi_{n_d}\\n_i\geq 2\text{ distinct}}
T(\ann{\xi_{n_1}\dots \xi_{n_d}\cc 1})\ \dsum\ F $$
where $F$ is a free $D$-comodule and $T(\an{\xi_{n_1}\dots \xi_{n_d}\cc 1})$ is
generated as a vector space by monomials of the
form $\partial^{\epsilon_1}(\xi_{n_1})\dots \partial^{\epsilon_d}(\xi_{n_d})$ for
$\epsilon_i\in \{ 0,1 \}$. I claim the surjection $f:B\to
k[\xi_2^3,\xi_3,\xi_4,\dots]$ takes $F$ to another free summand: this
map preserves the direct sum decomposition into summands of the form $D$,
$M(1)$, and $k$, and the image of a free summand $D$ must be either 0 or another
free summand (just as there are no $D$-\emph{module} maps $k=k[x]/(x) \to D$ or
$M(1)=k[x]/(x^2)\to D$, there are no $D$-comodule maps $D\to k$ or $D\to M(1)$).

Furthermore, I claim that $f$ acts as zero on summands
$T(\an{\xi_{n_1}\dots \xi_{n_d}\cc 1})$ where some $n_i=2$, and is the identity
otherwise. In the first case, every basis element
$\partial^{\epsilon_i}(\xi_2)\prod_{j\neq i}\partial^{\epsilon_j}(\xi_{n_j})$ in
$T(\an{\xi_{n_1}\dots \xi_{n_d}\cc 1})$ has the form $\xi_2 \prod_{j\neq
i}\partial^{\epsilon_j}(\xi_{n_j})\in \xi_1^3\cdot k[\xi_2^3,\xi_3,\xi_4,\dots]$
or $\xi_1^3 \prod_{j\neq i}\partial^{\epsilon_j}(\xi_{n_j})\in \xi_1^3\cdot
k[\xi_2^3,\xi_3,\xi_4,\dots]$, and these are sent to zero under $f$. If
instead $n_i>2$ for every $i$, then every term
$\partial^{\epsilon_1}(\xi_{n_1})\dots \partial^{\epsilon_d}(\xi_{n_d})$ is in
$k[\xi_2^3,\xi_3,\xi_4,\dots]$ and so $f$ acts as the identity. So we
have shown that there is a $D$-comodule isomorphism
$$ B_{1,\iy} = \Big(\!\!\!\!\dsums_\attop{\xi_{n_1}\dots \xi_{n_d}\\n_i\geq 3\text{
distinct}} T(\an{\xi_{n_1}\dots \xi_{n_d}\cc 1})\ \dsum\ F'\Big)\ \ \tensor\ \
(k_1\dsum k_{\xi_2} \dsum k_{\xi_2^2}) $$
where $F'$ is a free $D$-comodule. So we have
\begin{align*}
b_{10}^{-1}\Ext^*_D(k,B_{1,\iy})  & \isom \dsums_\attop{\xi_{n_1}\dots \xi_{n_d}\\n_i\geq 3\text{
distinct}} b_{10}^{-1}\Ext_D^*\big(k,T(\an{\xi_{n_1}\dots \xi_{n_d}\cc 1})\tensor
k\{ 1,\xi_2,\xi_2^2 \}\big) 
\\ & \isom \dsums_\attop{\xi_{n_1}\dots \xi_{n_d}\\n_i\geq 3\text{
distinct}} b_{10}^{-1}\Ext_D^*\big(k,T(\an{\xi_{n_1}\dots \xi_{n_d}\cc 1})\big)\tensor
k\{ 1,\xi_2,\xi_2^2 \}.
\end{align*}

By Proposition \ref{e(x)e(y)},
$b_{10}^{-1}\Ext_D^d(k,T(\an{\xi_{n_1}\dots \xi_{n_d}\cc 1}))$ is
generated by $e_{n_1}\dots e_{n_d}$, where $$e_n = [\xi_1]\xi_n - [\xi_1^2]\xi_{n-1}^3\in b_{10}^{-1}\Ext_D^1(k,T(\an{\xi_n\cc 1}))$$
is primitive. The map $B\to B_{1,\iy}$ gives rise to a map of
MPASS's, and in particular a map $b_{10}^{-1}\Ext^*_D(k,B)\to
b_{10}^{-1}\Ext^*_D(k,B_{1,\iy})$ of Hopf algebras over $b_{10}^{-1}\Ext^*_D(k,k)$ sending $e_n\mapsto e_n$
for $n\geq 3$, and $e_2\mapsto h_{10}\cdot \xi_2$. In particular, we have
\begin{equation}\label{E_1(B_{1,infty})} b_{10}^{-1}\Ext^*_D(k,B_{1,\iy}) \isom E[h_{10},e_3,e_4,\dots]\tensor
P[b_{10}^{\pm 1}]\tensor k\{ 1,\xi_2,\xi_2^2 \} \end{equation} and $e_n\in
b_{10}^{-1}\Ext^*_D(k,B_{1,\iy})$ is primitive. To find the coproduct on the
elements $\xi_2$ and $\xi_2^2$, use \eqref{MPASS-equals-filtration-sseq}, in
particular the fact that the (reduced) Hopf algebra diagonal corresponds to
$d_1$ in the filtration spectral sequence. In particular, $\xi_2\in
b_{10}^{-1}\Ext_D^*(k,B_{1,\iy})$ corresponds to the element $[\xi_2]\in
F^1/F^2 C^1_{D_{1,\iy}}(k,k)$, and
we have $\bar{d}_{\text{cobar}}([\xi_2]) = [\xi_1|\xi_1^3]$
which is zero in $C_{D_{1,\iy}}^*(k,k)$, so $\xi_2$ is primitive. Similarly, the
cobar differential on $C_{D_{1,\iy}}^*(k,k)$ shows $\bar{\Delta}(\xi_2^2) =
2\xi_2\tensor \xi_2$. Thus the tensor factor $k\{ 1,\xi_2,\xi_2^2 \}$ is, as a
coalgebra, a truncated polynomial algebra. This finishes the determination of
the coalgebra structure of $b_{10}^{-1}\Ext_D^*(k,B_{1,\iy})$ in
\eqref{E_1(B_{1,infty})}.
\end{proof}

The $E_2$ page \eqref{E_2_MPASS} of the MPASS is the cohomology of the Hopf algebroid
$$(b_{10}^{-1}\Ext^*_D(k,k),\ b_{10}^{-1}\Ext^*_D(k,B_{1,\iy})) = (E[h_{10}]\tensor
P[b_{10}^{\pm 1}],\ E[h_{10},e_3,e_4,\dots]\tensor P[b_{10}^{\pm 1}]\tensor
D[\xi_2])$$
so we have:
\begin{corollary}\label{D_{1,iy}-E_iy-page}
The MPASS $E_2$ page is:
$$ E_2^{**}\isom E[h_{10},h_{20}] \tensor P[b_{10}^{\pm
1},b_{20},w_3,w_4,\dots].$$
\end{corollary}

\begin{proposition}\label{theta-prime-iso-E_2}
The map $\theta'$ induces an isomorphism of $E_2$ pages after inverting $b_{10}$.
\end{proposition}
\begin{proof}
We first show that the $E_2$ pages of the filtration spectral sequences on
$C_P^*(k,k)$ and $\til{U}^*\tensor W'$ are abstractly isomorphic after inverting
$b_{10}$. By the
machinery of Section \ref{section:main-tool}, it suffices to calculate the
$E_2$ page for $\til{U}^*\tensor W'$ and check that it coincides with the $E_2$ page of the
MPASS from Corollary \ref{D_{1,iy}-E_iy-page}. Then we show that the
map $\theta'$ induces this isomorphism.

In the associated graded, there is a differential $d_0(a_2)= -a_1^2$, but the
corresponding differential on $q_2$ is a $d_1$.
So the filtration spectral sequence $\EU_r$ computing
$H^*(b_{10}^{-1}\til{U}^*\tensor W')$ has $E_0$ page
$$ \EU_0\isom b_{10}^{-1}UL^*(\xi_1)\tensor UL^*(\xi_2)\tensor W' $$
with differential $d_0(u_1 \tensor u_2 \tensor w) = d(u_1)\tensor u_2 \tensor w$. So $$\EU_1 \isom
H^*(b_{10}^{-1}UL^*(\xi_1))\tensor UL^*(\xi_2)\tensor W' \isom E[h_{10}]\tensor
P[b_{10}^{\pm 1}]\tensor UL^*(\xi_2)\tensor W'$$
and the only remaining differential is generated by $d_1(q_2) = -q_1^2$, so
$$ \EU_2 \isom E[h_{10}] \tensor P[b_{10}^\pm]\tensor H^*(UL^*(\xi_2))\tensor W'
= E[h_{10},h_{20}]\tensor P[b_{10}^{\pm 1},b_{20}]\tensor W'.$$
Then $E_r\isom E_2$ for $r\geq 2$.

To show that $\theta'$ is an isomorphism, it suffices to show that
$\theta'(h_{10})=h_{10}$, $\theta'(b_{10})=b_{10}$, $\theta'(h_{20})=h_{20}$,
$\theta'(b_{20})=b_{20}$, and $\theta'(w_n)=b_{10}w_n$ for $n\geq 3$. We use
the fact that $\theta'$ extends $\theta$ multiplicatively using the
concatenation product in the cobar complex. So $\theta'([a_1|\dots|a_n]) = \prod
\theta(a_i)$, and we have:
\begin{align*}
\theta'(h_{10}) & = \theta'([\xi_1]) = \theta(\xi_1) = a_1
\\\theta'(b_{10}) & = \theta'([\xi_1|\xi_1^2]+[\xi_1^2|\xi_1]) =
\theta(\xi_1)\theta(\xi_1^2) + \theta(\xi_1^2)\theta(\xi_1) = a_1a_2 + a_2a_1 =
b_{10} 
\\\theta'(h_{20}) & = \theta'([\xi_2]) =\theta(\xi_2) = q_1
\\\theta'(b_{20}) & = \theta'([\xi_2|\xi_2^2]+[\xi_2^2|\xi_2]) =
\theta(\xi_2)\theta(\xi_2^2) + \theta(\xi_2^2)\theta(\xi_2) = q_1q_2 + q_2q_1 =
b_{20}
\\\theta'(w_n) & = \theta'([\xi_1|\xi_n] - [\xi_1^2|\xi_{n-1}^3])
= a_1a_2w_n + a_2a_1w_n = b_{10}w_n.
\qedhere
\end{align*}
\end{proof}
\begin{proof}[Proof of Theorem \ref{D_{1,iy}}]
In Section \ref{section:def-theta} we constructed a map
$\theta':C^*_{D_{1,\iy}}(k,k)\to \til{U}^*\tensor W'$ which is
filtration-preserving, where $C^*_{D_{1,\iy}}(k,k)$ has the filtration
associated to the MPASS and $\til{U}^*\tensor W'$ has the filtration constructed
in Section \ref{section:theta-qis}. 
By Proposition \ref{theta-prime-iso-E_2}, $\theta'$ induces an isomorphism of
spectral sequences after inverting $b_{10}$, and so it induces an isomorphism in
cohomology. Thus $$ b_{10}^{-1} \Ext_{D_{1,\iy}}(k,k) = b_{10}^{-1}
H^*(C_{D_{1,\iy}}(k,k)) \isom b_{10}^{-1} H^*(b_{10}^{-1} \til{U}^*\tensor
W').$$
The result follows from \eqref{coh-UW'}.
\end{proof}

\section*{Appendix A: Convergence of localized spectral sequences}
\renewcommand{\thesection}{A}
\setcounter{subsection}{0}
\setcounter{theorem}{0}
\setcounter{equation}{0}
In this appendix, we study the convergence of two $b_{10}$-localized spectral
sequences, the $b_{10}$-localized MPASS (the main subject of this paper) and the
$b_{10}$-localized ISS (introduced in Section \ref{section:d_4}). In each case,
the non-localized spectral sequences converges for straightforward reasons.

In general, there are two possible ways in which a localization of a convergent
spectral sequence can fail to converge.
\begin{enumerate} 
\item There could be a $b_{10}$-tower $x$ in $E_\iy$ that does not appear in
$b_{10}^{-1}E_\iy$ because it is broken into a series of $b_{10}$-torsion towers
connected by hidden multiplications.
\item There could be a $b_{10}$-tower $x$ in $b_{10}^{-1}E_\iy$ that is not a
permanent cycle in $E_\iy$ because in the non-localized spectral sequence it
supports a series of increasing-length differentials to $b_{10}$-torsion
elements (so these differentials would be zero in $b_{10}^{-1}E_r$).
\end{enumerate}
(The reverse of (2), where a sequence of torsion elements
supports a differential that hits a $b_{10}$-tower, cannot happen: if
$d_r(x)=y$ and $b_{10}^nx = 0$ in $E_r$, then $0=d_r(b_{10}^nx) =
b_{10}^nd_r(x)=b_{10}^ny$.)

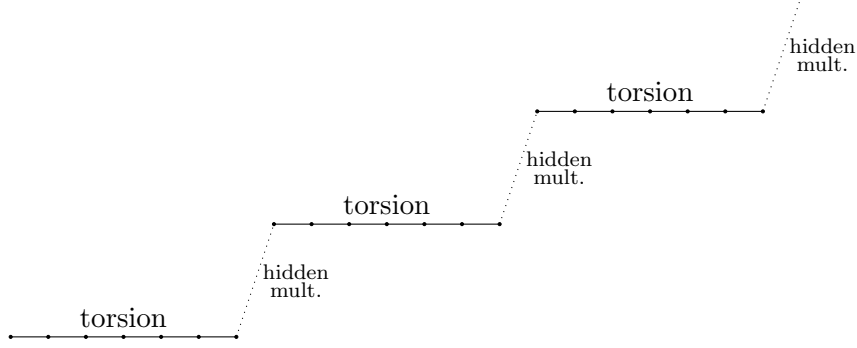
\begin{figure}[H]
\begin{center}
\begin{tikzpicture}
\foreach \i in {0,...,6}{
    \draw[fill] (\i*0.50,0) circle [radius=0.02];
}
\foreach \i in {0,...,6}{
    \draw[fill] (\i*0.50+7*0.50,1.5) circle [radius=0.02];
}
\foreach \i in {0,...,6}{
    \draw[fill] (\i*0.50+14*0.50,3) circle [radius=0.02];
}
\draw (0,0) -- (6*0.50, 0);
\draw (7*0.50,1.5) -- (13*0.50, 1.5);
\draw (14*0.50,3) -- (20*0.50, 3);
\draw[dotted] (6*0.50, 0) -- (7*0.50,1.5);
\draw[dotted] (13*0.50, 1.5) -- (14*0.50,3);
\draw[dotted] (20*0.50, 3) -- (21*0.50,4.5);
\node[align=center,above] at (1.5, 0) {torsion};
\node[align=center,above] at (5, 1.5) {torsion};
\node[align=center,above] at (8.5, 3) {torsion};
\node[align=center,right] at (3.5,0.75) {\hspace{-8pt}$\substack{\text{hidden}\\\text{mult.}}$};
\node[align=center,right] at (7,2.25) {\hspace{-8pt}$\substack{\text{hidden}\\\text{mult.}}$};
\node[align=center,right] at (10.5,3.75) {\hspace{-8pt}$\substack{\text{hidden}\\\text{mult.}}$};
\end{tikzpicture}
\end{center}
\caption{Illustration of (1): this represents a $b_{10}$-tower in ``homotopy''}
\end{figure}

\subsection{Convergence of the $K(\xi_1)$-based MPASS}\label{section:convergence}
In this section we prove convergence of the $B_\Gamma$-based MPASS of Theorem
\ref{K(xi_1)-MPASS-M} in the case that
$\Gamma$ is a quotient of $P$ (in fact, the only property of $\Gamma$
that is used is that $u(x)\geq u(\xi_1^3)$ for $u\in \Gamma$).
The convergence argument will only rely on the form of the $E_1$ page.

\begin{proposition}\label{convergence}
For any non-negatively graded $\Gamma$-comodule $M$, the $b_{10}$-localized
$K(\xi_1)$-based MPASS
\begin{equation}\label{MPASS-general} E_1^{s,*} = b_{10}^{-1}\Ext_D(k,
\bar{B}_\Gamma^{\tensor s}\tensor M) \implies
b_{10}^{-1}\Ext_\Gamma(k,M)\end{equation}
converges.
\end{proposition}
The proof is a slight modification of \cite[Proposition 4.4.1, Proposition
4.2.6]{palmieri-book}.

Recall our grading convention: $x \in E_1^{s,t,u}$ is an element in
$\Ext^t_\Gamma(k,B_\Gamma\tensor\bar{B}_\Gamma^{\tensor s})$ with internal degree $u$.

\begin{lemma} \label{2016-10-27_0-count}
Let $M$ be a bounded-below graded $D$-comodule and suppose $u_M = \min\{ u(x)\st
x\in M \}$. If $x\in \Ext^*_D(k,M)$ is a nonzero element of degree $(s,t,u)$ and
$x\neq 0$, then $u\geq u_M+6t-2$.
\end{lemma}
\begin{proof}
It suffices to check the cases $M = k$, $M(1) = k[\xi_1]/\xi_1^2$, and $D$.
In the case $M =
k$, we have $\Ext_D^*(k,k\{ y \})=E[h_{10}]\tensor k[b_{10}]\tensor k\{ y \}$.
In the case $M = M(1)$, write $M = k\{ y,\partial y \}$; then $\Ext_D^*(k,M) =
k[b_{10}]\tensor k\{ \partial y, e(y) \}$ where $e(y) = [\xi_1]y -
[\xi_1^2](\partial y)$. In the case $M = D$, $\Ext_D^0(k,D)\isom k$ is
concentrated in homological degree zero.
In each of these cases, we verify the desired statement, using
the fact that $b_{10}\in E_1^{0,2,12}$ and $h_{10}\in E_1^{0,1,4}$. 
\end{proof}
%\begin{proof} 
%First we check the cases when $M \isom k,M(1)$, or $D$.

%\emph{Case 1: $M \isom k$.} Let $y$ be the generator of $M$, in degree $(t,u)=
%(0,u(y))$. We have $\Ext^*_D(k,k\{ y \})= E[h_{10}]\tensor P[b_{10}]\tensor k\{
%y \}$ where $h_{10}$ is in degree $(t,u) = (1,|\xi_1|) = (1,4)$ and $b_{10}$ is
%in degree $(t,u)=(2,12)$. The minimum degree element is $y$, so $u_M = u(y)$.
%Every element has the form $h_{10}b_{10}^ny$ or $b_{10}y$ for $n\geq 0$, and
%both of these satisfy $u\geq u_M+6t-2$.

%\emph{Case 2: $M\isom M(1)$.} Write $M = k\{ y,\partial y \}$, where $\partial
%y$ is in degree $(0, u(\partial y))$ and $\partial y$ is in degree $(0,
%u(\partial y)+4)$.
%By Lemma \ref{e(x)e(y)}(1), $\Ext_D^*(k,M) = \F_p[b_{10}]\tensor k\{ \partial
%y,e(y) \}$ where $e(y)$ is in degree $(t,u) = (1,u(\partial y)+8)$.
%The minimum degree element is $\partial y$, and all the elements satisfy $u\geq
%u_M+6t$.

%\emph{Case 3: $M\isom D$.} Here, $\Ext^0_D(k,M)\isom k$ has degree
%$(t,u)=(0,u_M)$ and $\Ext^t_D(k,M)=0$ for $t>0$.

%In general, a homogeneous element $x\in M$ is a sum $\sum x_i$ for $x_i\in M_i$
%where $M_i$ is a summand of the above type, and by definition, $u_{M_i}\geq
%u_M$. So $u(x) = u(x_i) \geq u_{M_i}+6t-2\geq u_M+6t-2$.
%\end{proof}
\begin{proposition} \label{vanishing-plane-E_1}
There is a vanishing plane in the $E_1$ page of \eqref{MPASS-general}:
$E^{s,t,u}_1=0$ if $u<12s+6t-2$.
\end{proposition}
\begin{proof} 
Recall $E_1^{s,t,*} = \Ext^t_\Gamma(k,\Gamma\cotensor_D(\bar{B_\Gamma}^{\tensor s}\tensor M))\isom
\Ext_D(k,\bar{B}_\Gamma^{\tensor s}\tensor M)$. Since $\Gamma$ is a quotient of
$P$, if $x\in \bar{B}_\Gamma$ is nonzero then $u(x)\leq u(\xi_1^3) = 12$. Therefore a nonzero element $x\in
\bar{B}_\Gamma^{\tensor s}\tensor M$ has $u\geq 12s$. By Lemma \ref{2016-10-27_0-count}, if
$x\in E_1^{s,t,u}$ has degree $(s,t,u)$, then $u\geq 12s+6t-2$.
\end{proof}
\begin{corollary} \label{2016-10-27_0_d_r}
The differential $d_r:E_r^{s,t,u}\to E_r^{s+r,t-r+1,u}$ is zero if $r> {1\over 6}(u-12s-6t-4)$.
\end{corollary}
\begin{proof} 
Given $x\in E_r^{s,t,u}$, $d_r(x)\in E_r^{s',t',u'}=E_r^{s+r,t-r+1,u}$ will be zero
because of the vanishing plane if $12s'+6t'-2-u'>0$. But
$$ 12s'+6t'-2-u' = 12(s+r)+6(t-r+1)-2-u = (12s+6t+4-u) +6r $$
which is $>0$ for $r$ as indicated.
\end{proof}
\begin{corollary}
There is a vanishing line in $\Ext_\Gamma^*(k,M)$: if $x\in \Ext_\Gamma^{t',u}(k,M)$ and
$u-6t'+2 < 0$ then $x=0$.
\end{corollary}
\begin{proof}
Permanent cycles in $E_1^{s,t,u}$ converge to elements in $\Ext_\Gamma^{s+t,u}(k,M)$. 
Any such $x$ would then be represented by a permanent cycle in $E_1^{s,t,u}$
with $u-6(s+t)+2<0\leq 6s$ (since Adams filtrations are non-negative), which
falls in the vanishing region of Proposition \ref{vanishing-plane-E_1}.
\end{proof}
Note that $b_{10}\in \Ext_\Gamma^{2,12}(k,M)$ acts parallel to this vanishing line.

\begin{proof} [Proof of Proposition \ref{convergence}]
Convergence of the non-localized MPASS follows from a
general result by Palmieri \cite[Proposition 1.4.3]{palmieri-book}. 

For convergence problem (1), suppose $x$ has degree $(s_x,t_x,u_x)$. If there were no multiplicative
extensions, then $b_{10}^ix$ would have degree $(s_x,t_x+2i,u_x+12i)$. But
multiplicative extensions cause it to have the expected internal degree $u$ and
stem $s+t$, but higher $s$. That is, $b_{10}^ix$ has degree
$(s_x+n_i,t_x+2i-n_i,u_x+12i)$ for some $n_i>0$, and because this scenario
involves the existence of infinitely many multiplicative extensions, the
sequence $(n_i)_i$ is increasing and unbounded above. This causes us to run
afoul of the vanishing plane (Proposition \ref{vanishing-plane-E_1}) for
sufficiently large $i$:
\begin{align*}
12s+6t-2-u  & = 12(s_x + n_i) + 6(t_x+2i-n_i) -2 -
(u_x+12i) \\ & = 12s_x + 6t_x-2-u_x + 6n_i
\end{align*}
which is $>0$ for $i\gg 0$.

For convergence problem (2), the scenario is, more precisely, as follows: we have a
$b_{10}$-periodic element $x\in \Ext^*_\Gamma(k,k)$, and a sequence of differentials
$d_{r_i}(b_{10}^ix) = y_i\neq 0$, where every $y_i$ is $b_{10}$-torsion. The
sequence $(r_i)_i$ must be increasing and bounded above: if $b_{10}^{n_i}y_i=0$
then $d_{r_i}(b_{10}^{n_i}x)=b_{10}^{n_i}y_i=0$, and so if $b_{10}^{n_i}x$ is to
support a differential $d_{r_{n_i}}$, we must have $r_{n_i}>r_i$.
% Why can't all the $d_r$'s be the same? If $d_r(x)=y_1$ and $b_{10}^ny_1=0$,
% then $d_r(b_{10}^nx) = b_{10}d_r(b_{10}^{n-1}x)=b_{10}\cdot b_{10}^{n-1}y=0$
% and we need a higher differential to kill $b_{10}^nx$. Repeat process\dots
Note that the condition on $r$ in Corollary \ref{2016-10-27_0_d_r} is the same for
all $b_{10}^ix$. So some of the $r_i$'s will be greater than this bound,
contradicting the assumption that $d_{r_i}(b_{10}^ix)\neq 0$.
\end{proof}

\subsection{Convergence of the $b_{10}$-local ISS}\label{section:ISS-convergence}
In this section, we consider the $b_{10}$-local ISS computing $b_{10}^{-1}
\Ext_{P_n}^*(k,k)$. As discussed in Section \ref{section:ISS-computation}, this
is obtained by $b_{10}$-localizing a filtration spectral sequence on the cobar
complex for $P_n$, where the filtration is defined by taking powers of the
augmentation ideal. Let $E_r^{ISS}$ denote the $E_r$ page of the non-localized
ISS and $b_{10}^{-1}E_r^{ISS}$ denote the $E_r$ page of the localized ISS.

\begin{lemma}\label{ISS-vanishing-line}
There is a slope ${1\over 4}$ vanishing line in $E_1^{ISS}$ in $(u,s)$
coordinates. That is, if $x\in E_1^{ISS}$ has $s(x) > {1\over 4}u(x)$ then
$x=0$.
\end{lemma}
\begin{proof}
In Section \ref{section:ISS-computation} we computed the $E_1$ page:
$$ E_1^{ISS} =  \tensors_{(i,j)\in I}E[h_{ij}]\tensor k[b_{ij}] $$
where $I = \{ (1,0),(1,1),(2,0),(n-2,0),(n-2,1),
(n-2,2),(n-1,0),(n-1,1),(n,0) \}$.
These generators occur in the following degrees:
\begin{center}\renewcommand{\arraystretch}{1.1}
\begin{tabular}{|c|c|c|c|}
\hline element & $u$ & $s$ & $u/s$
\\\hline\hline $h_{ij}$ & $2(3^i-1)3^j$ & $3^j$ & $2(3^i-1)$
\\\hline $b_{ij}$ & $2(3^i-1)3^{j+1}$ & $3^{j+1}$ & $2(3^i-1)$
\\\hline
\end{tabular}
\end{center}
So we have ${u\over s}\geq 2(3^1-1)=4$, which proves the lemma. Note that
$b_{10}$, in degree $(u=12, s=3)$, acts parallel to the vanishing line.
\end{proof}
Here is a picture:
\begin{center}
\begin{tikzpicture} [xscale=0.25]
\draw[fill,teal!10] (0,0) -- (12.8,3.2) -- (17.5,3.2) -- (17.5,0) -- (0,0);
\draw[<->] (0,3.2) -- (0,0) -- (17.8,0);
\draw [help lines] (0,0) grid (17.5, 3.2);
\node[below right] at (17.5,0) {$u$};
\node[left] at (0,3.2) {$s$};
\draw[fill] (4,1) ellipse(5pt and 5/4pt);
\node[below right] at (4,1) {$h_{10}$};
\draw[fill] (12,3) ellipse(5pt and 5/4pt);
\node[below right] at (12,3) {$b_{10}$};
\draw[fill] (16,1) ellipse(5pt and 5/4pt);
\node[below] at (16,1) {$h_{20}$};
\draw[dotted,thick] (0,0) -- (13,13/4);
\node[left] at (0,1) {\tiny 1};
\node[left] at (0,2) {\tiny 2};
\node[below] at (4,0) {\tiny 4};
\node[below] at (8,0) {\tiny 8};
\node[below] at (12,0) {\tiny 12};
\node[below] at (16,0) {\tiny 16};
\node[below left] at (0,0) {\tiny 0};
\end{tikzpicture}
\end{center}
Differentials are vertical: $d_r$ takes elements in degree $(u,s)$ to degree $(u,s+r)$.
\begin{proposition}
The $b_{10}$-localized ISS converges to $b_{10}^{-1}\Ext_{P_n}(k,k)$.
\end{proposition}
\begin{proof}
The non-localized ISS converges because it is based on a decreasing filtration
of the cobar complex that clearly satisfies both $\intss_s F^sC_{P_n}(k,k) = \{
0 \}$ and $\unions_s F^sC_{P_n}(k,k) = C_{P_n}(k,k)$.

The two convergence problems are illustrated below:
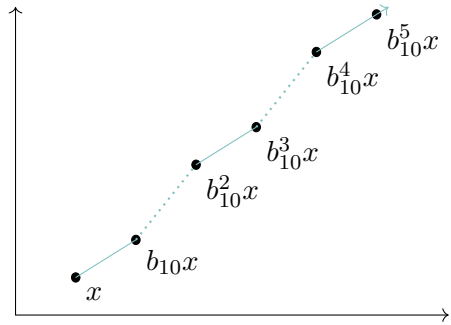
\begin{figure}[H]
\begin{center}
\begin{tikzpicture} [xscale=0.8,yscale=0.5]
\draw[<->] (0,8.2) -- (0,0) -- (7.2,0);
\draw[fill] (1,1) ellipse (2pt and 7/2pt);
\draw[fill] (2,2) ellipse (2pt and 7/2pt);
\draw[fill] (3,4) ellipse (2pt and 7/2pt);
\draw[fill] (4,5) ellipse (2pt and 7/2pt);
\draw[fill] (5,7) ellipse (2pt and 7/2pt);
\draw[fill] (6,8) ellipse (2pt and 7/2pt);
\node[below right] at (1,1) {$x$};
\node[below right] at (2,2) {$b_{10}x$};
\node[below right] at (3,4) {$b_{10}^2x$};
\node[below right] at (4,5) {$b_{10}^3x$};
\node[below right] at (5,7) {$b_{10}^4x$};
\node[below right] at (6,8) {$b_{10}^5x$};
\draw[teal!50] (1,1) -- (2,2);
\draw[teal!50,dotted,thick] (2,2) -- (3,4);
\draw[teal!50] (3,4) -- (4,5);
\draw[teal!50,dotted,thick] (4,5) -- (5,7);
\draw[teal!50,->] (5,7) -- (6.2,8.2);
\end{tikzpicture}
\end{center}
\caption{Convergence problem (1) for the ISS}
\end{figure}

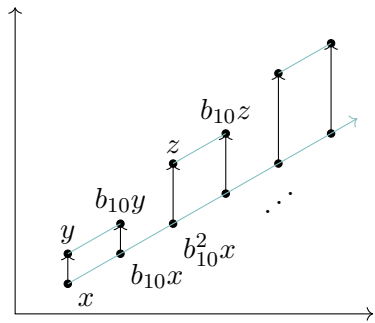
\begin{figure}[H]
\begin{center}
\begin{tikzpicture} [xscale=0.7,yscale=0.4]
\draw[<->] (0,10.2) -- (0,0) -- (6.8,0);
\draw[fill] (1,1) ellipse (2pt and 7/2pt);
\draw[fill] (1,2) ellipse (2pt and 7/2pt);
\draw[->] (1,1) -- (1,2);
\draw[fill] (2,2) ellipse (2pt and 7/2pt);
\draw[fill] (2,3) ellipse (2pt and 7/2pt);
\draw[->] (2,2) -- (2,3);
\draw[fill] (3,3) ellipse (2pt and 7/2pt);
\draw[fill] (3,5) ellipse (2pt and 7/2pt);
\draw[->] (3,3) -- (3,5);
\draw[fill] (4,4) ellipse (2pt and 7/2pt);
\draw[fill] (4,6) ellipse (2pt and 7/2pt);
\draw[->] (4,4) -- (4,6);
\draw[fill] (5,5) ellipse (2pt and 7/2pt);
\draw[fill] (5,8) ellipse (2pt and 7/2pt);
\draw[->] (5,5) -- (5,8);
\draw[fill] (6,6) ellipse (2pt and 7/2pt);
\draw[fill] (6,9) ellipse (2pt and 7/2pt);
\draw[->] (6,6) -- (6,9);
\draw[->,teal!50] (1,1) -- (6.5,6.5);
\draw[teal!50] (1,2) -- (2,3);
\draw[teal!50] (3,5) -- (4,6);
\draw[teal!50] (5,8) -- (6,9);
\node[below right] at (1,1) {$x$};
\node[below right] at (2,2) {$b_{10}x$};
\node[below right] at (3,3) {$b_{10}^2x$};
\node[right] at (4,4) {\hspace{10pt}$\iddots$};
\node[above] at (1,2) {$y$};
\node[above] at (2,3) {$b_{10}y$};
\node[above] at (3,5) {$z$};
\node[above] at (4,6) {$b_{10}z$};
\end{tikzpicture}
\end{center}
\caption{Convergence problem (2) for the ISS}
\end{figure}

In both of these cases, it is clear from the pictures that these cannot happen
if there is a vanishing line of slope equal to the degree of $b_{10}$, as
guaranteed by Lemma \ref{ISS-vanishing-line}.
\end{proof}
\begin{remark}
The same proof shows that the ISS for $b_{10}^{-1}\Ext_P(k,k)$ converges; in
particular, the vanishing line in Lemma \ref{ISS-vanishing-line} goes through
even with more $h_{ij}$'s and $b_{ij}$'s in the $E_1$ page.
\end{remark}

\section*{Appendix B: MPASS charts}
\renewcommand{\thesection}{B}
\setcounter{theorem}{0}
\setcounter{equation}{0}
%% sidewaysfigure works by hijacking the figure environment. But I want to use
%% the float package, which puts the other figures in less stupid places, and
%% *that* redefinition of the figure environment confuses sidewaysfigure. So this
%% is copied out of amscls.sty .
%\makeatletter
%\renewenvironment{figure}{%
%  \@float{figure}%
%}{%
%  \end@float
%}
%\makeatother
\rotatebox{270}{\begin{minipage}{0.75\textheight}
    \includegraphics[width=\textwidth]{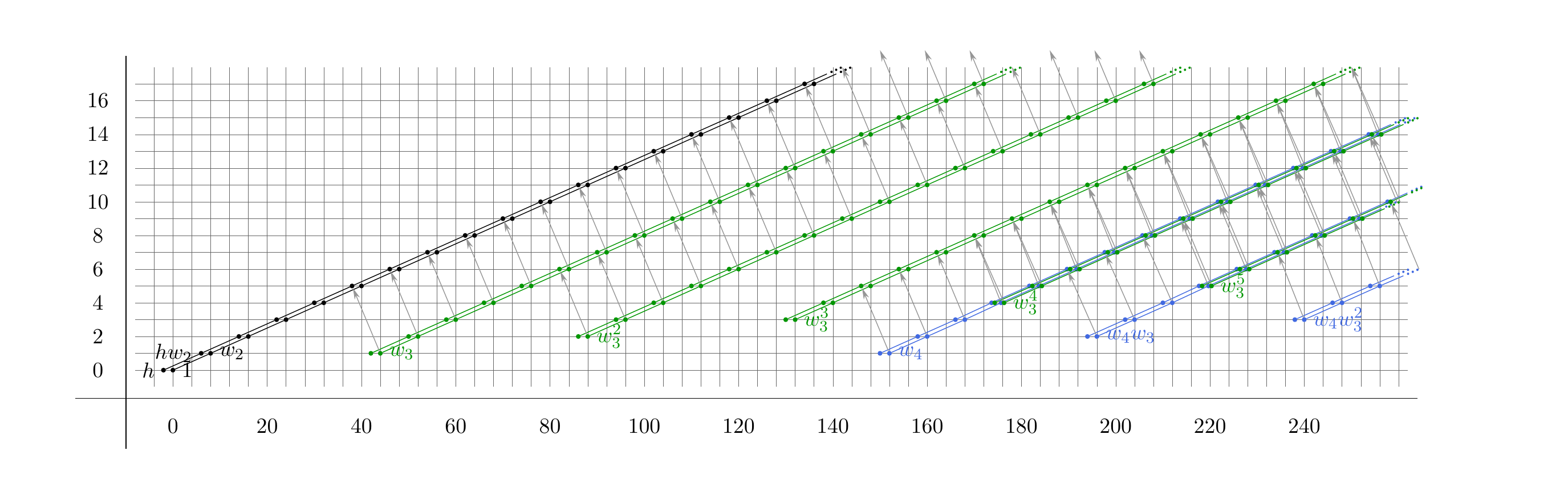}
    \captionof{figure}{$E_4$ page of the $K(\xi_1)$-based MPASS, with $d_4$ differentials shown}
    \label{figure:MPASS-E_4}
    \includegraphics[width=\textwidth]{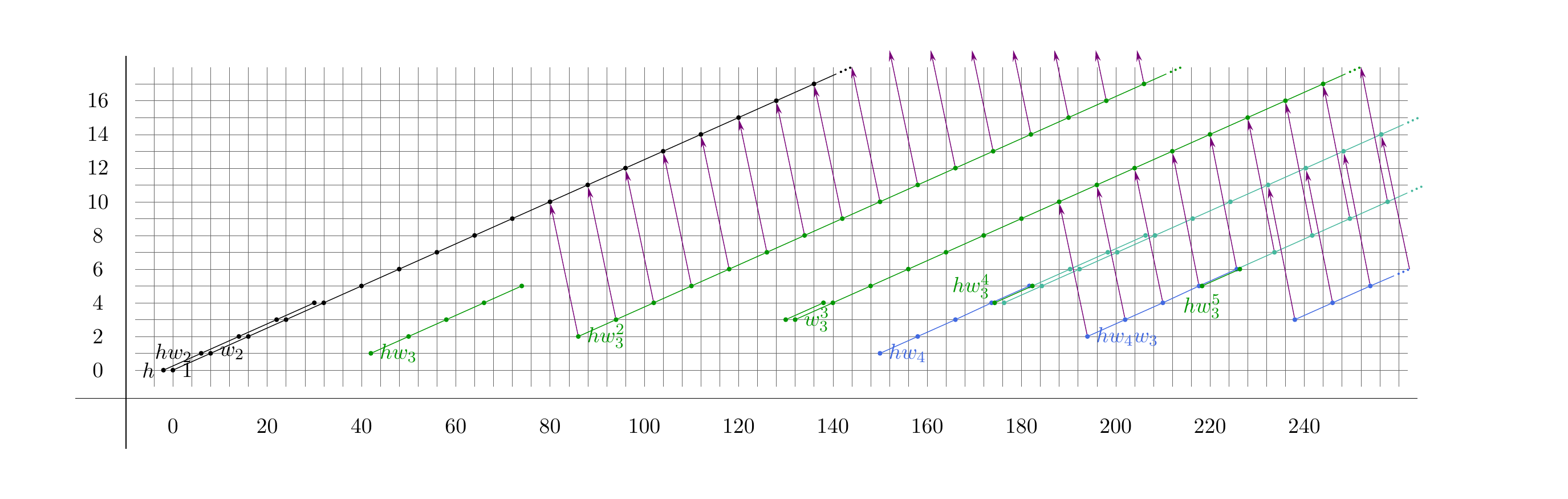}
    \captionof{figure}{$E_8$ page of the $K(\xi_1)$-based MPASS}
    \label{figure:MPASS-E_8}
\end{minipage}}
%\begin{sidewaysfigure}
%\vspace{425pt}
%    \includegraphics[width=\textheight]{../ust-grading1.pdf}
%    \caption{$E_4$ page of the $K(\xi_1)$-based MPASS, with $d_4$ differentials
%	shown}
%  \label{figure:MPASS-E_4}
%    \includegraphics[width=\textheight]{../ust-grading2.pdf}
%    \caption{$E_8$ page of the $K(\xi_1)$-based MPASS}
%  \label{figure:MPASS-E_8}
%\end{sidewaysfigure}

\rotatebox{270}{\begin{minipage}{0.75\textheight}
\vspace{-400pt}

    \includegraphics[width=\textwidth]{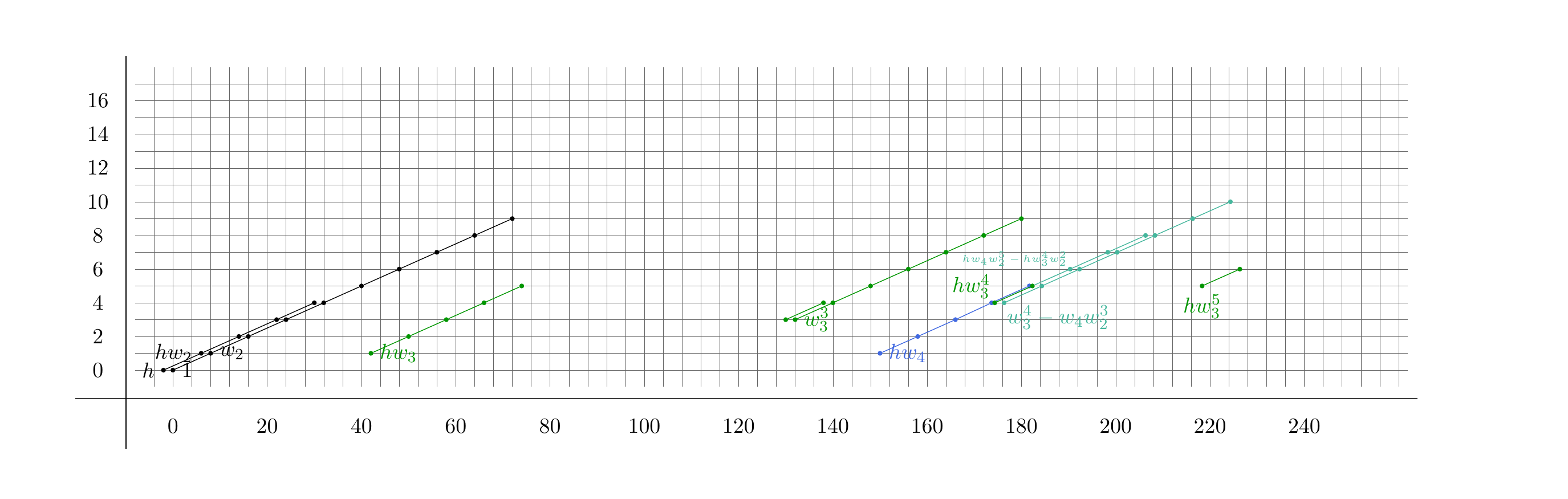}
    \captionof{figure}{$E_\iy$ page of the $K(\xi_1)$-based MPASS}
    \label{figure:MPASS-E_iy}
\end{minipage}}

%\begin{sidewaysfigure}
%\vspace{400pt}
%    \includegraphics[width=\textheight]{../ust-grading3.pdf}
%    \caption{$E_\iy$ page of the $K(\xi_1)$-based MPASS}
%  \label{figure:MPASS-E_iy}
%\end{sidewaysfigure}

\newpage
\bibliography{b10.bib}
\bibliographystyle{alpha}

\end{document}